\documentclass[12pt]{amsart}
\usepackage{amsbsy}
\usepackage{graphicx,epsfig,subfigure,psfrag}
\begin{document}

%%%%%%%%%%%%%%%%%%%%%%%%%%%%%%%%%%%%%%%%%%%%%%%%%%%%%%%%%%%%%%%%%%%%
% Theorem, definition, lemma, proposition, corollary and proof
%%%%%%%%%%%%%%%%%%%%%%%%%%%%%%%%%%%%%%%%%%%%%%%%%%%%%%%%%%%%%%%%%%%%
%%%%%%%%%%%%%%%%%%%%%%%%%%%%%%%%%%%%%%%%%%%%%%%%
\newtheorem{theorem}{Theorem}%[section]
\newtheorem{proposition}{Proposition}%[section]
\newtheorem{lemma}{Lemma}%[section]
\newtheorem{corollary}{Corollary}%[section]
\newtheorem{definition}{Definition}%[section]
\newtheorem{remark}{Remark}%[section]
%%%%%%%%%%%%%%%%%%%%%%%%%%%%%%%%%%%
%%%%%%%%%%%%%%%%%%%%%%%%%%%%%%%%%%%%%%%%%%%%%%                  NEW
%%\newcommand{\be}{\begin{equation}}
%%\newcommand{\ee}{\end{equation}}
%%%%%%%%%%%%%%%%%%%%%%%%%%%%%%%%%%%%%%%%%%%%%%%
%%%%%%%%%%%%%%%%%%%%%%%%%%%%%%%%%%%%%%%%%%%%
\newcommand{\tex}{\textstyle}
%%%%%%%%%%%%%%%%%%%%%%%%%%%%%%%%%%%%%%%%%%%%
%%%%%%%%%%%%%%%%%%%%%%%%%%%%%%%%%%%%%%%%%%%%
\numberwithin{equation}{section} \numberwithin{theorem}{section}
\numberwithin{proposition}{section} \numberwithin{lemma}{section}
\numberwithin{corollary}{section}
\numberwithin{definition}{section} \numberwithin{remark}{section}
%%%%%%%%%%%%%%%%%%%%%%%%%%%%%%%%%%%%%%%%%%%%
%%%%%%%%%%%%%%%%%%%%%%%%%%%%%%%%%%%%%%%%%%%%
\newcommand{\ren}{\mathbb{R}^N}
\newcommand{\re}{\mathbb{R}}
\newcommand{\n}{\nabla}
\newcommand{\iy}{\infty}
\newcommand{\pa}{\partial}
\newcommand{\fp}{\noindent}
\newcommand{\ms}{\medskip\vskip-.1cm}
\newcommand{\mpb}{\medskip}
%%%%%%%%%%%%%%%%%%%%%%%%%%%%%%%%%%%%%%%%%%%%%%%%%
\newcommand{\AAA}{{\bf A}}
\newcommand{\BB}{{\bf B}}
\newcommand{\CC}{{\bf C}}
\newcommand{\DD}{{\bf D}}
\newcommand{\EE}{{\bf E}}
\newcommand{\FF}{{\bf F}}
\newcommand{\GG}{{\bf G}}
\newcommand{\oo}{{\mathbf \omega}}
\newcommand{\Am}{{\bf A}_{2m}}
\newcommand{\CCC}{{\mathbf  C}}
\newcommand{\II}{{\mathrm{Im}}\,}
\newcommand{\RR}{{\mathrm{Re}}\,}
\newcommand{\eee}{{\mathrm  e}}
%%%%%%%%%%%%%%%%%%%%%%%%%%%%%%%%%%%%%%%%%%%%%%%%%%%%%%%%%%%%%%%%%%%%%%% L^2\rho...
\newcommand{\LL}{L^2_\rho(\ren)}
\newcommand{\LLL}{L^2_{\rho^*}(\ren)}
%%%%%%%%%%%%%%%%%%%%%%%%%%%%%%%%%%
%%%%%%%%%%%%%%%%%%%%%%%%%%%%%%%%%%%%%%%%%%%%%%%%%%%%
\renewcommand{\a}{\alpha}
\renewcommand{\b}{\beta}
\newcommand{\g}{\gamma}
\newcommand{\G}{\Gamma}
\renewcommand{\d}{\delta}
\newcommand{\D}{\Delta}
\newcommand{\e}{\varepsilon}
\newcommand{\var}{\varphi}
\renewcommand{\l}{\lambda}
\renewcommand{\o}{\omega}
\renewcommand{\O}{\Omega}
\newcommand{\s}{\sigma}
\renewcommand{\t}{\tau}
\renewcommand{\th}{\theta}
\newcommand{\z}{\zeta}
\newcommand{\wx}{\widetilde x}
\newcommand{\wt}{\widetilde t}
\newcommand{\noi}{\noindent}
 %%%%%%%%%%%%%%%%%%%%%%%%%%%%%%%%%%%%%%%%%%%
\newcommand{\uu}{{\bf u}}
\newcommand{\xx}{{\bf x}}
\newcommand{\yy}{{\bf y}}
\newcommand{\zz}{{\bf z}}
\newcommand{\aaa}{{\bf a}}
\newcommand{\cc}{{\bf c}}
\newcommand{\jj}{{\bf j}}
\newcommand{\ggg}{{\bf g}}
\newcommand{\UU}{{\bf U}}
\newcommand{\YY}{{\bf Y}}
\newcommand{\HH}{{\bf H}}
\newcommand{\GGG}{{\bf G}}
\newcommand{\VV}{{\bf V}}
\newcommand{\ww}{{\bf w}}
\newcommand{\vv}{{\bf v}}
\newcommand{\hh}{{\bf h}}
\newcommand{\di}{{\rm div}\,}
\newcommand{\ii}{{\rm i}\,}
%%%%%%%%%%%%%%%%%%%%%%%%%%%%%%%%%%
%%%%%%%%%%%%%%%%%%%%%%%%%%%%%%%%%%%%%   VAG, NEW
\newcommand{\inA}{\quad \mbox{in} \quad \ren \times \re_+}
\newcommand{\inB}{\quad \mbox{in} \quad}
\newcommand{\inC}{\quad \mbox{in} \quad \re \times \re_+}
\newcommand{\inD}{\quad \mbox{in} \quad \re}
\newcommand{\forA}{\quad \mbox{for} \quad}
\newcommand{\whereA}{,\quad \mbox{where} \quad}
\newcommand{\asA}{\quad \mbox{as} \quad}
\newcommand{\andA}{\quad \mbox{and} \quad}
\newcommand{\withA}{,\quad \mbox{with} \quad}
\newcommand{\orA}{,\quad \mbox{or} \quad}
\newcommand{\ef}{\eqref}
\newcommand{\ssk}{\smallskip}
\newcommand{\LongA}{\quad \Longrightarrow \quad}
%%%%%%%%%%%%%%%%%%%%%%%%%%%%%%%%
%%%%%%%%%%%%%%%%%%%%%%%%%%%%%%%%%%
\def\com#1{\fbox{\parbox{6in}{\texttt{#1}}}}
%%%%%%%%%%%%%%%%%%%%%%%%%%%%%%%%%%
%%%%%%%%%%%%%%%%%%% From Paper1
\def\N{{\mathbb N}}
\def\A{{\cal A}}
\newcommand{\de}{\,d}
\newcommand{\eps}{\varepsilon}
\newcommand{\be}{\begin{equation}}
\newcommand{\ee}{\end{equation}}
\newcommand{\spt}{{\mbox spt}}
\newcommand{\ind}{{\mbox ind}}
\newcommand{\supp}{{\mbox supp}}
\newcommand{\dip}{\displaystyle}
\newcommand{\prt}{\partial}
\renewcommand{\theequation}{\thesection.\arabic{equation}}
\renewcommand{\baselinestretch}{1.1}
%%%%%%%%%%%%%%%%%%%%%%%%%%%%%%%%%%%%%%%%%%%%%%%
\newcommand{\Dm}{(-\D)^m}

%%%%%%%%%%%%%%%%%%%%%%%%% VICTOR
\title
 {\bf Five types of blow-up in a semilinear\\
  fourth-order
 reaction-diffusion
 %%heat
 equation:\\
 an analytic-numerical approach}

\author {V.A.~Galaktionov}

\address{Department of Mathematical Sciences, University of Bath,
 Bath BA2 7AY, UK}
%%% and Keldysh Institute of Applied Mathematics,
%%%% Miusskaya Sq. 4, 125047 Moscow, RUSSIA}
\email{vag@maths.bath.ac.uk}

%%\address{Department of Mathematical Sciences, University of Bath,
%% Bath BA2 7AY, UK}
%%\email{ivk20@maths.bath.ac.uk}

%%\thanks{Research supported by  RTN network
%%HPRN-CT-2002-00274 and CERN-INTAS00-0136}

\keywords{4th-order semilinear parabolic equation,
 blow-up, self-similar solutions in $\ren$,
 non
self-similar blow-up.
%% $\sqrt{\ln|\ln(T-t)|}$-factor
%%% global existence, uniform bounds
%% discrete real spectrum, complete set of
%%%eigenfunctions, asymptotic behaviour.
 %% \\
 %% {\bf Submitted to:}
  }
 %%%%% J. Nonl. Sci.}

\date{\today}
%%% \quad {\bf Fila/4BlowN1.tex}}
%%%%%%%%%%%%%%%%%%%%%%%%%%%%%%%%%%%%%%%%%%%%%%%%%%%%%%%%

%%%%%%%%%%%%%%%%%%%%%% VICTOR
\begin{abstract}
%%%%%%%%%%%%%%%%%%%%%% VICTOR

Five types of {\em blow-up  patterns} that can occur for the
4th-order semilinear parabolic equation of reaction-diffusion type
 $$
  \tex{
u_t= -\D^2 u + |u|^{p-1} u \inB \ren \times (0,T), \,\, p>1, \quad
\lim_{t \to T^-}\sup_{x \in \ren} |u(x,t)|= +\iy,
 }
 $$
are discussed. For the semilinear heat equation $u_t= \D u+ u^p$,
various blow-up patterns were under scrutiny since 1980s, while
the case of higher-order diffusion was studied much less,
regardless a wide range of its application. The types of blow-up
include:

\noi (i) {\bf Type I(ss)}: various patterns of {\bf s}elf-{\bf
s}imilar single point  blow-up, including those, for which the
final time profile $|u(\cdot,T^-)|^{N(p-1)/4}$  is a measure;

 \noi (ii) {\bf Type I(log)}: self-similar non-radial blow-up
with angular {\bf log}arithmic TW swirl;

\noi (iii) {\bf Type I(Her)}: non self-similar blow-up close to
stable/centre subspaces of {\bf Her}mitian  operators obtained via
linearization about constant uniform blow-up pattern;

\noi (iv) {\bf Type II(sing)}: non self-similar blow-up on
stable/centre manifolds of a {\bf sing}ular steady state
%% and
%%matching
 in the supercritical Sobolev range $p \ge p_{\rm S}=
\frac{N+4}{N-4}$ for $N>4$; and

\noi (v) {\bf Type II(LN)}: non self-similar blow-up along the
manifold of stationary generalized {\bf L}oewner--{\bf N}irenberg
type explicit solutions in the critical Sobolev case $p=p_{\rm
S}$, when $|u(\cdot,T^-)|^{N(p-1)/4}$ contains a measure as a
singular component.

All proposed types of blow-up are very difficult to justify, so
formal analytic and numerical methods are key in  supporting some
theoretical judgements.
%%A
%%rigorous proof of some of the above blow-up scenario seems
%%illusive.

%%%%%%%%%%%%%%%%%%%%%% VICTOR
\end{abstract}
%%%%%%%%%%%%%%%%%%%%%%

%%%%%%%%%%%%%%%%%%%%%%%%%%%
\maketitle

%%%%%%%%%%%%%%%%%%%%%%%%%%%
\section{From second-order to higher-order blow-up R--D models: a PDE
 route  from XXth to XXIst century}
\label{S.1}

%%%%%%%%%%%%%%%%%%%%%%%%%%%%%%%%%%%%%%%%%%%%
\subsection{The RDE--4 and applications}

This paper is devoted to a description of blow-up patterns for the
%%fourth-order semilinear heat equation, or
{\em fourth-order reaction-diffusion equation} (the RDE--4 in
short)
\be
\label{m2}
  \tex{
u_t= -\D^2 u + |u|^{p-1} u \inB \ren \times (0,T) \whereA p>1,
 }
 \ee
 where $\D$ stands for the Laplacian in $\ren$.
 %% so that the equation contains the bi-harmonic
 %% diffusion operator .
 This has
  the
 bi-harmonic diffusion  $-\D^2$ and is a higher-order counterpart
 of classic second-order  PDEs, which we begin our discussion with.
 For
 applications of such higher-diffusion models, see short surveys and references in
 \cite{BGW1, GW1}. In general,
 higher-order semilinear parabolic equations arise in many physical
applications such as thin film theory, convection-explosion
theory, lubrication theory, flame and wave propagation (the
Kuramoto-Sivashinsky equation and the extended Fisher-Kolmogorov
equation), phase transition at critical Lifschitz points,
bi-stable systems and applications to structural mechanics; the
effect of fourth-order terms on self-focusing problems in
nonlinear optics are also well-known in applied and mathematical
literature. For a systematic treatment of extended KPPF-equations,
see Peletier--Troy  \cite{PelTroy}.

Note  that another related fourth-order one-dimensional semilinear
parabolic equation
\begin{equation}
\label{JS}
 u_t =
   - u_{xxxx} - [(2-(u_x)^2)u_x]_x - \a u + q {\mathrm e}^{su},
\end{equation}
where $\a $, $q$ and $s$ are positive constants obtained from
physical parameters,
 occurs in the Semenov-Rayleigh-Benard problem
\cite{JMS}, where the equation is derived  in studying the
interaction between natural convection and the explosion of an
exothermically-reacting fluid confined between two isothermal
horizontal plates. This is an evolution equation for the
temperature fluctuations in the presence  of natural convection,
wall losses and chemistry. It can be considered as a formal
combination of the equation derived in \cite{GS} (see also
\cite{CP}) for the Rayleigh-Benard problem and of the Semenov-like
energy balance \cite{Sem, F-K} showing that  natural convection
and the explosion mechanism may reinforce each other;
 see more details on physics
and mathematics of blow-up in \cite{GW1}.
 In a special limit, \ef{JS} reduces to the {\em generalized
 Frank-Kamenetskii equation} (see \cite{BGW1} for blow-up stuff)
 \begin{equation}
\label{u4e} u_t = -u_{xxxx} + {\mathrm e}^u,
\end{equation}
 which is a natural extension of the classic
 Frank-Kamenetskii equation; see below.

  Equation  \ef{m2} can be considered as a
 non-mass-conservative counterpart of the well-known {\em limit unstable
 Cahn--Hilliard equation} from phase transition,
 %%% and other theories
  \be
   \label{CH1}
   u_t= - u_{xxxx} - (|u|^{p-1}u)_{xx} \inB \re \times \re_+,
    \ee
 which is known to admit various families of blow-up solutions; see \cite{EGW1}
 for a long list of references. Somehow, \ef{m2} is related to the
 famous {\em Kuramoto--Sivashinsky equation} from flame
 propagation theory
  \be
  \label{KS1}
  u_t= - u_{xxxx} -u_{xx} + u u_x \inB \re \times \re_+,
   \ee
 which always admits  global solutions, so no blow-up for \ef{KS1} exists.

%%%%%%%%%%%%%%%%%%%%%%%%%%%%%%%%%%%%%%%%%
\subsection{On second-order reaction-diffusion (R--D) equations: a training ground
of blow-up PDE research in the XXth century}

 Blow-up phenomena, as examples of extremely
nonstationary behaviour of nonlinear mechanical and physical
systems, become more natural in PDE theory since  a systematic
developing combustion theory in the 1930s.
 This essential combustion influence began with the derivation
 of the semilinear parabolic
reaction-diffusion PDE such as
 the classic {\em
Frank-Kamenetskii equation} (1938) \cite{Fr-K}
 \be
  \label{FK1}
u_t= \D u + {\mathrm e}^u \inA,
 \ee
which occurs in combustion theory  of solid fuels and is often
also called the {\em solid fuel model}. First blow-up results in
related ODE models are due to Todes in 1933; see the famous
monograph \cite{ZBLM} for details of the history and applications.
The related model with a power superlinear source term takes the
form (also available among various nonlinear combustion models
\cite{ZBLM})
 \be
  \label{FK1p}
u_t= \D u + |u|^{p-1}u \inA \whereA p>1.
%%% \quad (u=u(x,t) \ge 0).
 \ee
 Thus, for such typical models, {\em blow-up} means that in the
 Cauchy problem\footnote{For simplicity, we avoid using initial-boundary
 value problems, where boundary conditions can affect some
 manipulations and speculations around; though can be included.},
the classic bounded solution $u=u(x,t)$ exists in $\ren \times
(0,T)$, while
 \be
 \label{Bl1}
 \tex{
  \sup_{x \in \ren} \, |u(x,t)| \to +\iy \asA t \to T^-,
  }
  \ee
  where $T \in \re_+=(0,+\iy)$ is then called the {\em blow-up time} of
  the solution $u(x,t)$.

During last fifty years of very intensive research starting from
seminal Fujita results in 1966 (on what is now called {\em Fujita
exponents}), we have currently got rather complete understanding
of the types of blow-up for the semilinear \ef{FK1}, \ef{FK1p} and
other models. This is very well explained in a number of
monographs; see
 \cite{BebEb, SGKM, GSVR, Pao, MitPoh, AMGV, GalGeom, QSupl}.

However, one should remember that even for simple R--D equation
such as \ef{FK1} and \ef{FK1p}, there are blow-up scenarios in the
multi-dimensional geometries, which still did not get a proper
rigorous mathematical justification. For instance, there are a
number surprises even in the radial geometry for \ef{FK1p}, which
 reads for $r=|x|>0$ as
 \be
 \label{p2}
  \tex{
  u_t= \frac 1{r^{N-1}}\, \big(r^{N-1} u_r\big)_r+ |u|^{p-1} u
  \inB \re_+ \times \re_+ \quad (u_r|_{r=0}=0, \,\,\, \mbox{
  symmetry}),
   }
   \ee
   in the supercritical Sobolev range
    \be
    \label{S1}
    \tex{
    p > p_{\rm S}= \frac{N+2}{N-2} \whereA N>2.
     }
     \ee
Several critical exponents, which may  essentially change blow-up
evolution, appear for \ef{p2} in the range \ef{S1}, among those
let us mention the most amazing ones:
 \be
 \label{S2}
  \tex{
  p_{\rm JL}= 1+ \frac 4{N-4-2\sqrt{N-1}}, \, p_{\rm L}= 1+ \frac 6{N-10} \,\, (N \ge 11);
  \,\, p_{\rm M}=1+ \frac 7{N-11} \,\,(N \ge 12);
  \,\,\, \mbox{etc.}
  }
  \ee
In particular, this shows that, in the parameter range
 \be
 \label{S3}
 N \ge 11 \andA p \ge p_{\rm S},
  \ee
 new principal issues of blow-up evolution for \ef{p2} essentially  take place.
 Note that, in \cite{GV}, some critical blow-up exponents were shown to
 exist for the quasilinear combustion equation with a porous
 medium diffusion:
  \be
  \label{mm1}
  u_t= \D u^{m} + u^p \whereA  p>m>1 \quad (u(x,t) \ge 0).
   \ee
   This shows certain universality of formation of blow-up
   singularities for a wider class of R--D equations, which now we
   are going to extend to the RDE--4 \ef{m2}.

 We do not plan to give any detailed enough  review of such a variety of these delicate and becoming
 diverse   (rather surprisingly)
   in the XXIst century
  mathematical
 results, which quite  recently attracted the attention of several remarkable
 mathematicians
  from various areas of PDE theory. We refer to \cite{Vel, HVsup, GV} for earlier results since 1980s and 90s,
    and to more recent papers
  \cite{Fila05, Mat04} and
 \cite{Miz04}--\cite{Miz07} as a guide to the research, which
 was essentially intensified last few years. Further results can
 be traced out by the {\tt MathSciNet}, using most recent papers
 of the authors mentioned above.

It is worth mentioning that most of these results
%%mentioned above
have been obtained for nonnegative blow-up solutions of \ef{FK1},
\ef{FK1p}, and \ef{S1}, since the positivity property is naturally
supported by the Maximum Principle (the MP) for such second-order
parabolic equations.
 For instance, a full classification of such
nonnegative blow-up patterns for \ef{FK1p} (all of them belong to
the family Type I(Her)) in the subcritical range $1<p<p_{\rm S}$
was obtained in \cite{Vel}. For \ef{FK1}, this happens in
dimension $N=1$ and 2. In other words, the family of blow-up
patterns for
 \ef{FK1} and subcritical \ef{FK1p} first formally introduced in \cite{GHPV}
  is {\em evolutionary
 complete} (a notion from \cite{CompG}, where further
 references can be found).
 In the range $p \ge p_{\rm S}$ for \ef{FK1p} and from $N=3$ for \ef{FK1}, there occur self-similar
patterns of Type I(ss) and many others being non-self-similar,
which makes the global blow-up flow much more complicated.

%% more than twenty five years ago.
For $p \ge p_{\rm S}$, such a complete classification for
\ef{FK1p} is far from being complete. E.g.,
 \be
 \label{S4}
 \mbox{for $p \ge p_{\rm S}$ in (\ref{FK1p}), nonsymmetric blow-up patterns are
practically unknown}.
 \ee
 Moreover,
 for solutions of changing sign, the results
are much more rare and are essentially incomplete. It is worth
mentioning surprising blow-up patterns of changing sign
constructed in \cite{Fil00}, with the structure to be used later
on for \ef{m2}, where we comment on this Type II(LN) blow-up
patterns for \ef{m2} in greater detail.

%%%%%%%%%%%%%%%%%%%%%%%%%%%%%%%%%%%%%%%%%%%%%%%%
\subsection{Back to the RDE--4: five types of blow-up patterns and layout of the paper}

We are going to discuss possible types of blow-up behaviour for
the RDE--4 \ef{m2}. In what follows, we are using the auxiliary
classification from Hamilton
 \cite{Ham95}, where Type I blow-up
 means the solutions satisfying, for some constant $C >0$ (depending
 on $u$),
  \be
  \label{S5}
  \begin{matrix}
 \mbox{Type I:} \quad
 (T-t)^{\frac 1{p-1}} |u(x,t)| \le C  \asA t \to T^-,
  \,\,\,\mbox{and, otherwise}\ssk\ssk\\
  \mbox{Type II:} \quad \limsup_{t \to T^-} (T-t)^{\frac 1{p-1}}
  \sup_x |u(x,t)|=+\iy \qquad\qquad\,\,\,
  %% \quad(\mbox{Type I}),
  \end{matrix}
   \ee
   (Type II also called {\em slow} blow-up in
   \cite{Ham95}).
   In R--D theory, blow-up with the dimensional estimate \ef{S5}
   was usually called of {\em self-similar} rate, while Type II was
   referred to as {\em fast} and {\rm non self-similar}; see
   \cite{AMGV} and \cite{SGKM}.

   \ssk

Thus, we plan to describe the following five types of blow-up with
an extra classification issues in each of them (this list also
shows the overall layout of the paper):

\ssk

 \noi (i) {\bf Type
I(ss)}: various patterns of {\bf s}elf-{\bf s}imilar single point
blow-up mainly in radial geometry, including those, for which
$|u(\cdot,T^-)|^{N(p-1)/4}$  is a measure (Section \ref{S.2};
almost nothing is known for non-radial similarity blow-up patterns
for $N \ge 2$, which are expected to exist);

\ssk

 \noi (ii) {\bf Type I(log)}: non radial self-similar blow-up with
 angular
{\bf log}arithmic travelling wave ({\bf log}TW) swirl, which in
the similarity rescaled variables corresponds to periodic orbits
as $\o$-limit sets (Section \ref{S.3});

\ssk

 \noi (iii) {\bf Type I(Her)}: non self-similar blow-up close to
stable/centre subspaces of {\bf Her}mitian  operators obtained via
linearization about constant uniform blow-up and matching with a
Hamilton--Jacobi region (Section \ref{S.4});

\ssk

 \noi (iv) {\bf Type II(sing)}: non self-similar blow-up on
stable/centre manifolds of {\bf sing}ular steady state (SSS) in
the supercritical Sobolev (and Hardy) range $p > p_{\rm S}=
\frac{N+4}{N-4}$ for $N>4$ with matching to a central
quasi-stationary region (Section \ref{S.5}); and

\ssk

 \noi (v) {\bf Type II(LN)}: non self-similar blow-up along the
manifold of stationary generalized
 {\bf L}oewner--{\bf N}irenberg type
explicit solutions in the critical Sobolev case $p=p_{\rm S}$,
when final time profiles  $|u(\cdot,T^-)|^{N(p-1)/4}$ contain
measures in the singular component (Section \ref{S.6}).

\ssk

We must admit that the analysis of all the blow-up type indicated
above is very difficult mathematically, so we do not present
practically no rigorous results. Recall that, even for the
second-order equation \ef{FK1}, all these types excluding Type
I(Her)  still did not have not only any complete classification,
but some of them were not detected at all.
%% It is worth mentioning,
For \ef{m2},  the best known  critical exponent is obviously
Sobolev's one
 \be
 \label{S6}
 \tex{
 p_{\rm S}= \frac{N+4}{N-4} \whereA N\ge 5, \,\,\, \mbox{and also} \,\,\, p_{\rm *}= \frac
 N{N-4},
 }
 \ee
 while the others, as counterparts of those in \ef{S2},
 need further study and understanding.
 %%as ue what is the meaning of other
 %%to be introduced need further
 %%and currently does not have a cl etc.
  However,
  %%% most plausibly,
 many  critical exponents for \ef{m2}  cannot
be explicitly calculated. Overall, we aim that our approaches to
blow-up patterns can be extended to $2m$th-order parabolic
equations such as
 \be
 \label{mm98}
  u_t= -(-\D)^m u + |u|^{p-1}u,
   \ee
 though the case $m=2$ (the first even $m$'s) already contains some surprises.

Nevertheless, it seems that, at some stage  of struggling for
developing new concepts, it is inevitable to attempt to perform
  a formal classification under the clear danger of a lack of
any rigorous justification\footnote{Actually following
Kolmogorov's legacy from the 1980s sounding not completely
literally as:
%%%%%s\footnote{%%This well corresponds to:
 %%\be
%%  \label{Kol1}
%% \fbox{$
%%  \begin{matrix}
%% \mbox{
``The main goal of a mathematician is not
 proving a theorem,
 %%} \ssk\\
%% \mbox{
 but an effective investigation of the
 problem..."\,.}.
%%  \end{matrix}
%% $}
%%% \ee
 %% (the author apologizes for a
%%%non-literal translation from the Russian).
 In this rather paradoxical connection, it is also worth mentioning that
%%It is well-understood  (see e.g., \cite[Ch.~5]{Maj02}) that
  the most well-known nowadays and the fundamental open problem of fluid
 mechanics\footnote{The Millennium Prize Problem for the Clay Institute; see Fefferman
 \cite{Feff00}.}
  and  PDE theory on {\em global existence or nonexistence (blow-up) of
 bounded smooth $L^2$-solutions} of the
 {\em Navier--Stokes equations} (the NSEs)
 %%\footnote{Claude Louis Marie Henri Navier, 1785-1836, and George
%%Gabriel Stokes, 1819-1903.} {\em equations} (the NSEs),
  \be
  \label{NS1}
 \uu_t +(\uu \cdot \n)\uu=- \n p + \D \uu, \quad \di \uu=0
\inB \re^3 \times
 \re_+,
 \quad \uu|_{t=0}=\uu_0\in L^2 \cap L^\iy,
  \ee
%% with {\em arbitrary bounded divergence-free $L^2$-data} $\uu_0$, is
%% two-fold:
 from one side belongs to a ``blow-up configurational" type: {\em to
 predict possible swirling
``twistor-tornado" type of blow-up patterns}. Moreover, it seems
that the NSEs \ef{NS1} was the first model, for which J.~Leray  in
1934 \cite[p.~245]{Ler34} formulated the so-called {\em Leray's
scenario} of self-similar blow-up as $t \to T^-$ and a similarity
continuation beyond for $t>T$.
 %% J.~Leray in 1934 (the so-called {\em Leray's scenario of
%%blow-up})
%%  in studying blow-up
 %%precisely corresponds
 %%to  (1934) of blow-up
%%  in the Navier--Stokes
%% equations in $\re^3$; see \cite[p.~245]{Ler34} for the formulation and
%% \cite{GalJMP} for recent discussions.
Nonexistence of such similarity blow-up for the NSEs \ef{NS1} was
proved in 1996 in Ne\v{c}as--Ru\v{z}i\v{c}ka--\v{S}ver\'ak
\cite{Nec96}. However, for the semilinear heat equations \ef{FK1p}
and \ef{mm1}, the validity of Leray's scenario of blow-up was
rigorously established; see \cite{GV, Miz06} and references
therein.

In general, we observe certain similarities between these two
blow-up
 problems; see \cite{GalJMP}, where Type I(log) patterns were introduced for \ef{NS1} and
 \cite{GalADE} for
   more  details and references on other related exact blow-up
solutions. Overall, we claim that equations \ef{FK1}, \ef{FK1p},
\ef{m2}, and \ef{NS1} admit some similar principles of
constructing various families of blow-up patterns, though, of
course, for the last two ones,  the construction gets essentially
harder and many steps are made formally, without proper
justification.  Especially for the NSEs \ef{NS1}, which compose a
nonlocal solenoidal parabolic equation:
  \be
  \label{NS1S}
 \uu_t +{\mathbb P}(\uu \cdot \n)\uu= \D \uu,
 %% \quad \di \uu=0
%%\inB \re^3 \times
 %%\re_+,
 %%\quad \uu|_{t=0} \in L^2 \cap L^\iy,
  \ee
  where the integral operator ${\mathbb P}=I- \n \D^{-1} (\n \cdot)$ is
  Leray--Hopf's projector onto the solenoidal vector
    field.
 More precisely, the RDE--2 \ef{FK1p} obeying the MP is indeed too
simple to mimic any NSEs blow-up patterns, while \ef{m2}, which
similar to \ef{NS1S} traces no MPs, can be about right (possibly,
still   illusionary). Then \ef{m2} stands for an auxiliary
``training ground" to approach understanding of mysterious and
hypothetical blow-up  for \ef{NS1}.
 %%% though some will be derived.

In Appendix A, we present other families of PDEs, which expose a
similar open problem on existence/nonexistence of $L^\iy$-blow-up
of solutions from  bounded smooth initial data. Overall, it is
worth saying that  the problem of description of blow-up patterns
and their {\em evolution completeness} takes and shapes certain
universality features in general PDE theory of the twenty first
century.

%%%%%%%%%%%%%%%%%%%%%%%%%%%%%%%%%%%%%%%%%%%%%%%%%%%%%
\section{{\bf Type
I(ss)}:
%%radially symmetric
 self-similar blow-up}
 \label{S.2}

This is the simplest and most natural type of blow-up for scaling
invariant equations such as \ef{m2}, where the behaviour as $t \to
T^-$ is given by a self-similar solution:
 \be
 \label{2.1}
 \tex{
 u_{\rm S}(x,t)=(T-t)^{-\frac 1{p-1}} f(y), \quad y= \frac
 x{(T-t)^{1/4}},%% \quad \mbox{and $f$ solves}
 }
  \ee
  where a non-constant function $f \ne 0$ is a proper solution of the elliptic
  problem:
  \be
  \label{2.2}
   \tex{
  {\bf A}(f) \equiv -\D^2 f- \frac 14 \, y \cdot \n f - \frac 1{p-1}\,
  f + |f|^{p-1}f=0 \inB \ren, \quad f(\iy)=0.
  }
  \ee
 We recall that, for \ef{FK1p}, such nontrivial self-similar Type
 I blow-up is nonexistent in the subcritical range $p \le
 \frac{N+2}{N-2}$. But this is not the case for the RDE--4
 \ef{m2}. Note that \ef{2.2} is a very difficult elliptic
 equation with the non-coercive and non-monotone operators, which are
  not variational in any weighted $L^2$-spaces. There are no still
  any sufficiently general results of solvability of \ef{2.2} in
  higher dimensions, so our research is a first
  attempt.
  %%% towards this.

In what follows, for any dimension $N \ge 1$, by $f_0(y)$ we will
denote the first monotone radially symmetric blow-up profile,
which, being on the lower $N$-branch (so $f_0$ is not unique, see
explanations below) is expected to be generic (i.e., structurally
stable in the rescaled sense). We also deal with the second
symmetric profile $f_1(y)$, which seems to be unstable, or, at
least, less stable than $f_0$. There are also other similarity
solutions concentrated about the singular SSS $U(y)$ (see Section
\ref{S.5}), but those, being adjacent to the unstable equilibrium
$U$ are expected to be unstable also.

It is worth mentioning that self-similar blow-up for \ef{m2} is
{\em incomplete}, i.e., blow-up solutions, in general,  admit
global extensions for $t>T$. Such principal questions
 are studied in \cite{GalBlExt} and will not be treated here.

%%%%%%%%%%%%%%%%%%%%%%%%%%%%%%%%%%%%%%%%
\subsection{One dimension: first examples of nonuniqueness}

\begin{figure}
%%\vskip -.3cm
\centering \subfigure[$p=\frac 32$ and $N=1$]{
\includegraphics[scale=0.52]{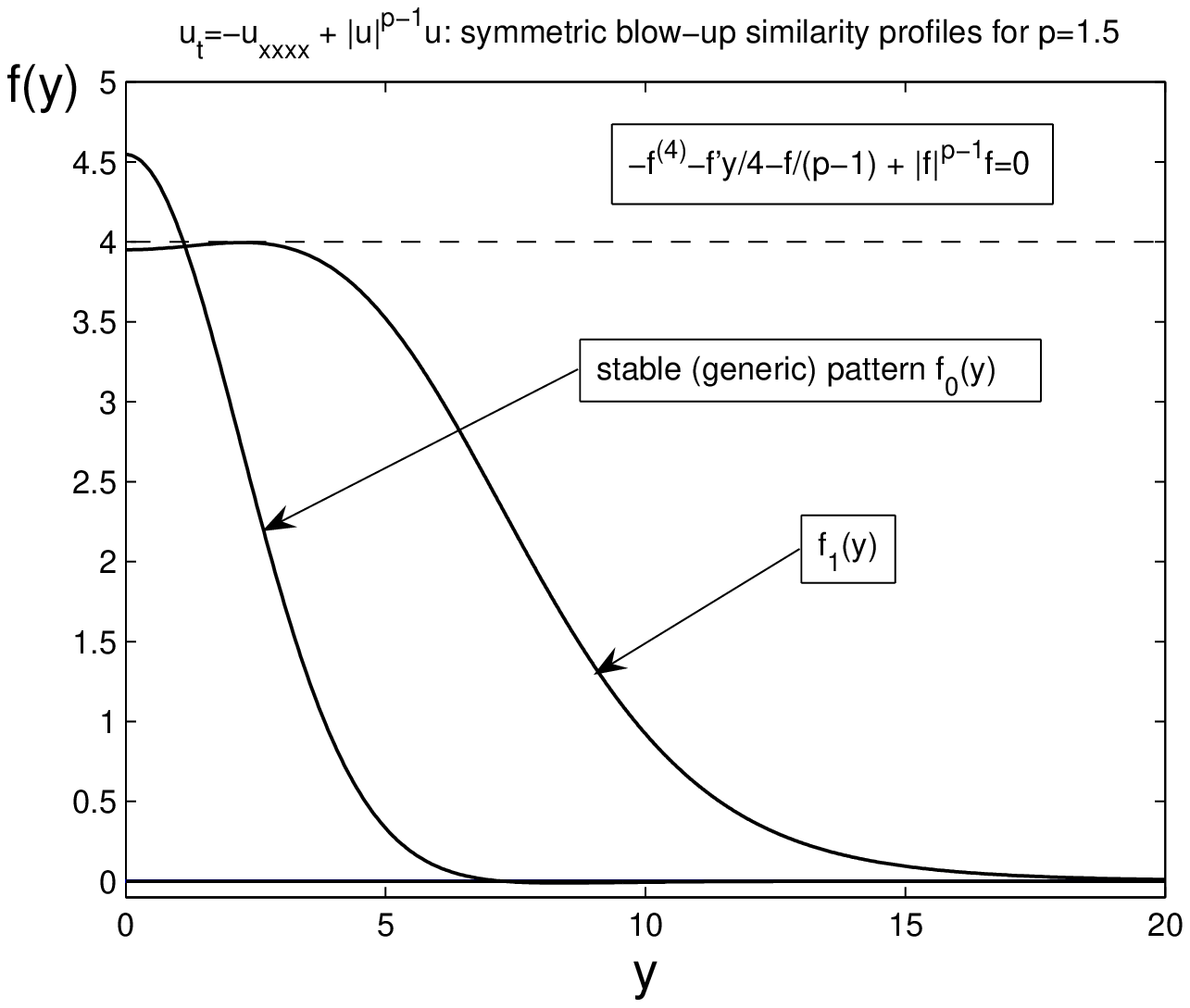}
} \subfigure[$p=2$ and $N=1$]{
\includegraphics[scale=0.52]{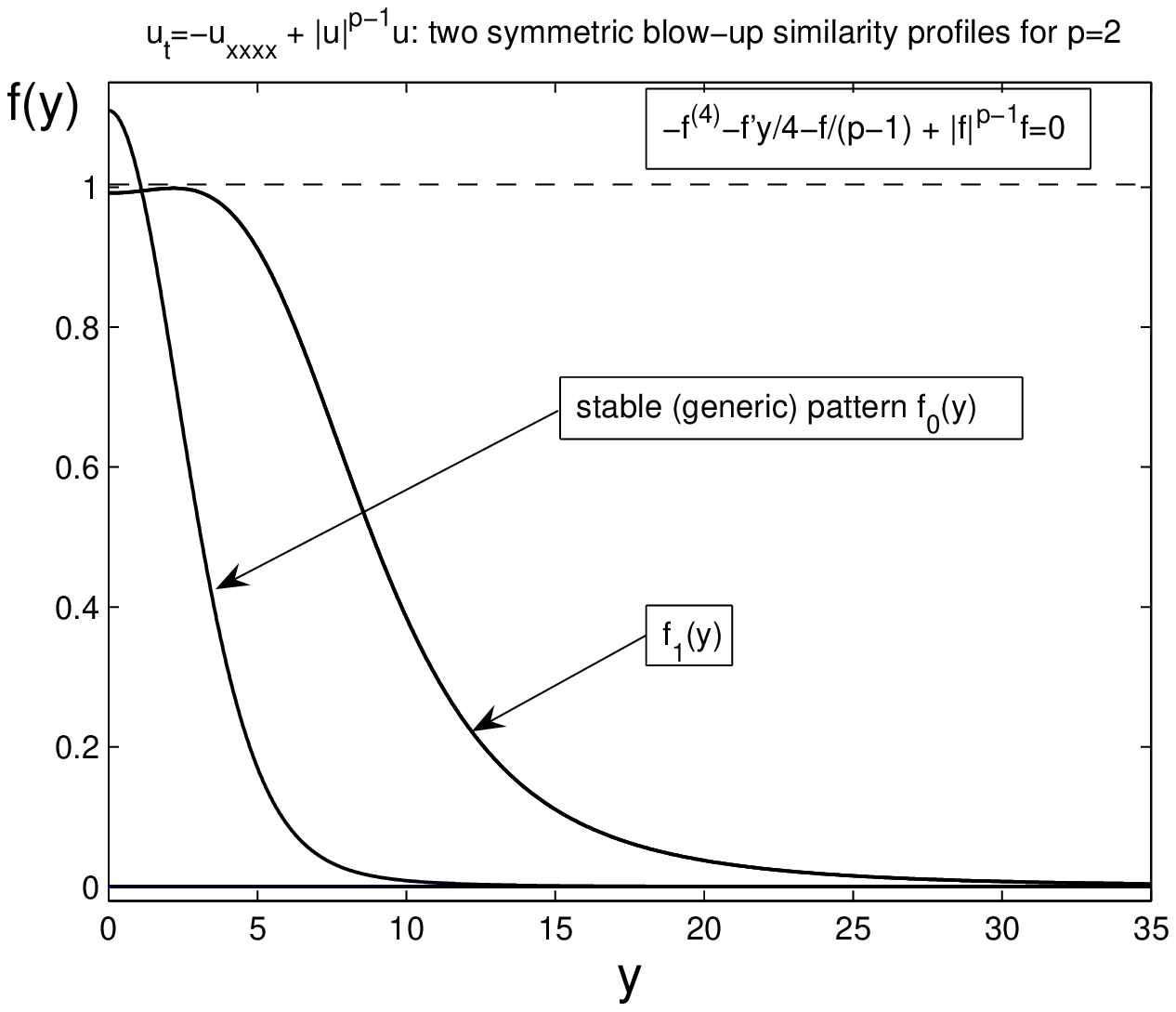}
}
 \vskip -.2cm
\caption{\rm\small Two self-similar blow-up solutions of \ef{2.3}:
$p=1.5$ (a) and $p=2$ (b).}
 %% \vskip -.3cm
 \label{F1}
\end{figure}

\begin{figure}
%  \vskip -.3cm
\centering
\includegraphics[scale=0.75]{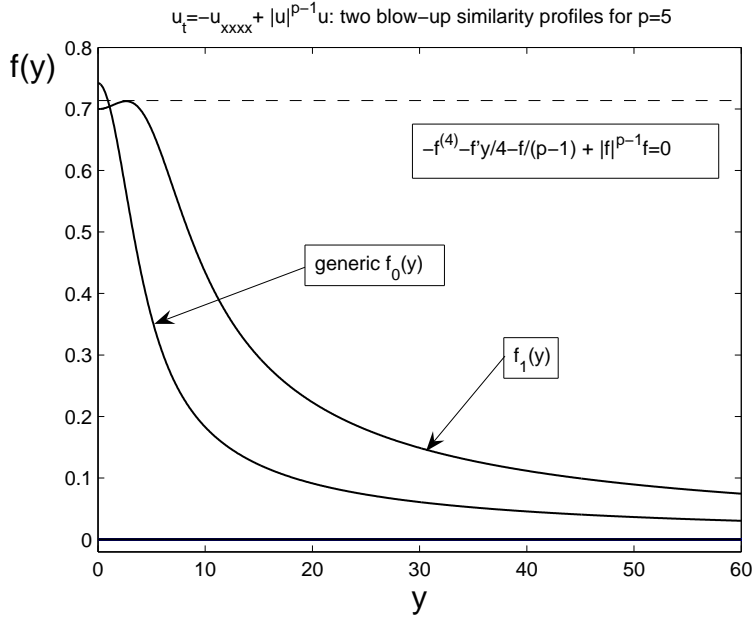} %%%%{AmFBP1.eps}  %%%%%%%{4het5.eps}   %%%%%{4het.eps}  Old
\vskip -.3cm \caption{\small  Two self-similar blow-up solutions
of \ef{2.3} for $p=5$.}
%%   \vskip -.3cm
 \label{F11}
\end{figure}

 Thus, for $N=1$, \ef{2.2} becomes the ODE
\be
  \label{2.3}
   \tex{
  {\bf A}(f) \equiv -f^{(4)}- \frac 14 \, y f' - \frac 1{p-1}\,
  f + |f|^{p-1}f=0 \inB \re, \quad f(\pm \iy)=0,
  }
  \ee
which was studied in \cite{BGW1} by a number of analytic-branching
and numerical methods. It was shown that \ef{2.3} admits at least
two different blow-up profiles with an algebraic decay at
infinity. See \cite[\S~3]{GW1} for further centre manifold-type
 arguments supporting this multiplicity result in a similar
4th-order blow-up problem.
 Without going into detail of such a study, we present a
few illustrations only and will address the essential dependence
of similarity profiles $f(y)$ on $p$. In Figure \ref{F1}, we
present those pairs of solutions of \ef{2.3} for $p=\frac 32$ and
$p=2$. All the profiles are symmetric (even), so satisfy the
symmetry condition
 \be
 \label{2.31}
 f'(0)=f'''(0)=0.
 \ee
No non-symmetric blow-up was detected in numerical experiments
(though there is no proof that such ones are nonexistent: recall
that ``moving plane" and Aleksandrov's Reflection Principle
methods do not apply to \ef{m2} without the MP). Figure \ref{F11}
shows similar two blow-up profiles for $p=5$.

%%%%%%%%%%%%%%%%%%%%%%%%%%%%%%%%%%%%%%%%%%%%%%%%%%%%%%%%%%%%%%%
%%FIG%%%%%%%%%%%%%%%%%%%%%%%

%%%%%%%%%%%%%%%%%%%%%%%%%%%%%%%%%%%%%%%%%%%%%%%%%%%%%%%%%%%%%%%%%%%

%%%%%%%%%%%%%%%%%%%%%%%%%%%%%%%%%%%%%%%%%%%%%%%%%%%%%%%%%%%%

%%%%%%%%%%%%%%%%%%%%%%%%%%%%%%%%%%%%%%%%%%%%%%%%%%%%%
\subsection{On existence of similarity profiles for $N=1$:
classification of blow-up and oscillatory
bundles}

We now provide extra details concerning  existence of at least a
single blow-up profile $f(y)$ satisfying \ef{2.3}, \ef{2.31}. We
perform shooting from $y=+\iy$
 by using the 2D bundle \ef{dd1} to $y=0$, where the symmetry condition \ef{2.31} are posed
 (or to $y=-\iy$, where the same bundle \ef{dd1} with $y \mapsto
 -y$ takes place). By $f=f(y;C_1,C_2)$, we denote the
 corresponding solution defined on some maximal interval
  \be
  \label{ma1}
   \tex{
y \in (y_0, +\iy) \whereA y_0=y_0(C_1,C_2)  \ge -\iy. }
 \ee
If $y_0(C_1,C_2)=-\iy$, then the corresponding solution
$f(y;C_1,C_2)$ is global and can represents a proper blow-up
profile (but not often, see below). Otherwise:
 \be
 \label{ma2}
 y_0(C_1,C_2) >-\iy \LongA f(y;C_1,C_2) \to \iy \asA y \to y_0^+.
  \ee
Note that ``oscillatory blow-up" for the ODE close to $y=y_0^+$:
 \be
 \label{ma3}
  f^{(4)}= |f|^{p-1}f(1+o(1)),
   \ee
where $\lim \sup f(y)=+\iy$ and $\lim \inf f(y)=-\iy$ as $y \to
y_0^+$, is nonexistent. The proof is easy and follows by
multiplying \ef{ma3} by $f'$ and integrating between two extremum
points $(y_1,y_2)$, where  the former one $y_1$ is chosen to be
sufficiently close to the blow-up value $y_0^+$, whence the
contradiction:
 $$
 \tex{
 0< \frac 12\, (f'')^2(y_1) \sim -\frac 1{p+1}\, |f|^{p+1}(y_1)<0.
 }
 $$

We first study this set of blow-up solutions. These results are
well understood for such fourth-order ODEs; see \cite{Gaz06}, so
we omit some details.

\begin{proposition}
 \label{Pr.Inf1}
  The set of blow-up solutions $\ef{ma3}$ is four-dimensional.
   \end{proposition}

\begin{figure}
%  \vskip -.3cm
\centering
\includegraphics[scale=0.65]{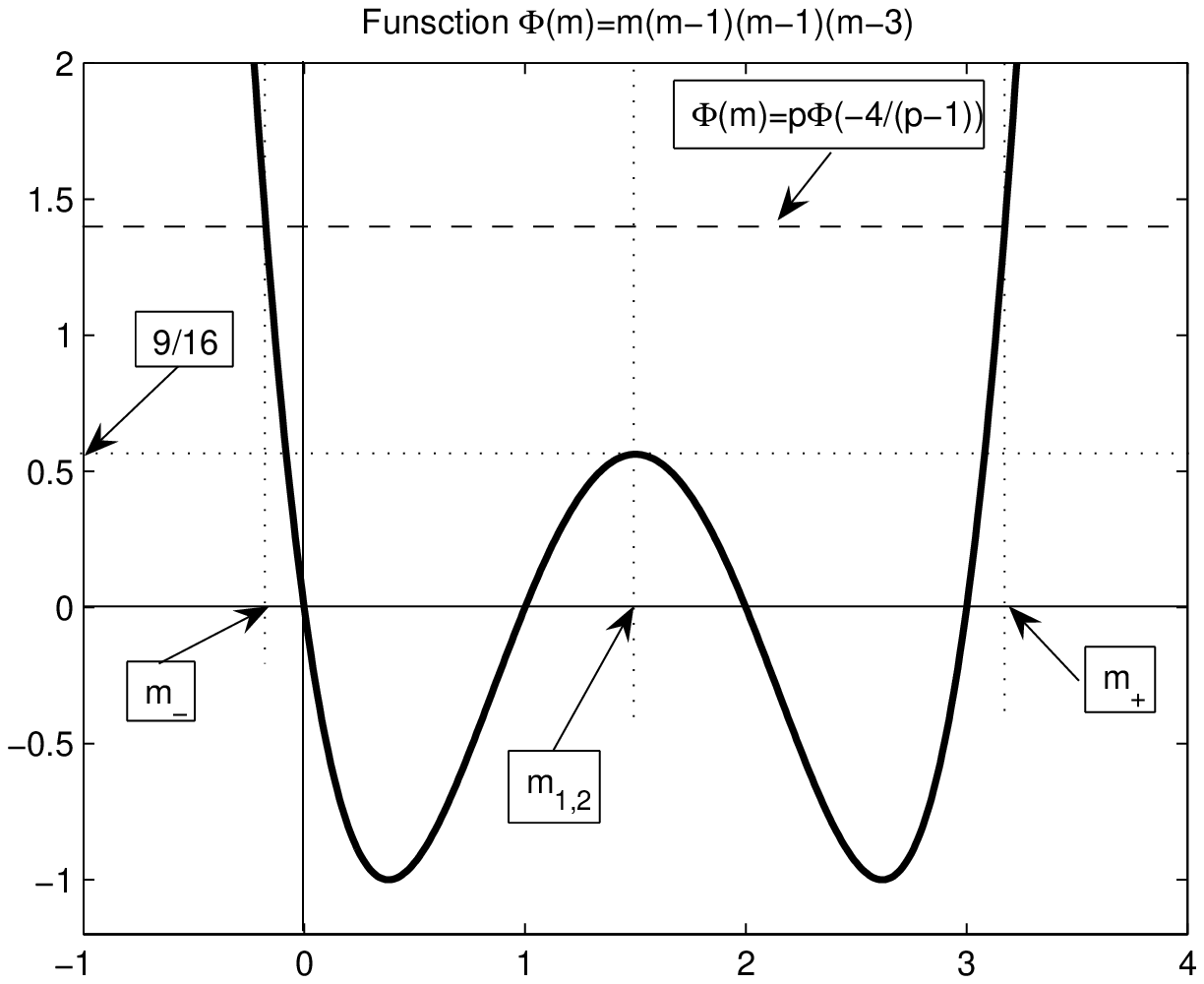} %%%%{AmFBP1.eps}  %%%%%%%{4het5.eps}   %%%%%{4het.eps}  Old
\vskip -.3cm \caption{\small The graph of function $\Phi(m)$ in
\ef{ma4}, and towards solutions of \ef{ma5}.}
%%   \vskip -.3cm
 \label{FPh}
\end{figure}

\noi{\em Proof.} The first parameter  is $y_0 \in \re$. Other are
obtained  from the principal part of the equation \ef{ma3}
describing blow-up via \ef{2.3} as $y \to y_0^+$.
%% \be
%% \label{ma3}
%%  f^{(4)}= |f|^{p-1}f(1+o(1)).
 %%  \ee
 We  apply a standard perturbation argument to \ef{ma3}. Omitting
 the $o(1)$-term and assuming that $f>0$, we find its explicit solution
 \be
 \label{ma4}
  \tex{
 f_0(y)=A_0(y-y_0)^{-\frac 4{p-1}}, \,\,\, A_0^{p-1}= \Phi(-
 \frac 4{p-1}), \,\,\, \Phi(m) \equiv m(m-1)(m-2)(m-3).
 }
 \ee
 For convenience, the graph of $\Phi(m)$ is shown in Figure
 \ref{FPh}. Note that it is symmetric relative to $m_0=\frac 32$,
 at which $\Phi(m)$ has a local maximum:
  \be
  \label{loc1}
   \tex{
   \Phi(\frac 32)= \frac 9{16}.
   }
   \ee

%%%%%%%%%%%%%%%%%%%%%%%%%%%%%%%%%%%%%%%%%%%%%%%%%%%%%%%%%%%%

By linearization, $f=f_0+Y$, we get Euler's ODE:
 \be
 \label{ma5}
  \tex{
 (y-y_0)^4 Y^{(4)}= p A_0^{p-1} \equiv p \Phi(-
 \frac 4{p-1}).
 }
  \ee
It follows that the general solution is composed from the
polynomial ones with the following characteristic equation:
 \be
 \label{ma6}
  \tex{
  Y(y)= (y-y_0)^m \LongA \Phi(m)=p \Phi(-
 \frac 4{p-1}).
  }
  \ee
 Since the multiplier $p>1$ in the last term in \ef{ma6} and $m=- \frac 4{p-1}$
 is a solution if this ``$p \times$" is omitted, this algebraic equation for $m$ admits a
 unique positive solution $m_+>3$, a negative one $m_- < - \frac 4{p-1}$, which is
 not acceptable by \ef{ma4}, ad two complex roots $m_{1,2}$ with ${\rm Re}\, m_{1,2}= \frac 32>0$.
 Therefore, the general solution of \ef{ma3} about the blow-up one
 \ef{ma4}, for any fixed $y_0$, has a 3D stable manifold. $\qed$
 %% \be
 %% \label{ma7}
 %%  f(y)= A_0(y_0-y)^{-\frac 4{p-1}} + A_*(y_0-y)^{m_+}+...\,
 %%  \whereA A_* \in \re,
 %%   \ee
 %%which completes the proof. $\qed$

 \ssk

Thus, according to Proposition \ref{Pr.Inf1}, the blow-up
behaviour with a fixed sign  \ef{ma2} (i.e., non-oscillatory) is
generic for the ODE \ef{2.3}. However, this 4D blow-up bundle
together with the 2D bundle of good solutions \ef{dd1} as $y \to
\pm \iy$ are not enough to justify the shooting procedure. Indeed,
by a straightforward dimensional estimate, an extra bundle at
infinity is missing.

\begin{figure}
%  \vskip -.3cm
\centering
\includegraphics[scale=0.65]{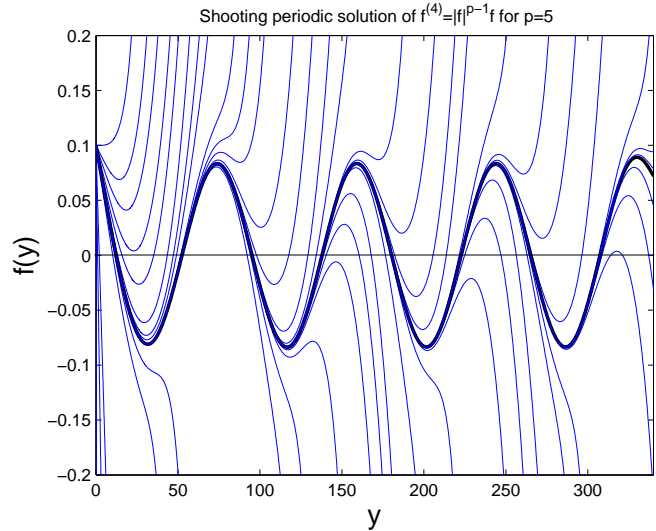} %%%%{AmFBP1.eps}  %%%%%%%{4het5.eps}   %%%%%{4het.eps}  Old
\vskip -.3cm \caption{\small A periodic solution of
$f^{(4)}=|f|^{p-1}f$ as a separatrix: $p=5$.}
%%   \vskip -.3cm
 \label{FVar}
\end{figure}

To introduce this new oscillatory bundle, we begin with the
simpler ODE \ef{ma3}, without the $o(1)$-term, and present in
Figure \ref{FVar} the results of shooting of a ``separatrix" that
lies between orbits, which blow-up to $\pm \iy$. Obviously, this
separatrix is a periodic solution of this equation with a
potential operator. Such variational problems are known to admit
periodic solutions of arbitrary period.

%%%%%%%%%%%%%%%%%%%%%%%%%%%%%%%%%%%%%%%%%%%%%%%%%%%%%%%%%%%%
%%%%%%%%%%%%%%%%%%%%%%%%%%%%%%%%%%%%%%%%%%%%%%%%%%%%%%%%%%%

\begin{figure}
%  \vskip -.3cm
\centering
\includegraphics[scale=0.65]{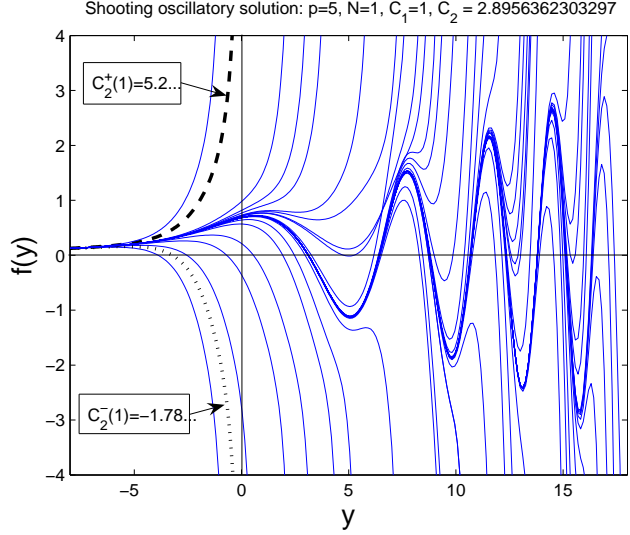} %%%%{AmFBP1.eps}  %%%%%%%{4het5.eps}   %%%%%{4het.eps}  Old
\vskip -.3cm \caption{\small Shooting an oscillatory  solution at
infinity  of \ef{2.3}: $p=5$.}
%%   \vskip -.3cm
 \label{FF1}
\end{figure}

\begin{figure}
%%\vskip -.3cm
\centering \subfigure[$p=7$]{
\includegraphics[scale=0.52]{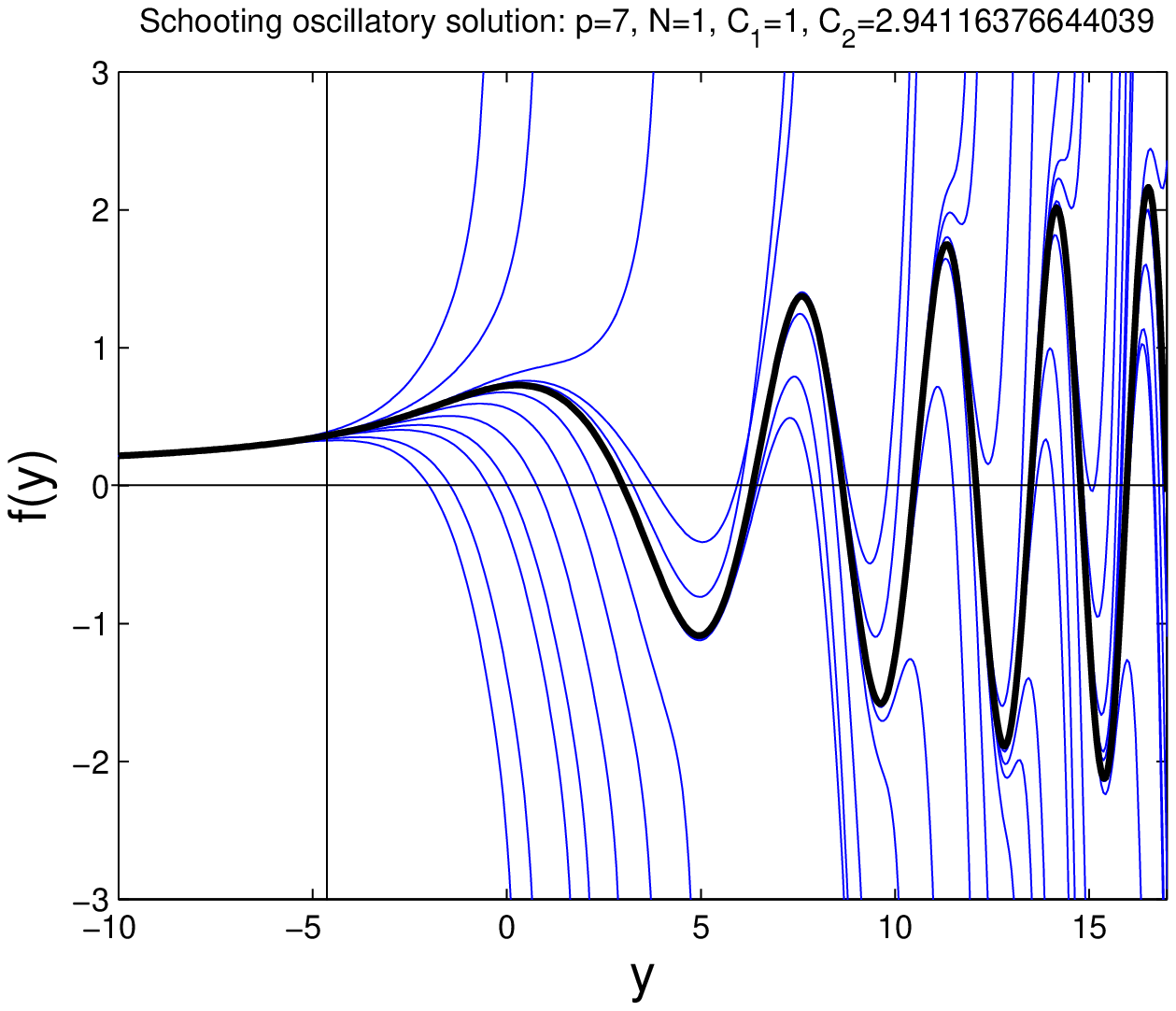}
} \subfigure[$p=2$]{
\includegraphics[scale=0.52]{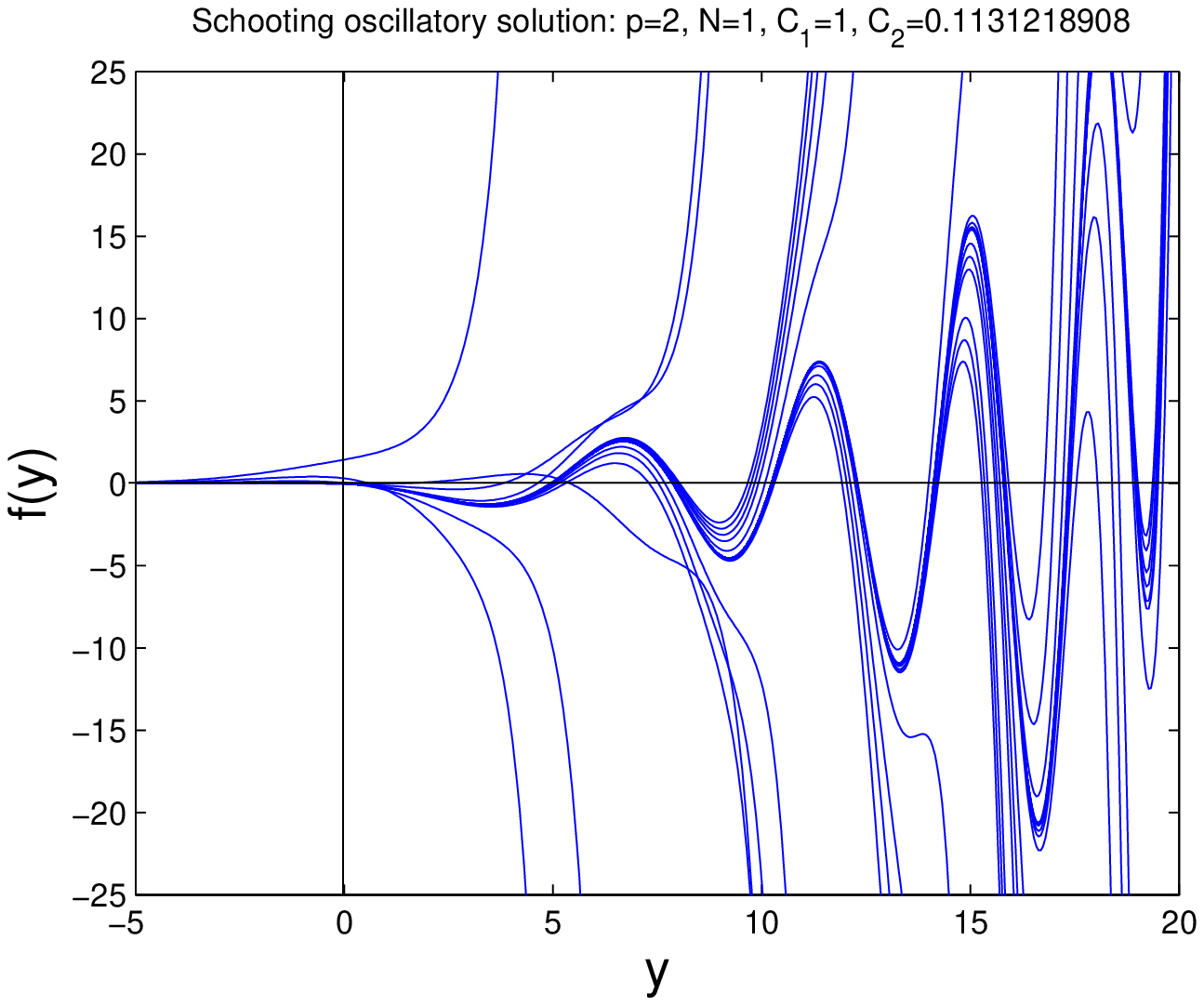}
}
 \vskip -.2cm
\caption{\rm\small  Shooting an oscillatory  solution at infinity
of (\ref{2.3}): $p=7$ (a) and $p=2$ (b).}
 %% \vskip -.3cm
 \label{FF2}
\end{figure}

Thus, Figure \ref{FVar} fixed a bounded oscillatory (periodic)
solution as $y \to +\iy$. When we return to the original equation
\ef{2.3}, which is not variational, we still are able to detect a
more complicated oscillatory structures at $y=\iy$. Namely, these
are generated by the principal terms in
 \be
 \label{pp1}
  \tex{
 f^{(4)} = - \frac 14\, f'y+ |f|^{p-1}f+... \asA y \to \iy.
 }
  \ee
Similar to Figure \ref{FVar}, in Figure \ref{FF1}, we present the
result of shooting (from $y=-\iy$, which is the same by symmetry)
of such oscillatory solutions of \ef{2.3} for $p=5$. It is easy to
see that such oscillatory solutions have increasing amplitude of
their oscillations as $y \to \iy$, which, as above, is proved by
multiplying \ef{pp1} by $f'$ and integrating over any interval
$y_1,y_2)$ between two extrema. Figure \ref{FF2} shows shooting of
similar oscillatory structures at infinity for $p=7$ (a) and $p=2$
(b). It is not very difficult to prove that the set of such
oscillatory orbits at infinity is 1D and this well corresponds to
the periodic one in Figure \ref{FVar} depending on the single
parameter being its arbitrary period.

 By $C_2^\pm(C_1)$ in Figure \ref{FF1}, we denote the values of
 the second parameters $C_2$ such that, for a fixed $C_1 \in \re$,
 the solutions $f(y;C_1,C_2^\pm)$ blow up to $\pm \iy$ respectively.
These values are necessary for shooting the symmetry conditions
\ef{2.31}.

%%%%%%%%%%%%%%%%%%%%%%%%%%%%%%%%%%%%%%%%%%%%%%%%%%%%%%%%%%%%
%%%%%%%%%%%%%%%%%%%%%%%%%%%%%%%%%%%%%%%%%%%%%%%%%%%%%%%%%%%

%%%%%%%%%%%%%%%%%%%%%%%%%%%%%%%%%%%%%%%%%%%%%%%%%%%%%%%%%%%%%%%
%%FIG%%%%%%%%%%%%%%%%%%%%%%%

%%%%%%%%%%%%%%%%%%%%%%%%%%%%%%%%%%%%%%%%%%%%%%%%%%%%%%%%%%%%%%%%%%%

Thus, overall, using two parameters $C_{1,2}$ in the bundle
\ef{dd1} for $y \gg 1$ leads to a well-posed problem of a 2D--2D
shooting:
 \be
 \label{w1}
  \mbox{find $C_{1,2}$ such that:} \quad
  \left\{
   \begin{matrix}
    y_0(C_1,C_2)=-\iy, \quad \mbox{and} \qquad\qquad\quad\\
    \mbox{no oscillatory behaviour as $y \to -\iy$}.
     \end{matrix}
     \right.
   \ee

Concerning the actual proof of existence via shooting of at least
a single blow-up patterns $f_0(y)$,  by construction and
oscillatory property of the equation \ef{2.3}, we first claim that
in view of continuity relative to the parameters,
 \be
 \label{ccl1}
 \mbox{for any $C_1>0$, there exists $C_2^*(C_1) \in (C_2^-(C_1),C_2^+(C_1))$ such that
 $f'''(0)=0$}.
  \ee
We next  change $C_1$ to prove that at this $C_2^*(C_1)$ the
derivative $f'(0)$ also changes sign. Indeed, one can see that
 \be
 \label{ccl2}
 f'(0;C_1,C_2^*(C_1)) > 0 \forA C_1\ll 1 \andA
f'(0;C_1,C_2^*(C_1)) < 0 \forA C_1\gg 1.
 \ee
 Actually, this means for such essentially different values of
 $C_1$, the solution $f(y;C_1,C_2^*(C_1)$ has first oscillatory ``humps"
 for $y>0$ and $y<0$ respectively. By continuity in $C_1$,
 \ef{ccl2} implies existence of a $C_1^*$ such that
  \be
  \label{ccl3}
  f'(0;C_1^*,C_2^*(C_1^*))=0,
  \ee
  which together with \ef{ccl1} induced the desired solution.
  Overall, the above geometric shooting well corresponds to that applied in the
  standard framework of classic ODE theory, so we do not treat
  this in greater detail. However, we must admit that proving analogously existence of
  the second solution $f_0(y)$ (detected earlier by
  not fully justified arguments of
   homotopy and branching theory and confirmed numerically) is an
   open problem. A more difficult open problem is to show why the
   problem \ef{w1} does not admit non-symmetric (non-even) solutions $f(y)$
   (or does it?).

%%%%%%%%%%%%%%%%%%%%%%%%%%%%%%%%%%%%%%%%%%%%%%%%%%%%%%%%%%%%
\subsection{Dimensions $N \ge 2$: on 2D shooting and analogous nonuniqueness}

In higher dimensions, it is easier to describe Type I(ss) blow-up
in radial geometry, where \ef{2.2} also becomes an ODE of the form
(now $y$  stands for $|y|>0$)
 \be
 \label{2.4}
 \tex{
%% {\bf A}(f) \equiv
-f^{(4)} - \frac {2(N-1)}y\, f'''- \frac
 {(N-1)(N-3)}{y^2}\, f'' + \frac{(N-1)(N-3)}{y^3}\, f' - \frac
 14\, y f'
 - \frac 1{p-1}\,
  f + |f|^{p-1}f=0
  %% \inB \re_+,
  %% f'(0)=f'''(0)=0
  %%\,\,\,(\mbox{symmetry)}.
  }
  \ee
  in $\re$,
with the same two symmetry condition \ef{2.31}.
 To explain the nature of difficulties in proving existence of
 solutions of \ef{2.4}, let us describe the admissible
 behaviour for $y \gg 1$. There exists a 2D bundle of such
 asymptotics (see details in \cite[\S~3.3]{Bl4}): as $y \to +\iy$,
  \be
  \label{dd1}
   \tex{
f(y) = \big[C_1 y^{-\frac 4{p-1}}+... \big]+ \big[C_2 y^{-\frac
23(N-\frac 2{p-1})}\, {\mathrm e}^{-a_0 y^{ 4/3}}+... \big]
\whereA a_0=3 \cdot 2^{-\frac 83}
 }
 \ee
and $C_1$ and $C_2$ are arbitrary parameters. This somehow reminds
a typical centre manifold structure of the origin $\{f=0\}$ at $y=
\iy$: the first term in \ef{dd1} is a node bundle with algebraic
decay, while the second one corresponds to ``non-analytic"
exponential bundle around any  of algebraic curves.
%% Finally, we
%%arrive at a
Thus, a dimensionally well-posed shooting is:
 \be
 \label{dd2}
  \fbox{$
 \mbox{{\bf Shooting:} \,\, using {\bf 2} parameters $C_{1,2}$ in (\ref{dd1})
to satisfy {\bf 2} conditions (\ref{2.31}).} $}
 \ee
In case of analytic dependence of solutions of \ef{2.4} on
parameters $C_{1,2}$ in the bundle \ef{dd1} (this is rather
plausible via standard trends of ODE theory, but difficult to
prove), the problem cannot have more than a countable set of
solutions. Actually, our numerics confirm that in wide parameter
ranges of $p>1$ and $N \ge 1$, there exist not more than two
solutions (up to other more unstable ones about the SSS; see
Section \ref{S.5}):
 \be
 \label{dd2N}
 f_0(y) \,\,\, \mbox{with} \,\,\, \{C_{10}(p,N),C_{20}(p,N)\},
 \andA f_1(y) \,\,\, \mbox{with} \,\,\, \{C_{11}(p,N),  C_{21}(p,N)\}.
  \ee
  The rest of this section is devoted to justify this.

The eventual similarity blow-up patterns can be characterized by
their {\em final  time profiles}: passing to the limit $t \to T^-$
in \ef{2.1} and using the expansion \ef{dd1} yields
 \be
 \label{dd3}
 \mbox{if $C_1(p,N) \not = 0$, then} \quad u(x,t) \to C_1
 |x|^{-\frac 4{p-1}} \asA t \to T^-
  \ee
  uniformly on any compact subset of $\ren \setminus \{0\}$.
If $C_1 =0$ in \ef{dd1}, i.e., $f(y)$ has an exponential decay at
infinity, then the limit is different: in the sense of
distributions,
\be
 \label{dd4}
  \tex{
 C_1(p,N)  = 0 \,\Longrightarrow \,
 \fbox{$
  |u(x,t)|^{\frac{N(p-1)}4} \to
 C_3 \d(x), \,\,
  t \to T^-; \,\,\, C_3= \int |f|^{\frac{N(p-1)}4} < \iy.
  $}
 }
  \ee
 It is very difficult to prove that \ef{dd4} actually takes place at some $p=p_\d(N)>1$
 (even for $N=1$),
 and we will justify this numerically for some not that large dimensions $N\le 11$.

\begin{figure}
%%\vskip -.3cm
\centering \subfigure[$p=2, \,\, N=3$]{
\includegraphics[scale=0.52]{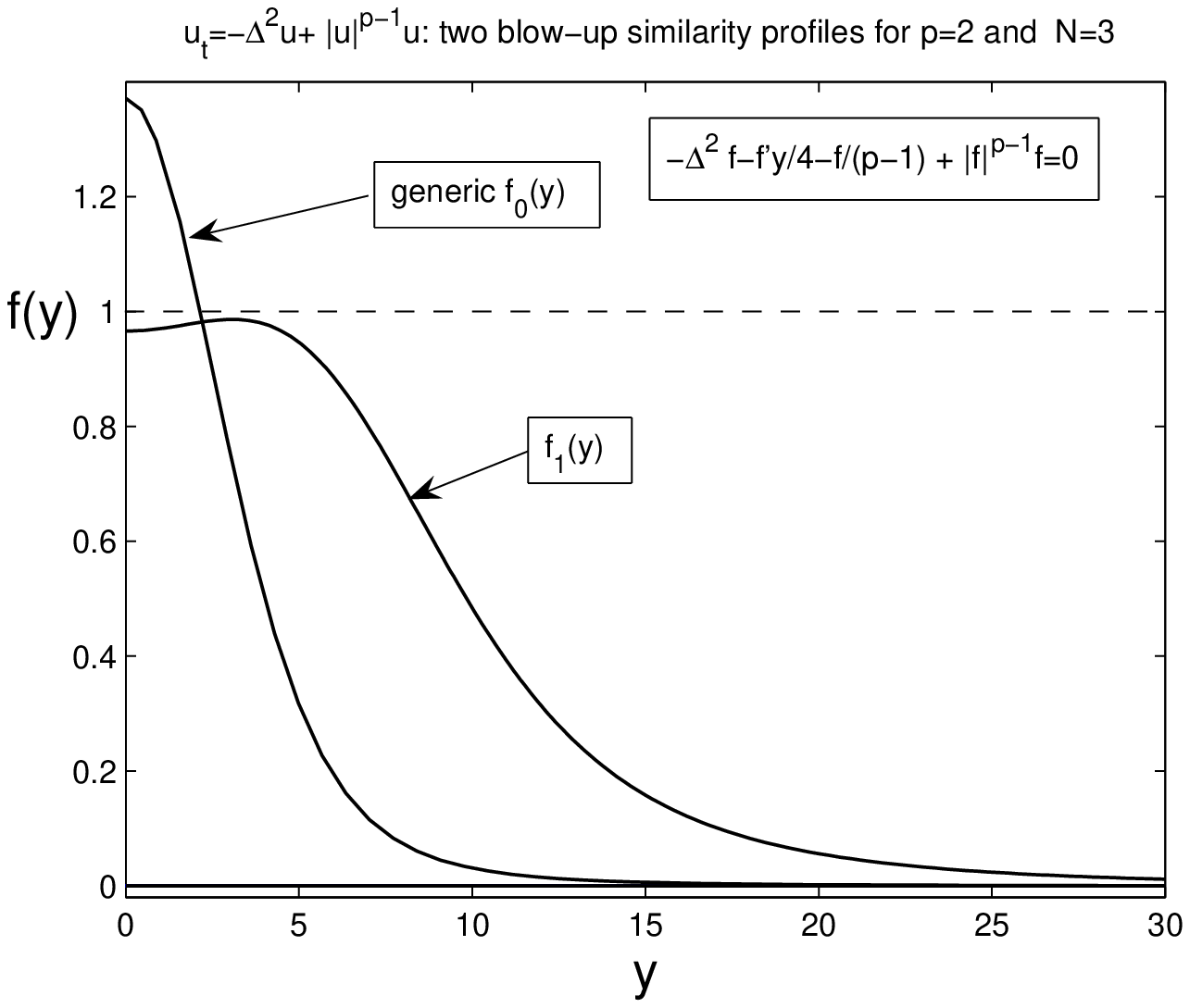}
} \subfigure[$p=2,\,\,N=10$]{
\includegraphics[scale=0.52]{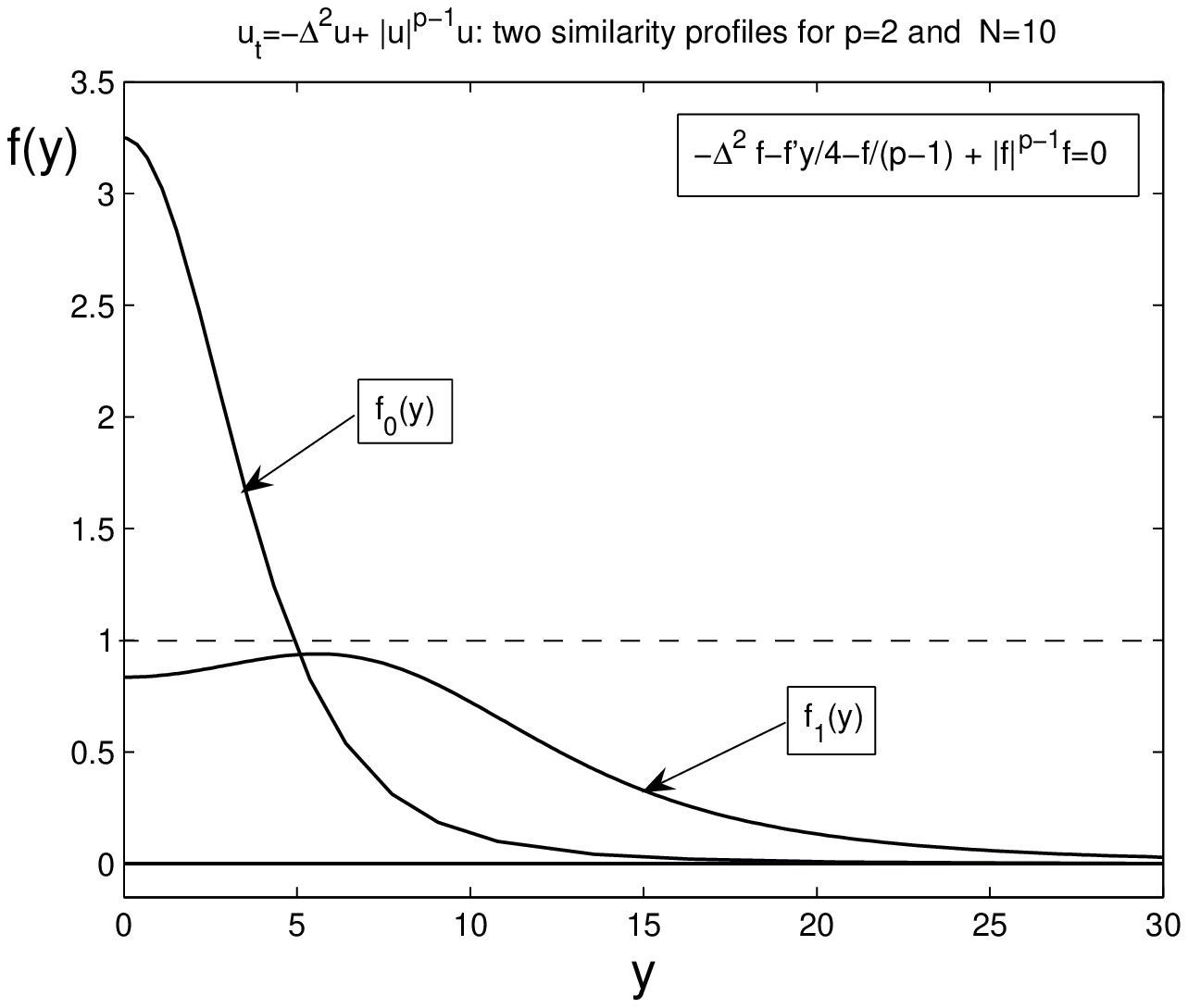}
}
 \vskip -.2cm
\caption{\rm\small Two self-similar blow-up solutions of \ef{2.4}
for $p=2$: $N=3$ (a) and $N=10$ (b).}
 %% \vskip -.3cm
 \label{F3}
\end{figure}

We now start describe various similarity blow-up profiles for $N
\ge 2$.
 As a first and
analogous to $N=1$ example, in Figure \ref{F3}, we construct
numerically first two profiles, $f_0(y)$ and $f_1(y)$, for the
three-dimensional case $N=3$ and $N=10$ for $p=2$, which look
rather similar to those in Figure \ref{F1} for $N=1$.
 %%%Figure
%%\ref{F31} shows similar profiles for $p=5$.

%%FIG%%%%%%%%%%%%%%%%%%%%%%%

%%%%%%%%%%%%%%%%%%%%%%%%%%%%%%%%%%%%%%%%%%%%%%%%%%%%%%%%%%%%%%%%%%%

\subsection{$N \ge 2$: $p$-branches of the  profile $f_0(y)$ and $f_1(y)$}

\begin{figure}
%%\vskip -.3cm
\centering \subfigure[$f_0$-deformation]{
\includegraphics[scale=0.52]{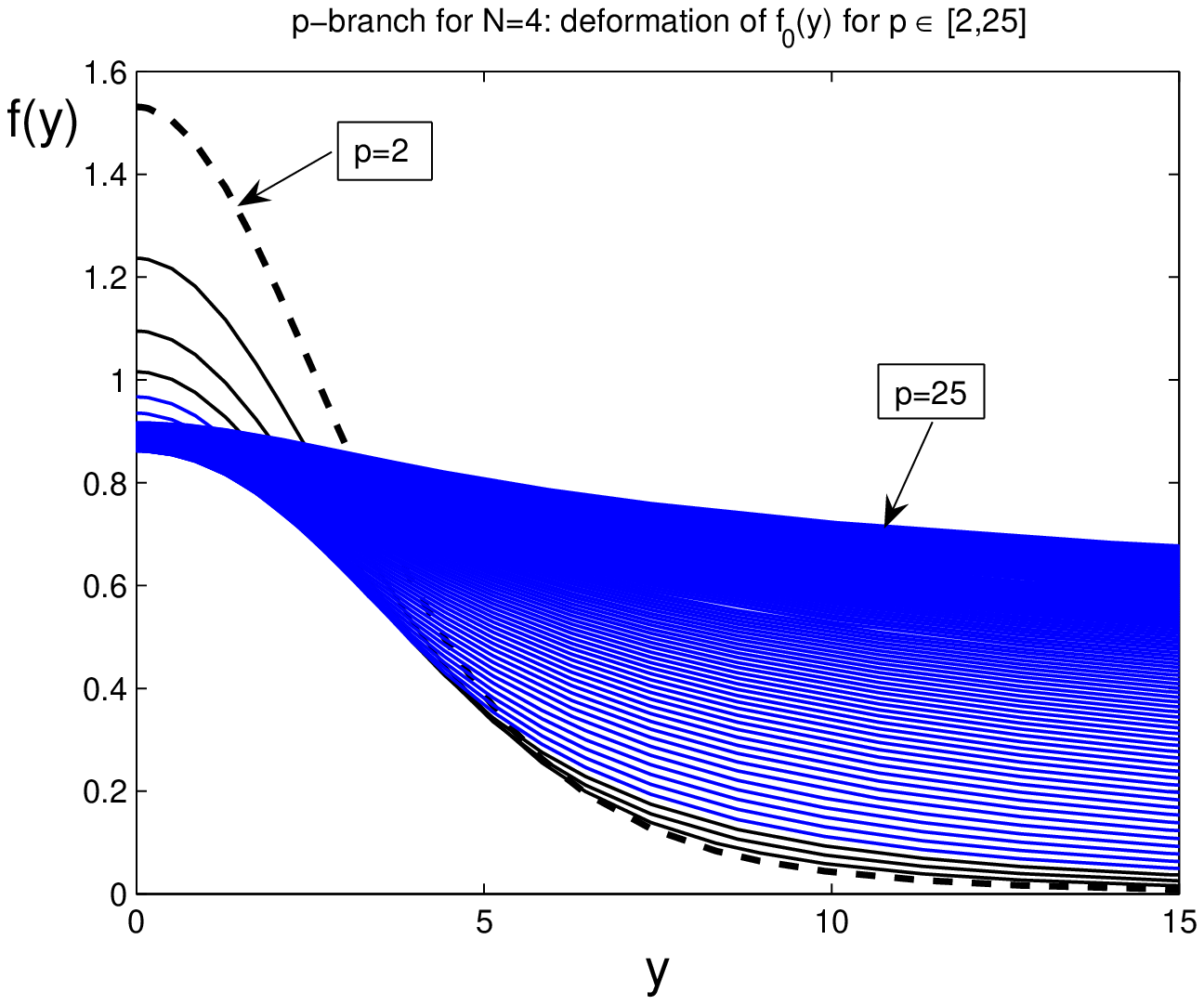}
} \subfigure[$p$-branch]{
\includegraphics[scale=0.52]{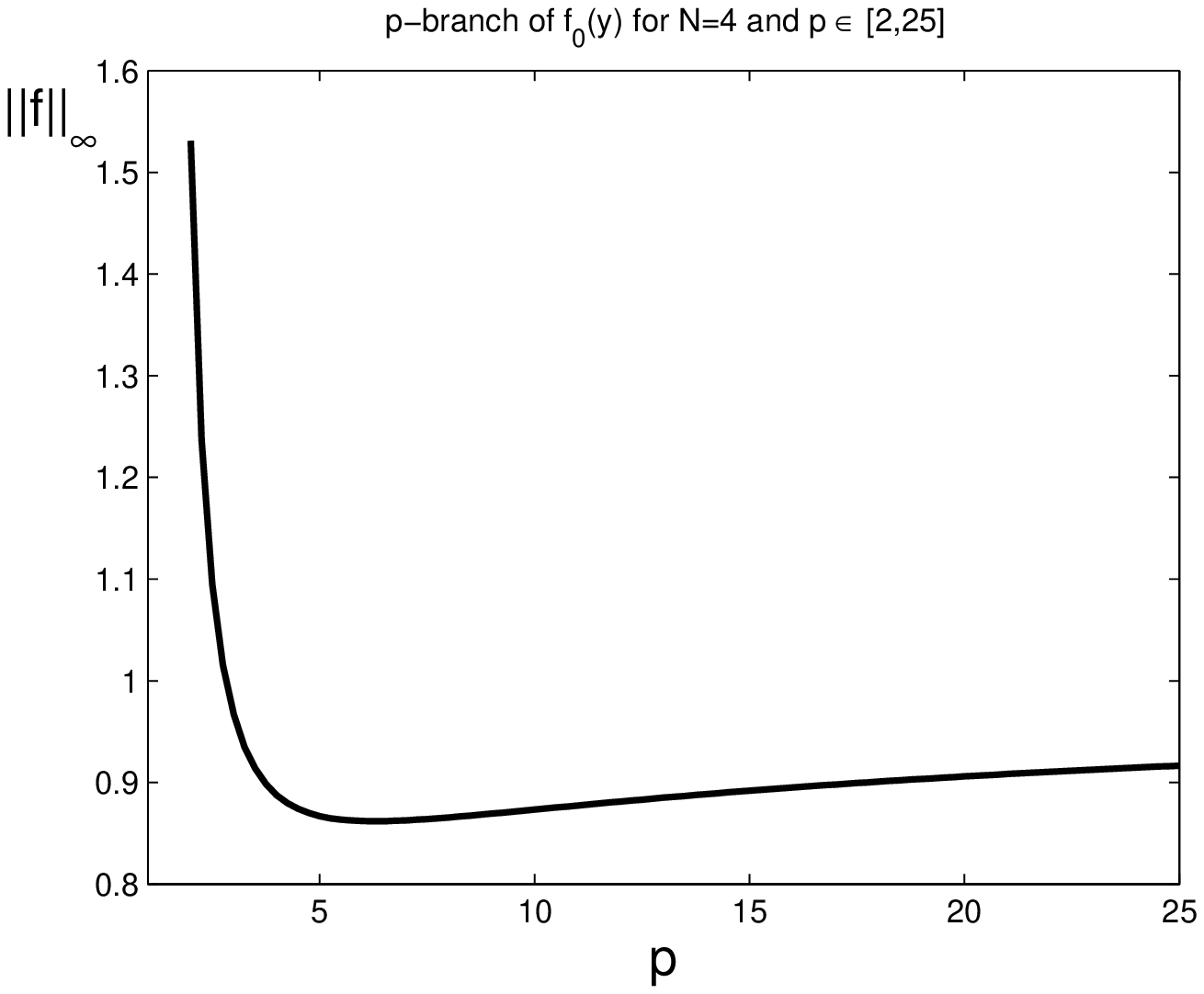}
}
 \vskip -.2cm
\caption{\rm\small The $p$-branch of blow-up self-similar profiles
$f_0(y)$ for $N=4$.}
 %% \vskip -.3cm
 \label{F31}
\end{figure}

\begin{figure}
%%\vskip -.3cm
\centering \subfigure[$f_0(y)$-deformation]{
\includegraphics[scale=0.52]{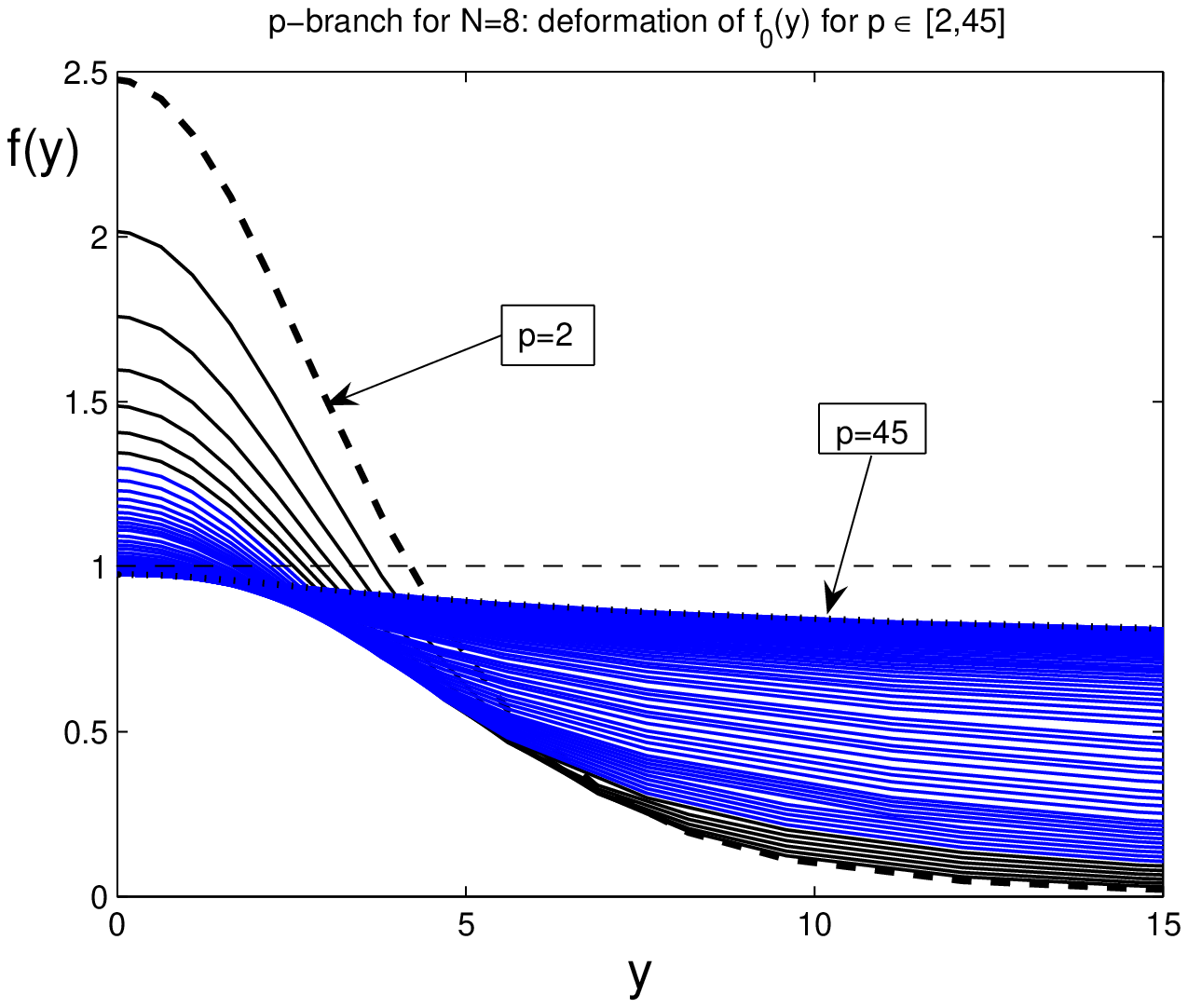}
} \subfigure[$p$-branch]{
\includegraphics[scale=0.52]{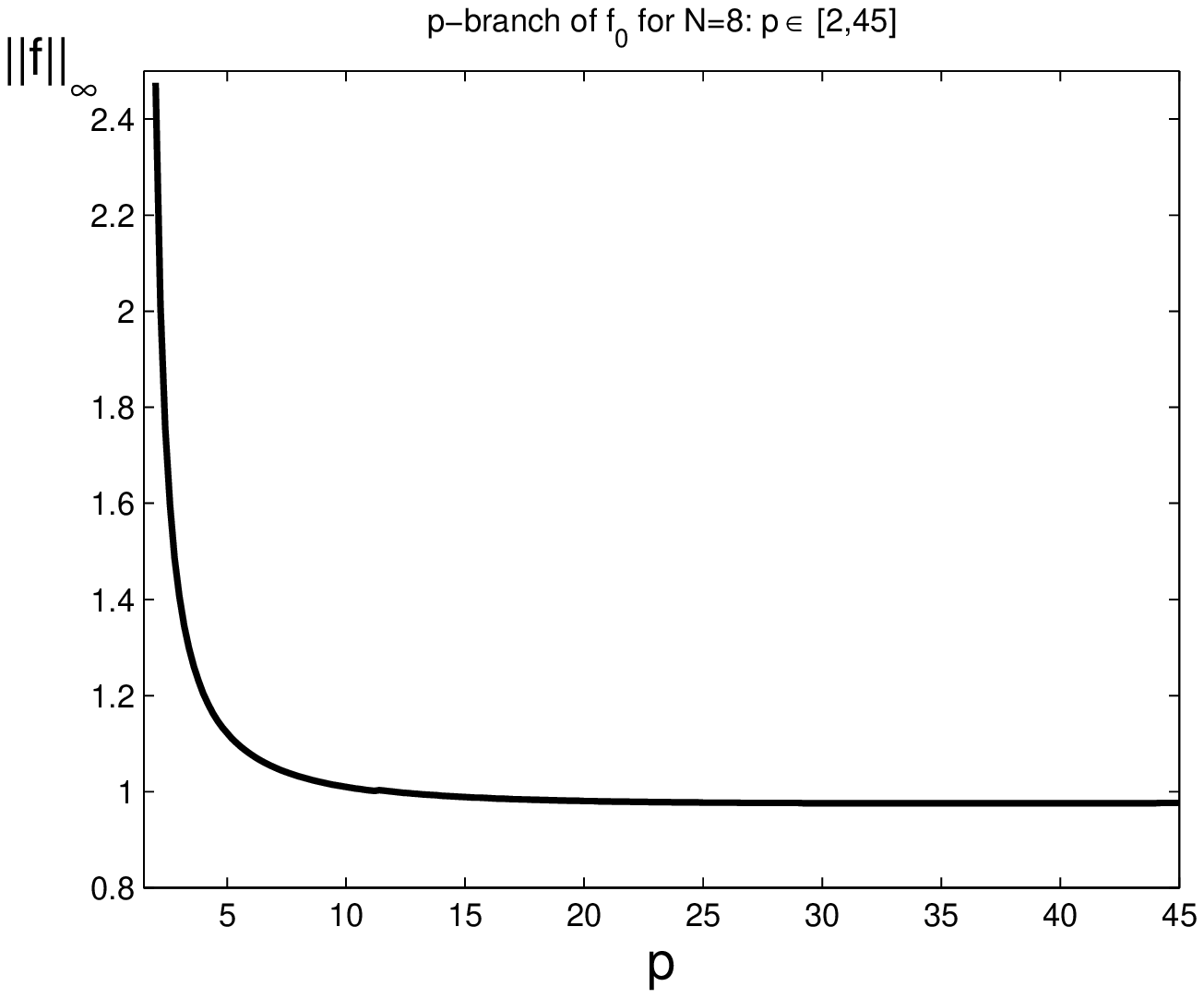}
}
 \vskip -.2cm
\caption{\rm\small The $p$-branch of blow-up self-similar profiles
$f_0(y)$ for $N=8$.}
 %% \vskip -.3cm
 \label{F32}
\end{figure}

\begin{figure}
%%\vskip -.3cm
\centering \subfigure[$p$-branch for $N=1$]{
\includegraphics[scale=0.52]{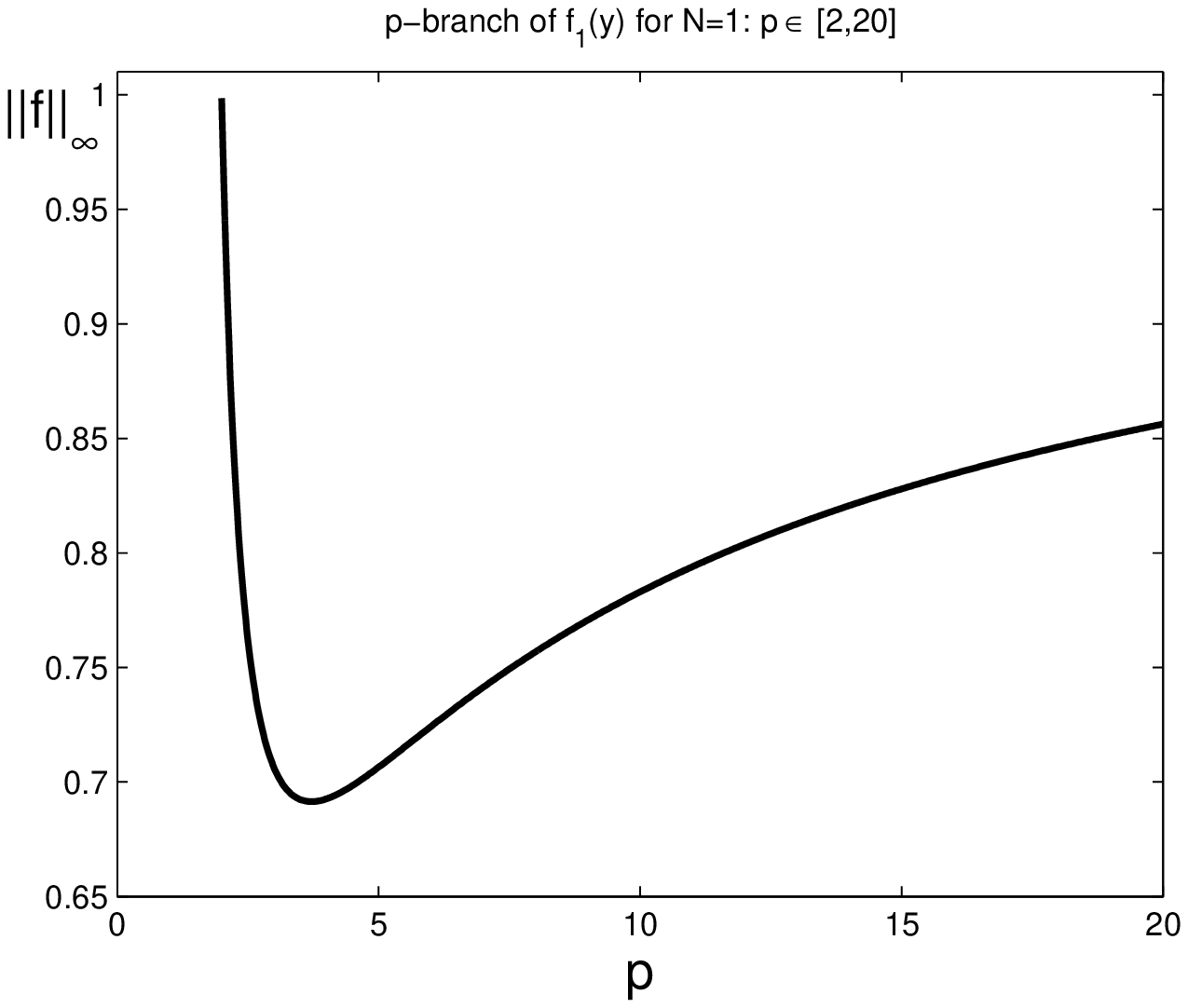}
} \subfigure[$p$-branch for $N=12$]{
\includegraphics[scale=0.52]{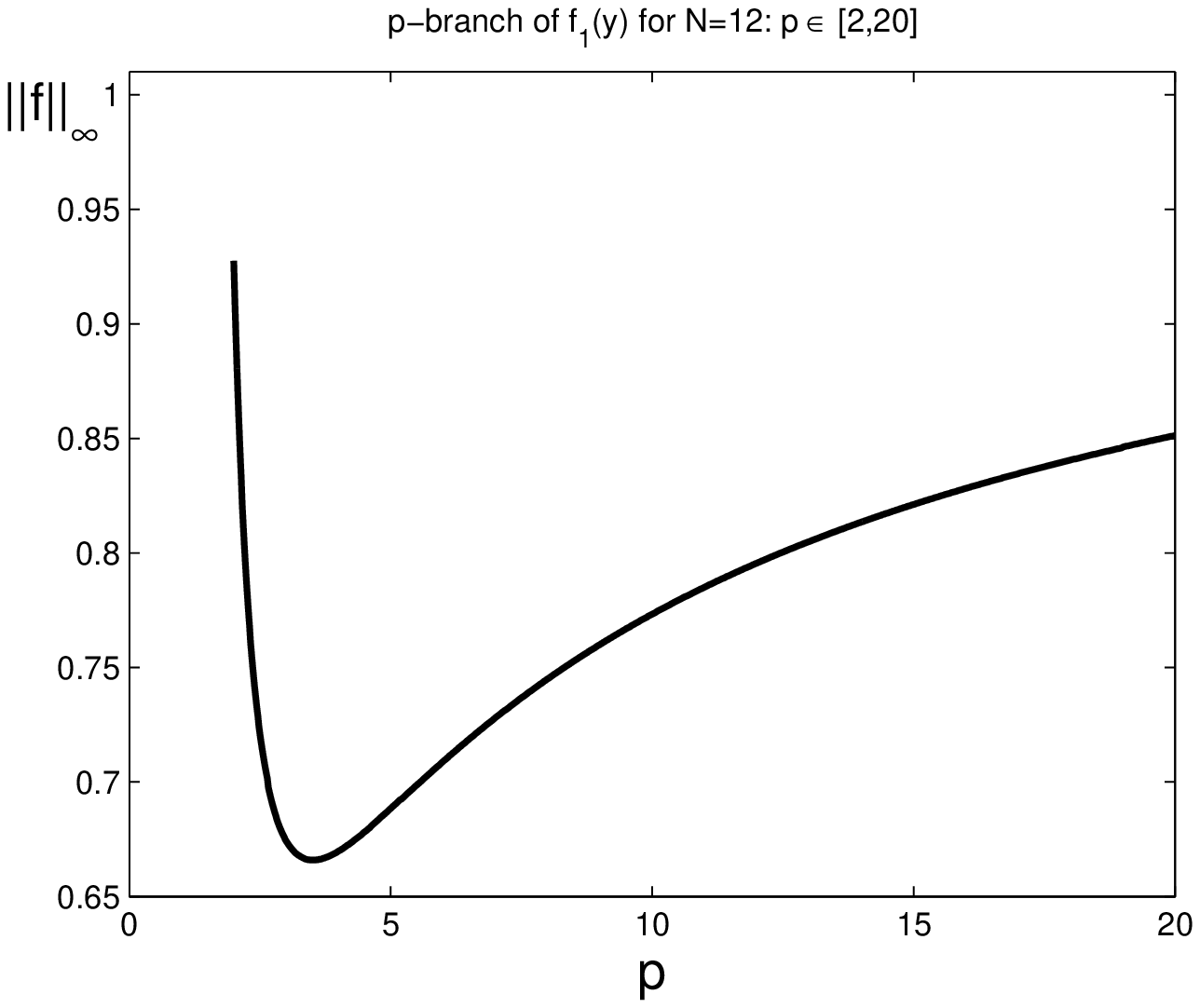}
}
 \vskip -.2cm
\caption{\rm\small The $p$-branches of blow-up self-similar
profiles $f_1(y)$ for $p=2$: $N=1$ (a) and $N=12$ (b).}
 %% \vskip -.3cm
 \label{Ff1s}
\end{figure}

Such $p$-branches of solutions are a convenient way to describe
families of profiles $f_0(y)$ depending on the exponent  $p$; cf.
\cite{GHUni, GW2}. In Figure \ref{F31}, we present such a branch
of $f_0$ for $N=4$, where  (a) shows the actual smooth deformation
of $f_0(y)$ with changing $p$, while (b) is the corresponding
$p$-branch. In Figure \ref{F32}, the same is done for $N=8$. Note
that both Figures  (b) show that $\|f\|_\iy=f(0)$ approaches 1 for
large $p$, which is a general phenomenon for such ODEs described
in \cite[\S~5]{GHUni}.  Similarly, Figure \ref{Ff1s} shows
$p$-branches of the second blow-up profile $f_1(y)$ for $p=2$ in
the cases $N=1$ (a) and $N=12$ (b) (the critical dimension, where
$p_{\rm S}=2$).

%%FIG%%%%%%%%%%%%%%%%%%%%%%%

%%FIG%%%%%%%%%%%%%%%%%%%%%%%

%%FIG%%%%%%%%%%%%%%%%%%%%%%%

\begin{figure}
%%\vskip -.3cm
\centering \subfigure[$f_1$-deformation]{
\includegraphics[scale=0.52]{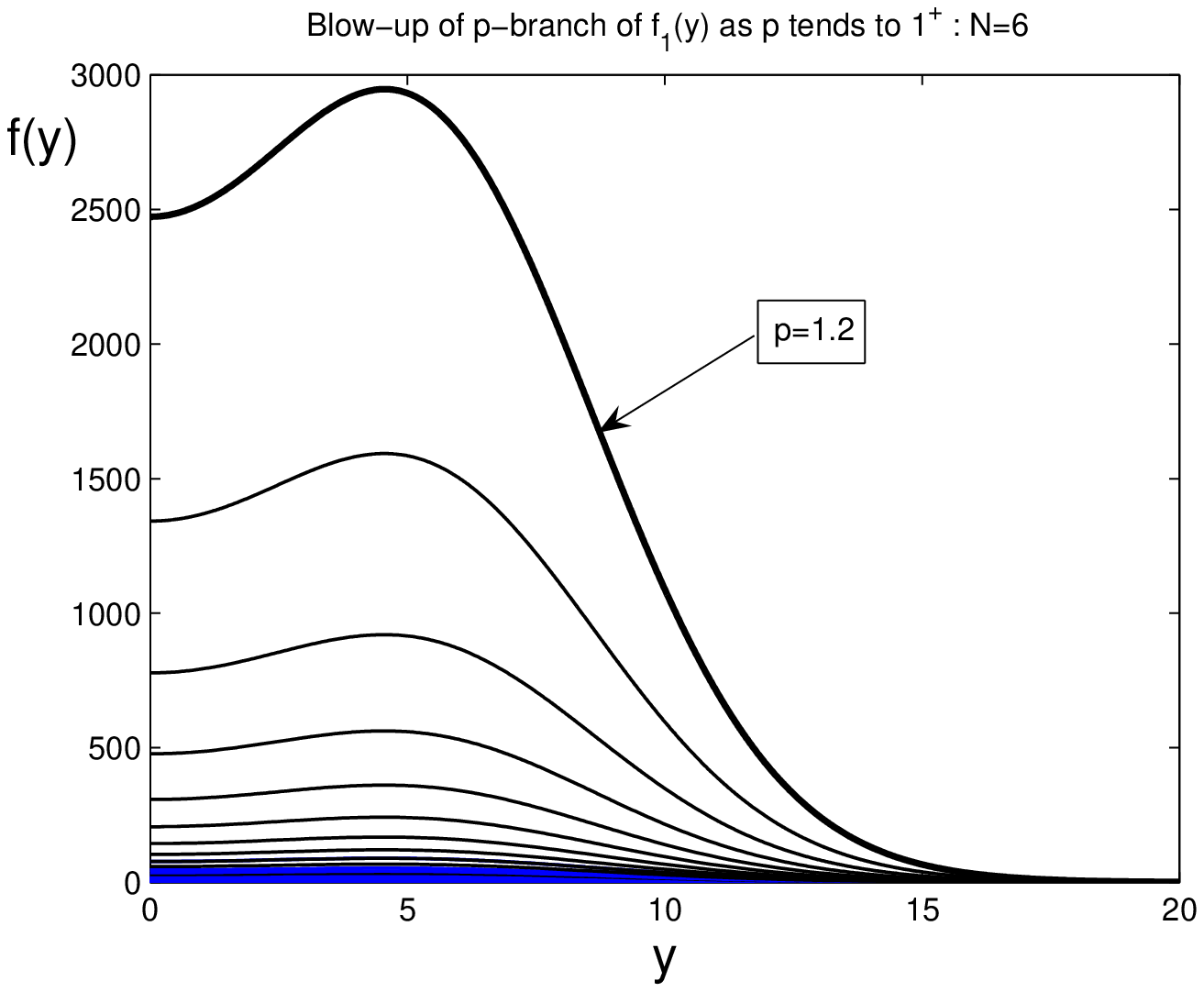}
} \subfigure[$p$-branch]{
\includegraphics[scale=0.52]{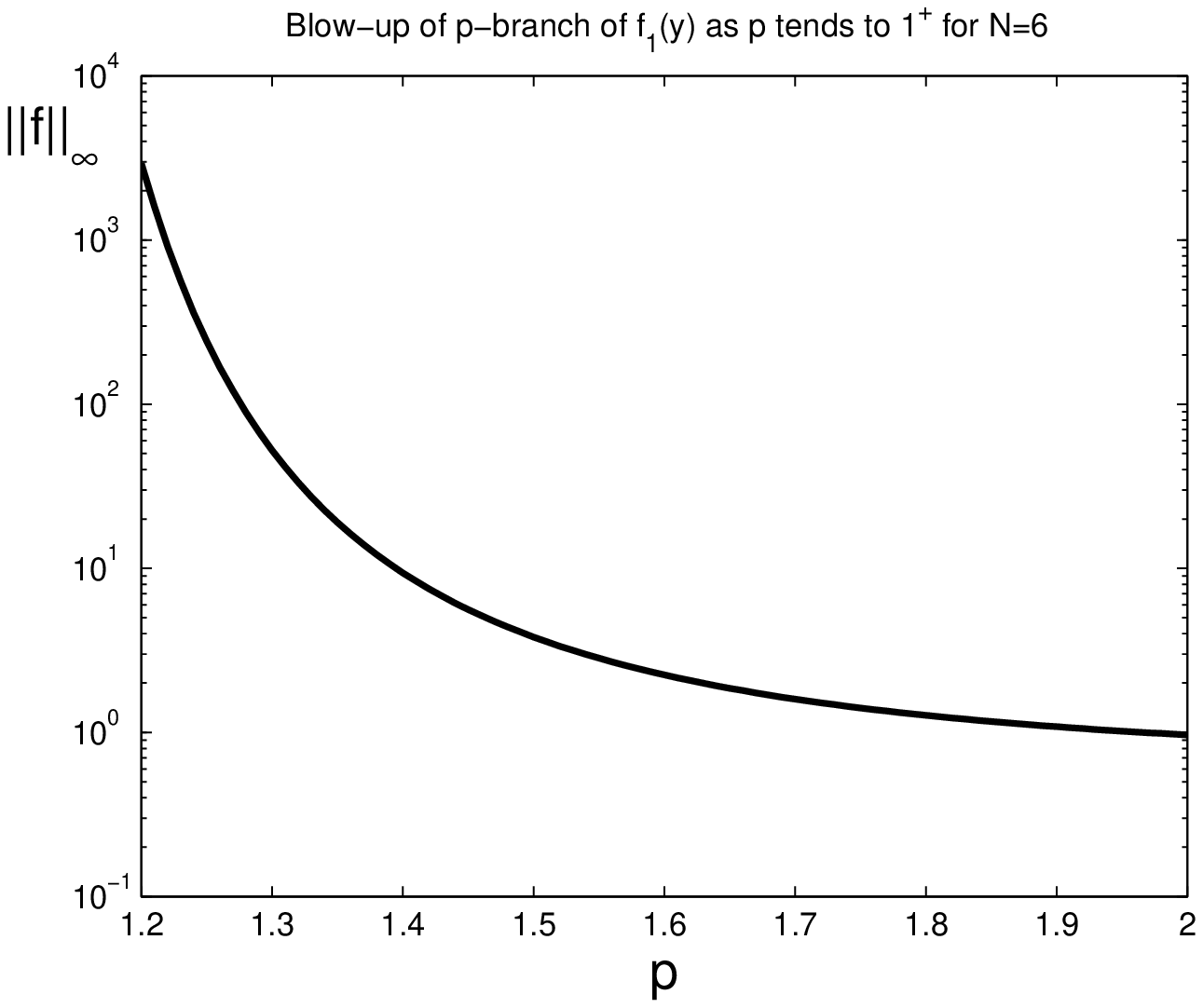}
}
 \vskip -.2cm
\caption{\rm\small Blow-up as $p \to 1^-$ of the $p$-branches of
profiles $f_1(y)$ for $N=6$: $p$-deformation of $f_1(y)$ (a) and
blow-up of $p$-branch in the log-scale (b).}
 %% \vskip -.3cm
 \label{FEx1}
\end{figure}

It is well understood that for equations such as \ef{2.4}, the
solutions $f(y)$ blow-up as $p \to 1^+$ with a super-exponential
rate $\sim (p-1)^{-1/(p-1)}$; see \cite{GHUni, GW2}. As an
example, in Figure \ref{FEx1}, we present such a blowing up
behaviour of the $p$-branch of $f_1(y)$ for $N=6$.

%%%%%%%%%%%%%%%%%%%%%%%%%%%%%%%%%%%%%%%%%%%%%%%%%%
%%FIG%%%%%%%%%%%%%%%%%%%%%%%

%%%%%%%%%%%%%%%%%%%%%%%%%%%%%%%%%%%%%%%%%%%%%%%%%%%%%%%%
\subsection{$N \ge 2$: $N$-branches of blow-up profiles}

\begin{figure}
%  \vskip -.3cm
\centering
\includegraphics[scale=0.75]{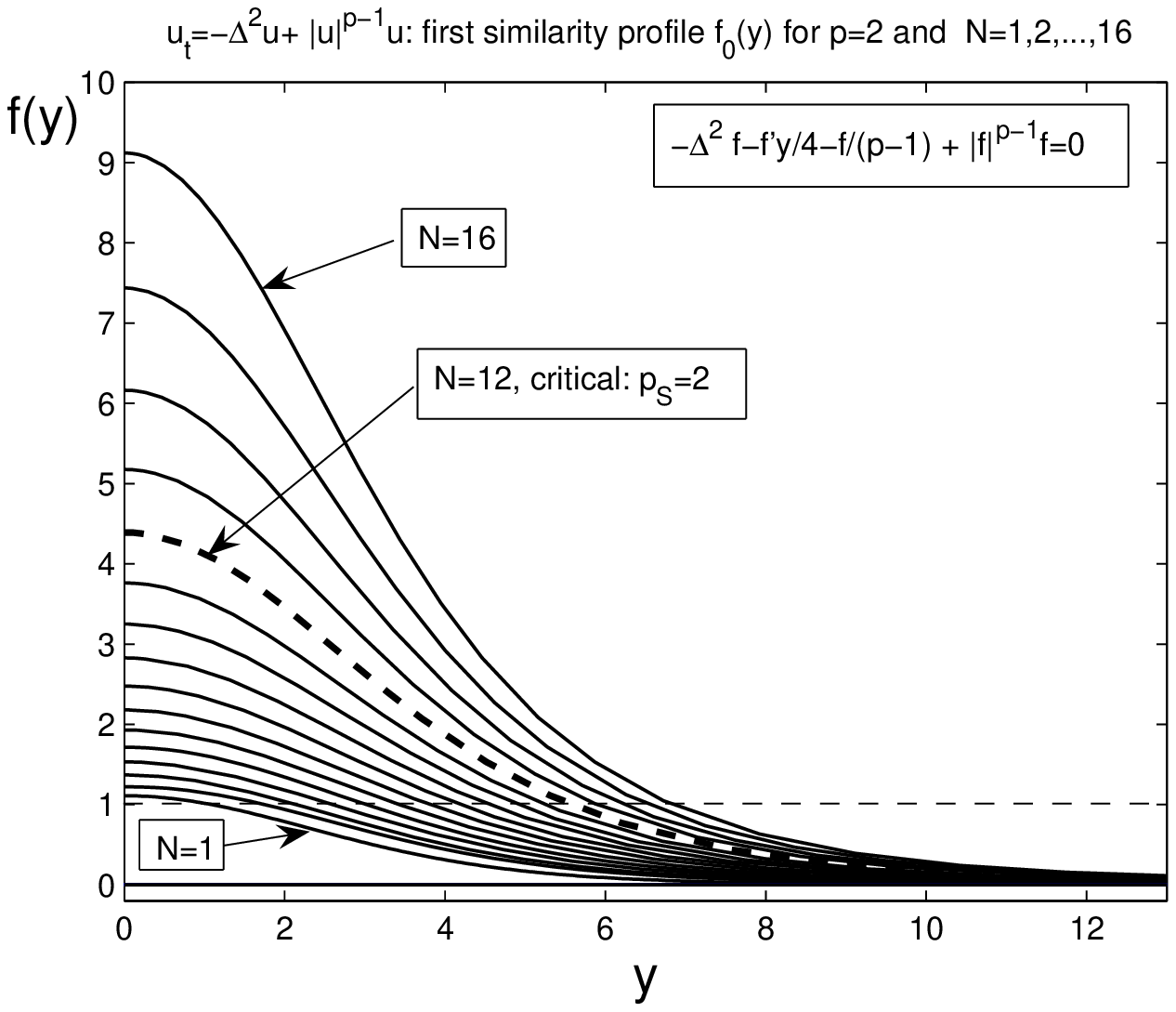} %%%%{AmFBP1.eps}  %%%%%%%{4het5.eps}   %%%%%{4het.eps}  Old
\vskip -.3cm \caption{\small Similarity blow-up solutions $f_0(y)$
of (\ref{2.4}) for $p=2$ in dimensions $N=1,\, 2,\, ... \, ,\,
16$.}
  %%% \vskip -.3cm
 \label{F4}
\end{figure}

\begin{figure}
%  \vskip -.3cm
\centering
\includegraphics[scale=0.75]{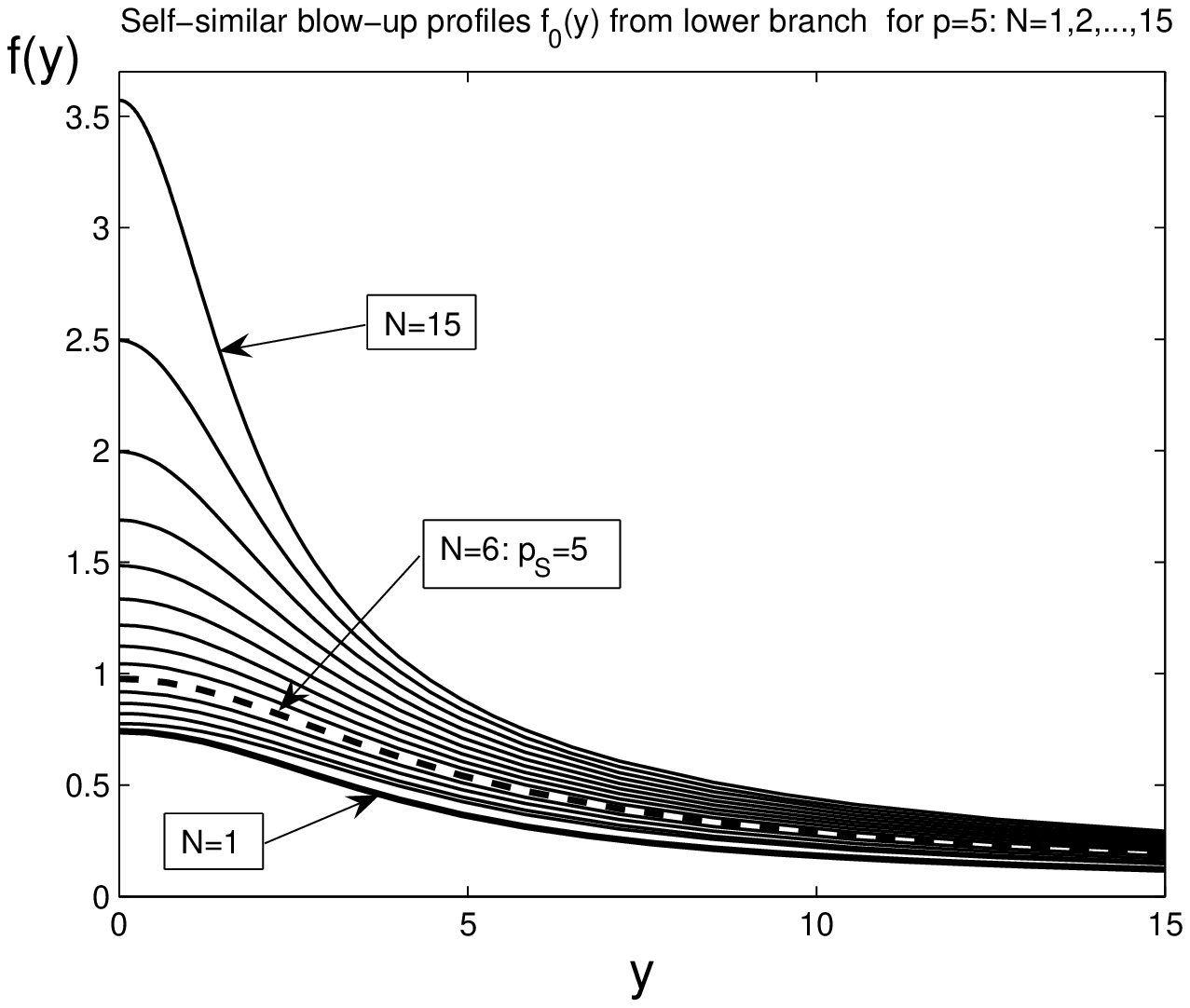} %%%%{ps4N1New.eps} %%%%{AmFBP1.eps}  %%%%%%%{4het5.eps}   %%%%%{4het.eps}  Old
\vskip -.3cm \caption{\small Similarity solutions $f_0$ of
(\ref{2.4}) for $p=5$ and  $N=1,\,2,\, ...\,$, $15$.}
  %% \vskip -.3cm
 \label{F2}
\end{figure}

Firstly, in some $N$-intervals, there is a continuous dependence
of $f_0(y)$ on the dimension, as Figure \ref{F4} clearly shows for
$p=2$ and Figure \ref{F2} for $p=5$ (those values of $p$ will be
constantly used later on for the sake of comparison).

%%%%%%%%%%%%%%%%%%%%%%%%%%%%%%%%%%%%%%%%%%%%%%%%%%%%%%%%%%%%

%%%%%%%%%%%%%%%%%%%%%%%%%%%%%%%%%%%%%%%%%%%%%%%%%%%%%%%%%%%%
%%%%%%%%%%%%%%%%%%%%%%%%%%%%%%%%%%%%%%%%%%%%%%%%%%%%%%%%%%%%%%

\begin{figure}
%%\vskip -.3cm
\centering \subfigure[$p=5$]{
\includegraphics[scale=0.52]{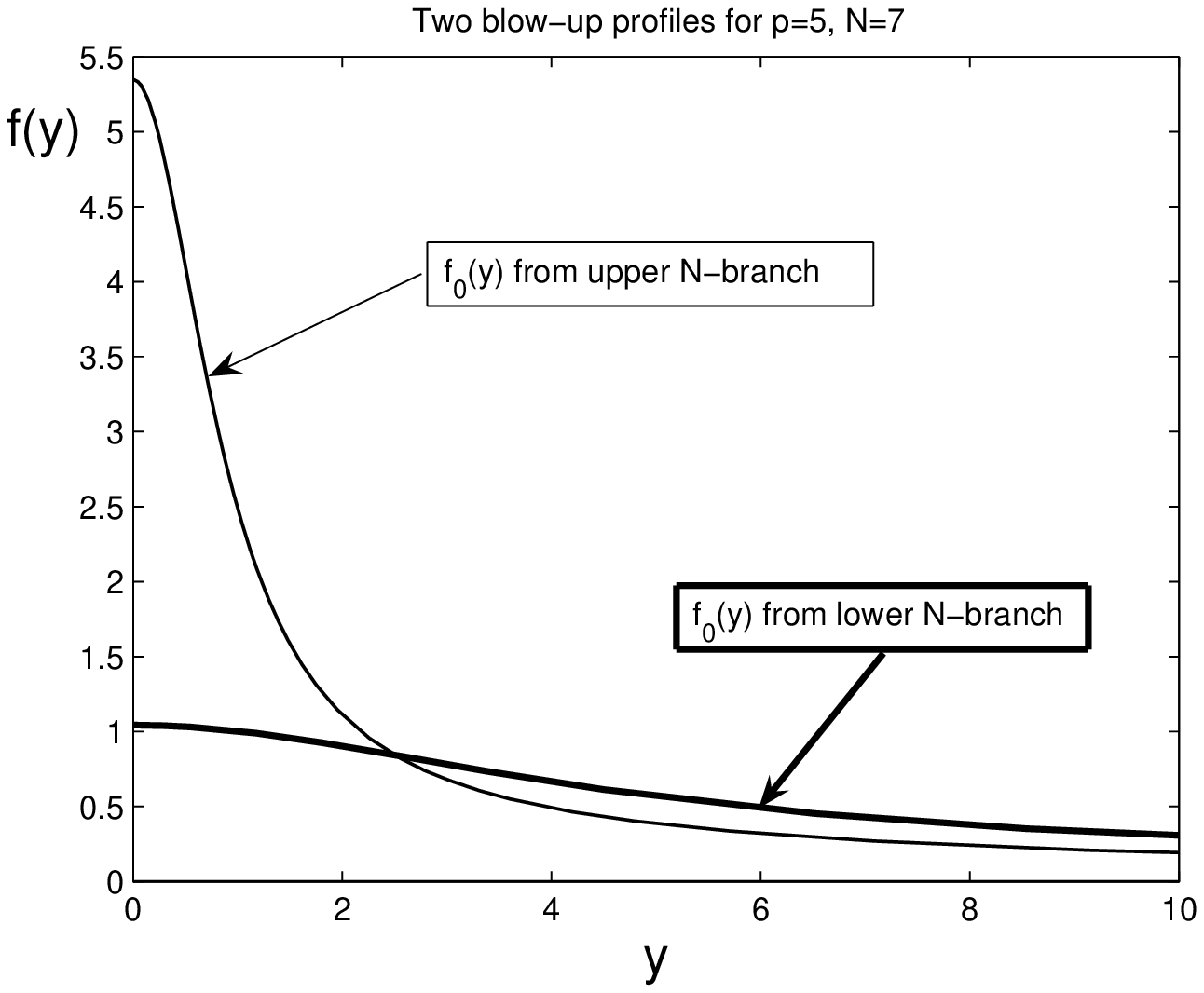}
} \subfigure[$p=2$]{
\includegraphics[scale=0.52]{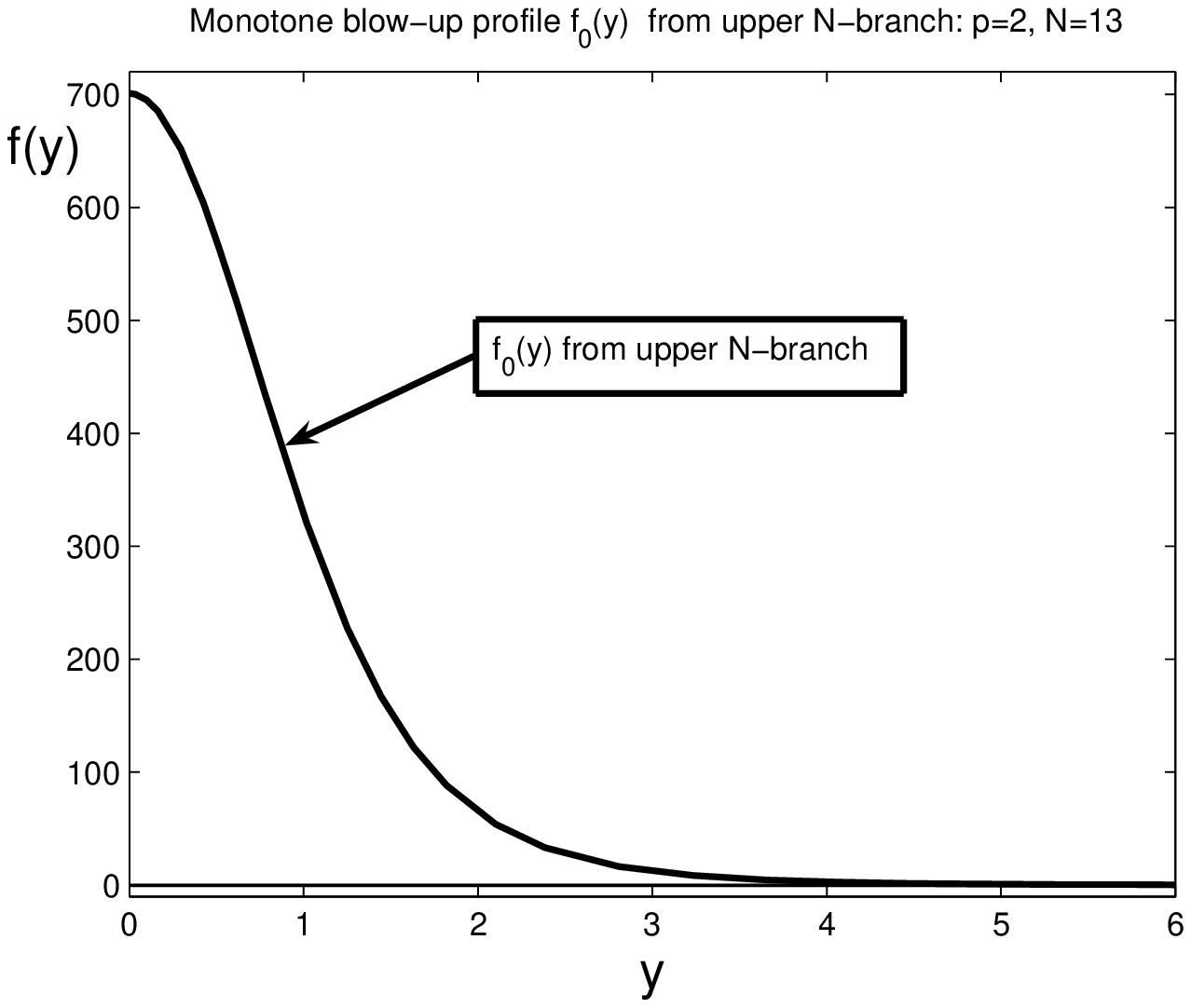}
}
 \vskip -.2cm
\caption{\rm\small Two monotone blow-up profiles of $f_0$-type:
for $p=5$ (a) and for $p=2$ (b).}
 %% \vskip -.3cm
 \label{F21}
\end{figure}

 However, we found that there are other solutions of the monotone type
$f_0$, which are shown in Figure \ref{F21} for $p=5$ (a) and $p=2$
(b), where in the latter one the profile from the lower $N$-branch
is not shown as being too relatively small.

%%FIG%%%%%%%%%%%%%%%%%%%%%%%

%%%%%%%%%%%%%%%%%%%%%%%%%%%%%%%%%%%%%%%%%%%%%%%%%%%%%%%%%%%%%%%%%%%

\begin{figure}
%%\vskip -.3cm
\centering \subfigure[$f_0$-deformation]{
\includegraphics[scale=0.52]{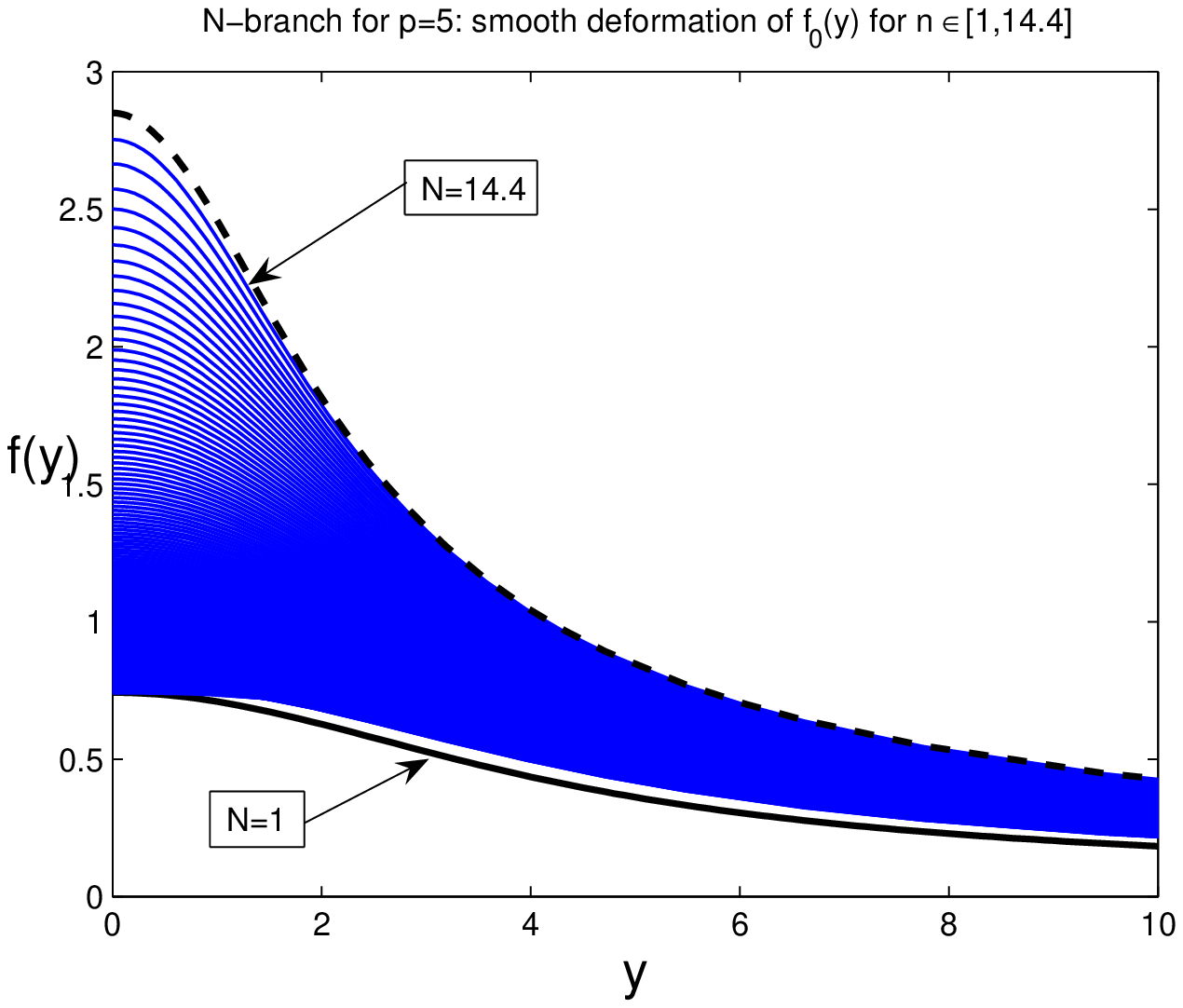}
} \subfigure[lower $N$-branch]{
\includegraphics[scale=0.52]{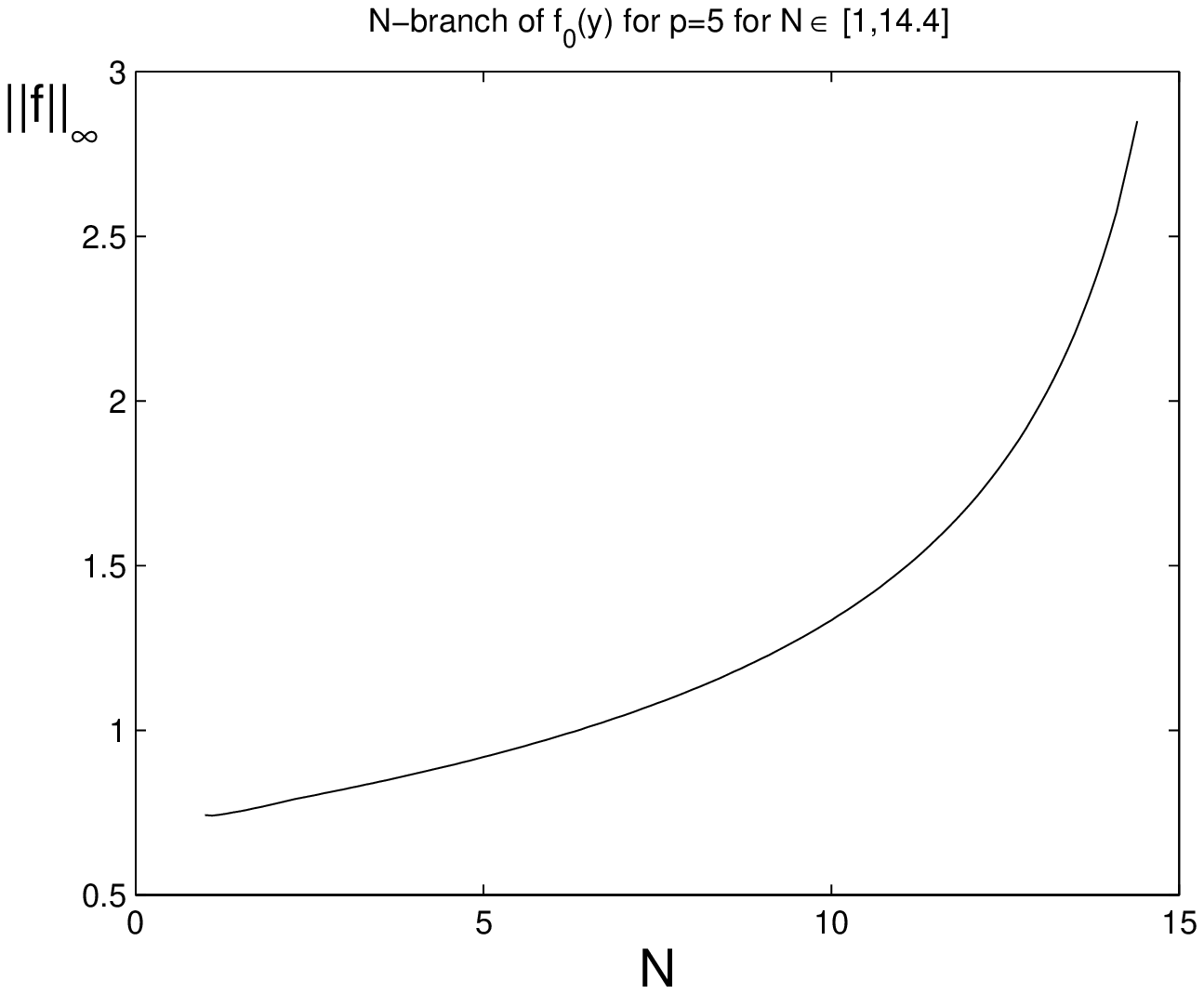}
}
 \vskip -.2cm
\caption{\rm\small Lower $N$-branch of $f_0$ for $p=5$.}
%%% for $p=5$ (a) and for $p=2$
%%(b).}
 %% \vskip -.3cm
 \label{F6}
\end{figure}

\begin{figure}
%%\vskip -.3cm
\centering \subfigure[$f_0$-deformation]{
\includegraphics[scale=0.52]{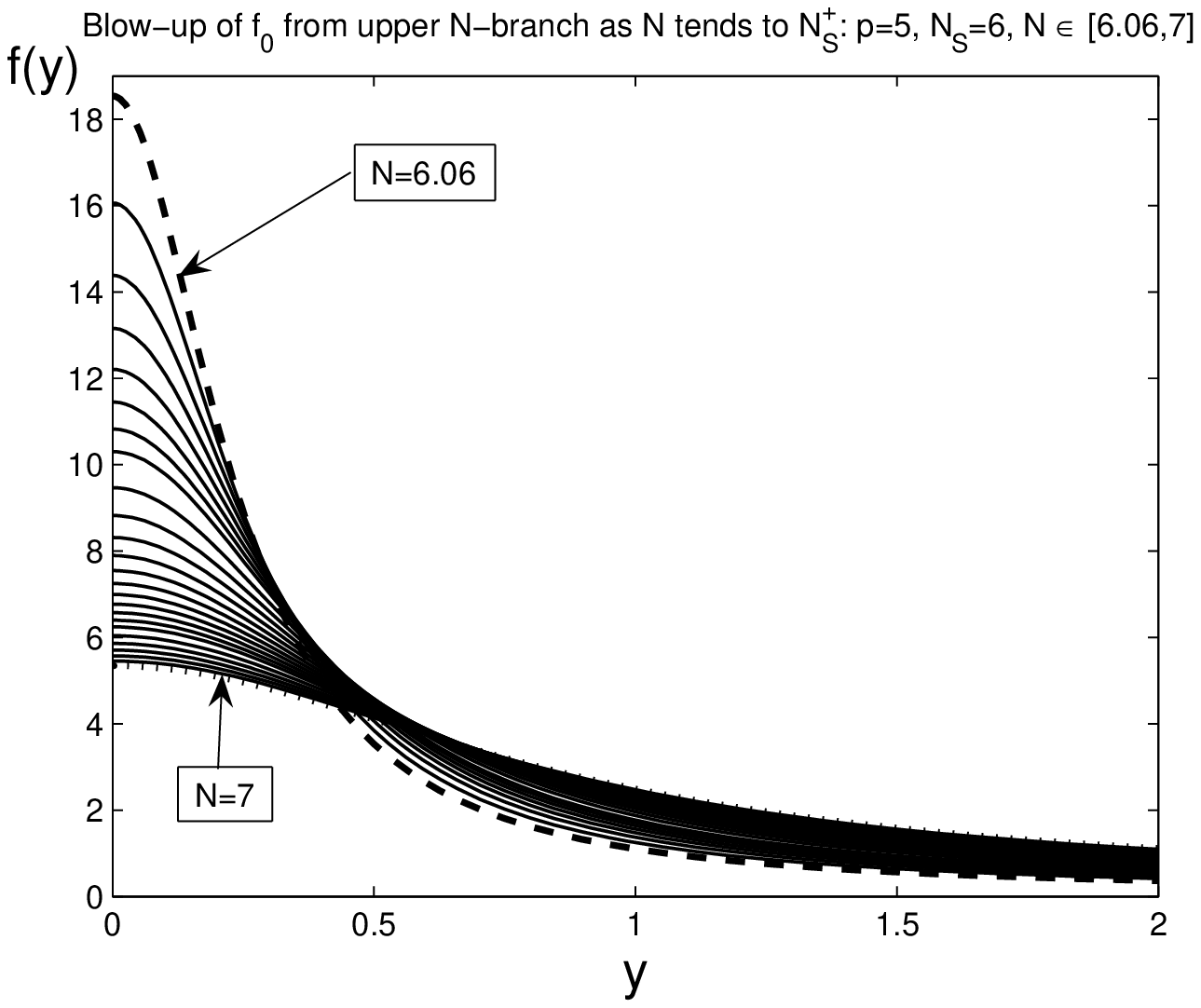}
} \subfigure[upper $N$-branch]{
\includegraphics[scale=0.52]{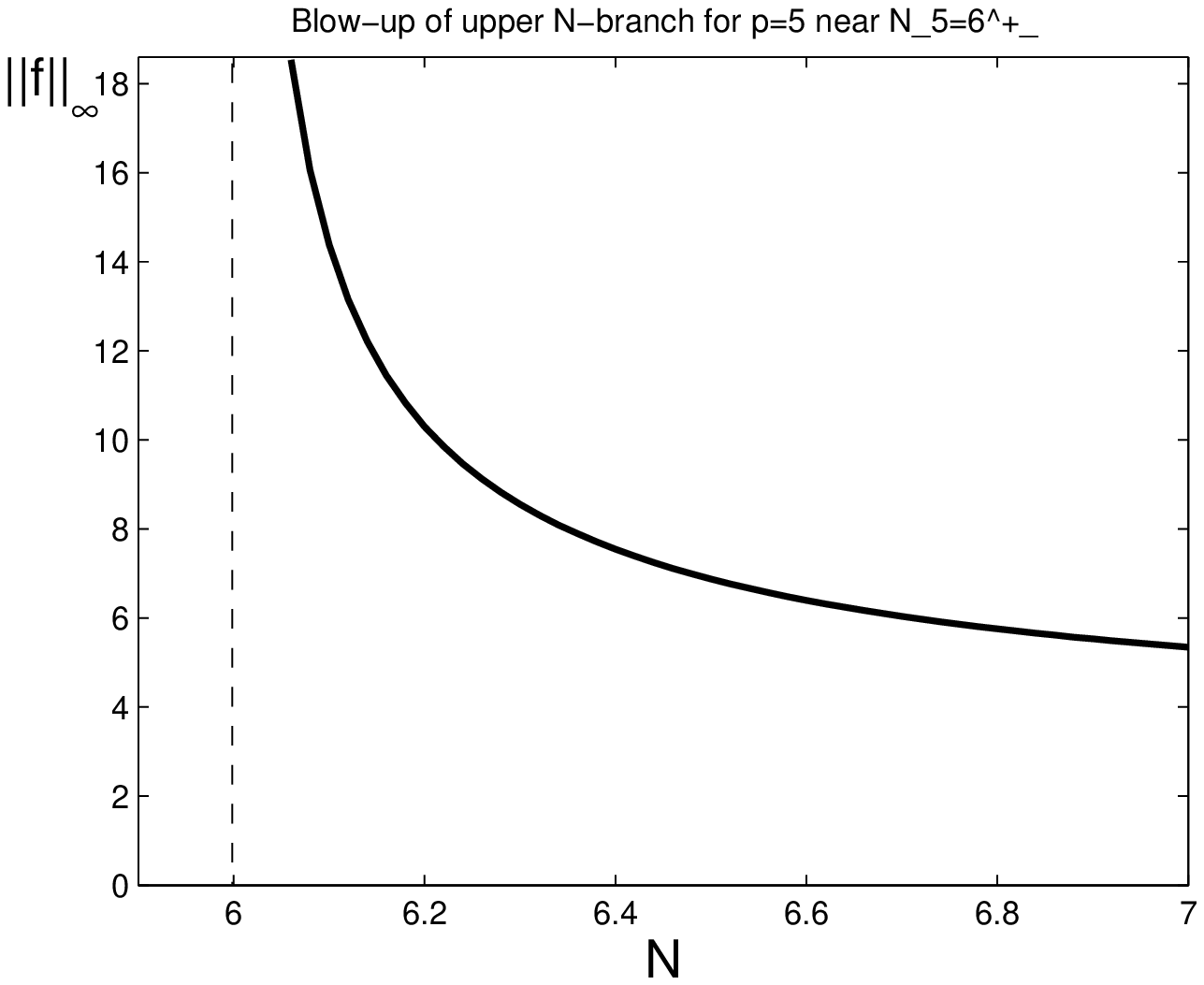}
}
 \vskip -.2cm
\caption{\rm\small Upper $N$-branch of $f_0$ for $p=5$ with
blow-up as $N \to 6^+$.}
 %% \vskip -.3cm
 \label{F7}
\end{figure}

\begin{figure}
 \centering
 \psfrag{||f||}{$\|f\|_\infty$}
 \psfrag{p}{$p$}
 \psfrag{ff}{$\|f\|_{\iy}$}
 \psfrag{NN}{$N$}
 \psfrag{Nsn}{$N_{\rm sn}$}
  \psfrag{v(x,t-)}{$v(x,T^-)$}
  \psfrag{final-time profile}{final-time profile}
   \psfrag{tapp1}{$t \approx 1^-$}
\psfrag{x}{$x$}
 \psfrag{0<t1<t2<t3<t4}{$0<t_1<t_2<t_3<t_4$}
  \psfrag{0}{$0$}
 \psfrag{l}{$l$}
 \psfrag{-l}{$-l$}
\includegraphics[scale=0.36]{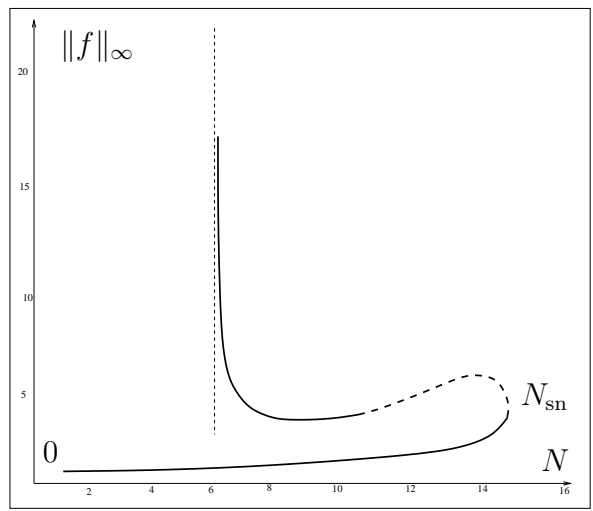}     %%%{F1EE1.pdf}
\caption{\small The whole $N$-branch of $f_0$ for
$p=5$.}
  %%   \vskip -.3cm
 \label{F8}
\end{figure}

Thus, secondly, this nonuniqueness demands another approach to
branching, namely, the $N$-branching that we perform next. In
Figure \ref{F6}, we  show the lower $N$-branch of solutions
$f_0(y)$ for $p=5$, where (a) describes the deformation of $f_0$
and (b) gives the actual $p$-branch.  Blow-up of the upper
$N$-branch as
 \be
 \label{NNpp}
  \tex{
   N \to N_5^+ \whereA N_p: \,\,\, \frac{N+4}{N-4}=p,
   }
   \ee
   so that $N_5=6$ for $p=5$, is shown in Figure \ref{F7}, with
   the same meaning of (a) and (b).
   A general view of the whole $N$-branch of $f_0$ for $p=5$ is schematically  explained
   in Figure \ref{F8}, where by dotted line we draw a possible
   expected but still
   hypothetical connection of the lower (stable) and the upper
   (more unstable, plausibly) $f_0$-branches, which we were not
   able to reconstruct  numerically. Numerical continuation in the
   parameter $N$ is quite a challenging problem in some $N$-ranges.
   %% especially for
   %%the author, who is not professional at all in this area.

%%\begin{figure}
%  \vskip -.3cm
%%\centering
%%\includegraphics[scale=0.75]{Nbrp5.eps} %%%%{AmFBP1.eps}  %%%%%%%{4het5.eps}   %%%%%{4het.eps}  Old
%\includegraphics{F1.1}
%%\vskip -.3cm
%%  \caption{\small The whole $N$-branch of $f_0$ for
%%$p=5$.}
%%   \vskip -.3cm
%% \label{F8}
%%\end{figure}
%%%%%%%%%%%%%%%%%%%%%%%%%%%%%%%%%%%%%%%%%%%%%%%%%%%%%%%%%%%%%%

 Thus we expect that there exists a saddle-node bifurcation at
 some
  \be
  \label{sn1}
  p=5: \quad N_{\rm sn} \in (14.979,15).
   \ee

%%FIG%%%%%%%%%%%%%%%%%%%%%%%

%%%%%%%%%%%%%%%%%%%%%%%%%%%%%%%%%%%%%%%%%%%%%%%%%%%%%%%%%%%%%%%%%%%

%%FIG%%%%%%%%%%%%%%%%%%%%%%%

%%%%%%%%%%%%%%%%%%%%%%%%%%%%%%%%%%%%%%%%%%%%%%%%%%%%%%%%%%%%%%%%%%%

%%%%%%%%%%%%%%%%%%%%%%%%%%%%%%%%%%%%%%%%%%%%%%%%%%%%%%%%%%%%
%%%%%%%%%%%%%%%%%%%%%%%%%%%%%%%%%%%%%%%%%%%%%%%%%%%%%%%%%
%%%%%%%%%%%%%%%%%%%%%%%%%%%%%%%%%%%%%%%%%%%%%%%%%%%%%%%%%%%%%%

\begin{figure}
%%\vskip -.3cm
\centering \subfigure[lower $N$-branch]{
\includegraphics[scale=0.52]{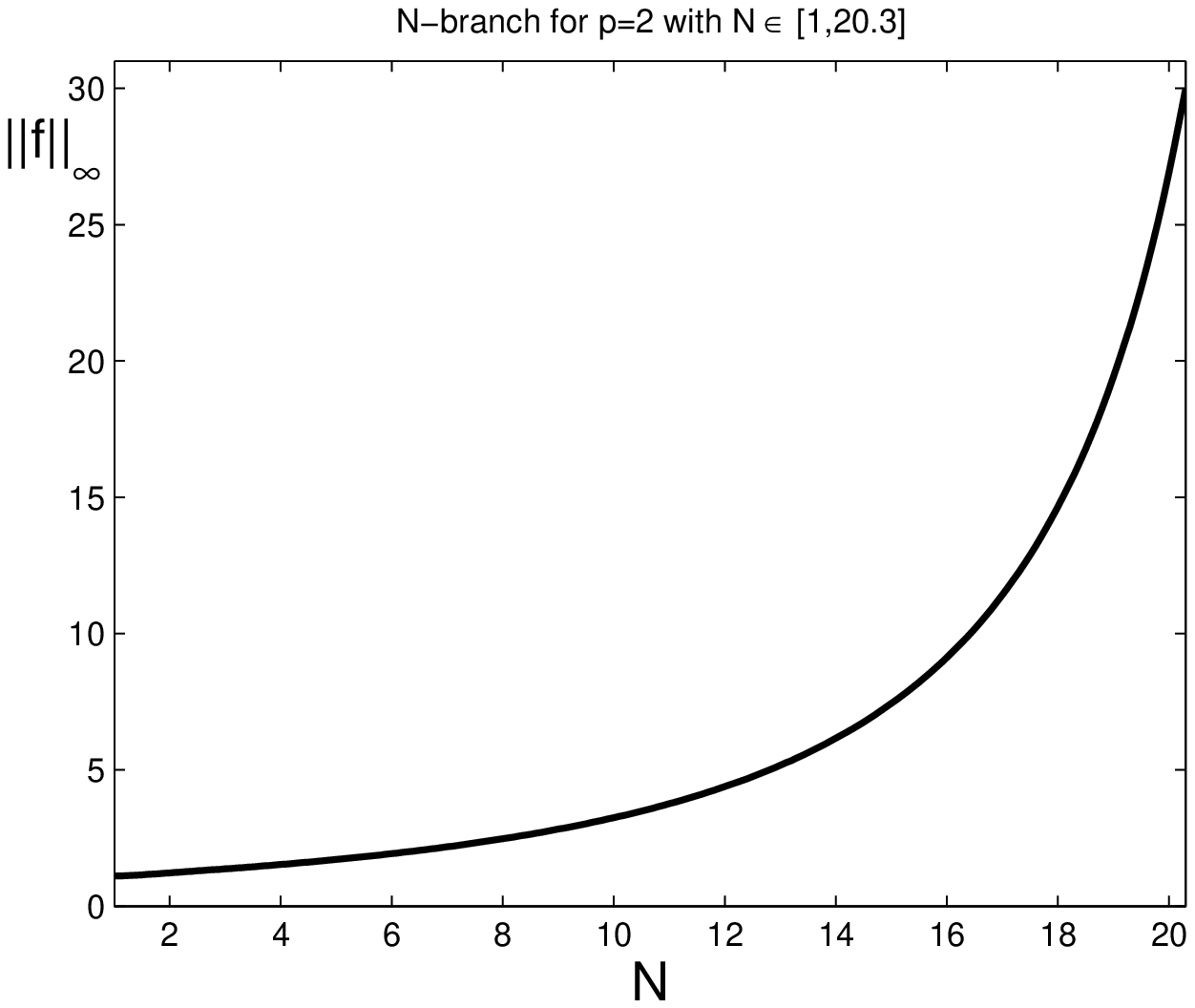}          %%%{p2N132.eps}
} \subfigure[upper $N$-branch]{
\includegraphics[scale=0.52]{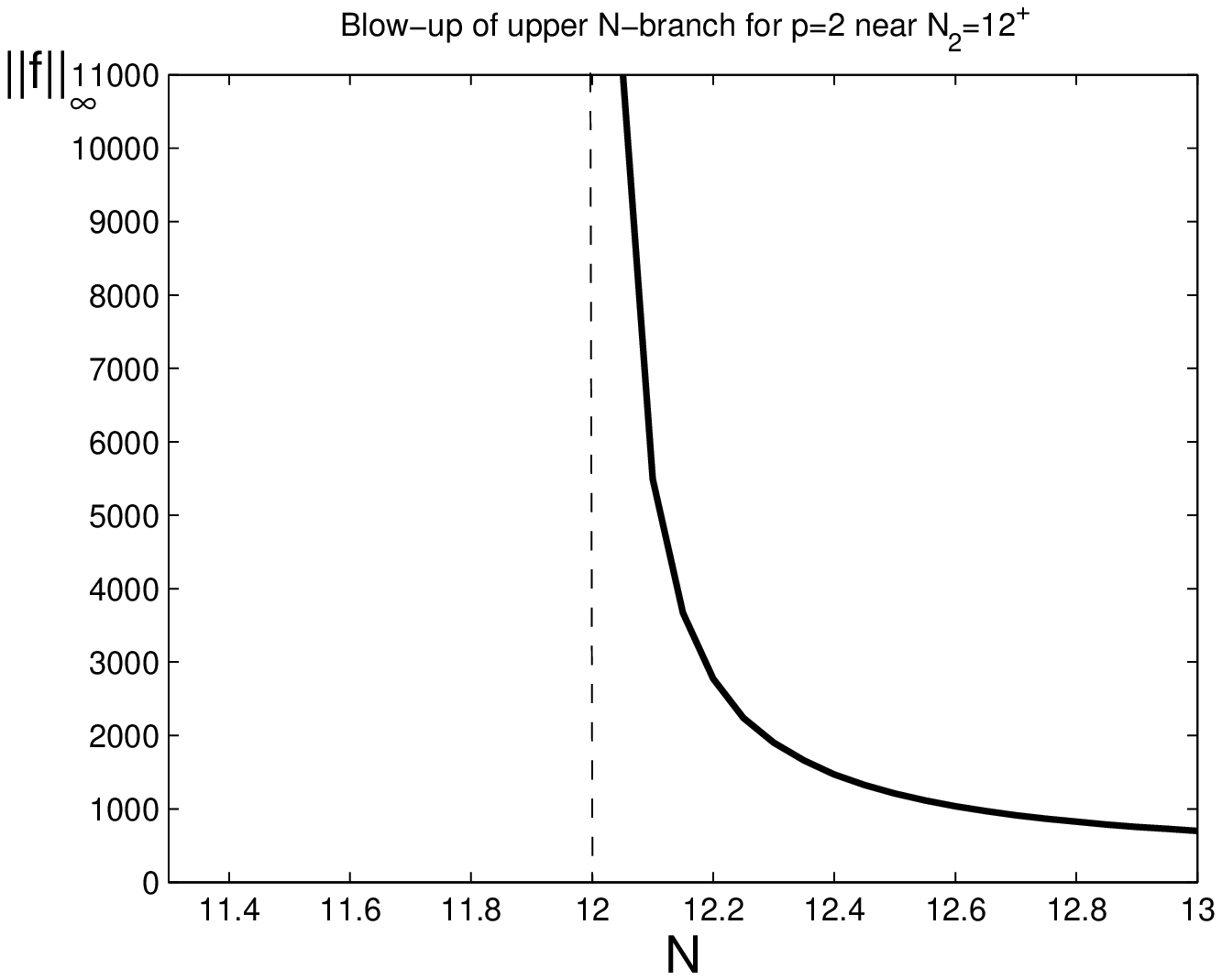}
}
 \vskip -.2cm
\caption{\rm\small Lower and upper $N$-branches of $f_0$ for
$p=2$: the lower branch for $N \in [1,20.3]$ (a) and  blow-up of
the upper one as $N \to 12^+$ (b).}
 %% \vskip -.3cm
 \label{F9}
\end{figure}

In Figure \ref{F9} for $p=2$, we show blow-up of the upper
$N$-branch as $N \to N_2=12^+$. We then expect that the lower and
upper branches have a turning (saddle-node) bifurcation point at
some
 $$
 p=2: \quad N_{\rm sn} \in (20.3,23].
 $$
 We
   hope that such an interesting saddle-node branching phenomenon
   will attract true experts in numerical methods, bearing in mind
   that numerical experiments might be for a long time the only
   tool of the study of such blow-up phenomena.

%%FIG%%%%%%%%%%%%%%%%%%%%%%%

%%%%%%%%%%%%%%%%%%%%%%%%%%%%%%%%%%%%%%%%%%%%%%%%%%%%%%%%%%%%%%%%%%%

\begin{figure}
%%\vskip -.3cm
\centering \subfigure[$f_0$-deformation]{
\includegraphics[scale=0.52]{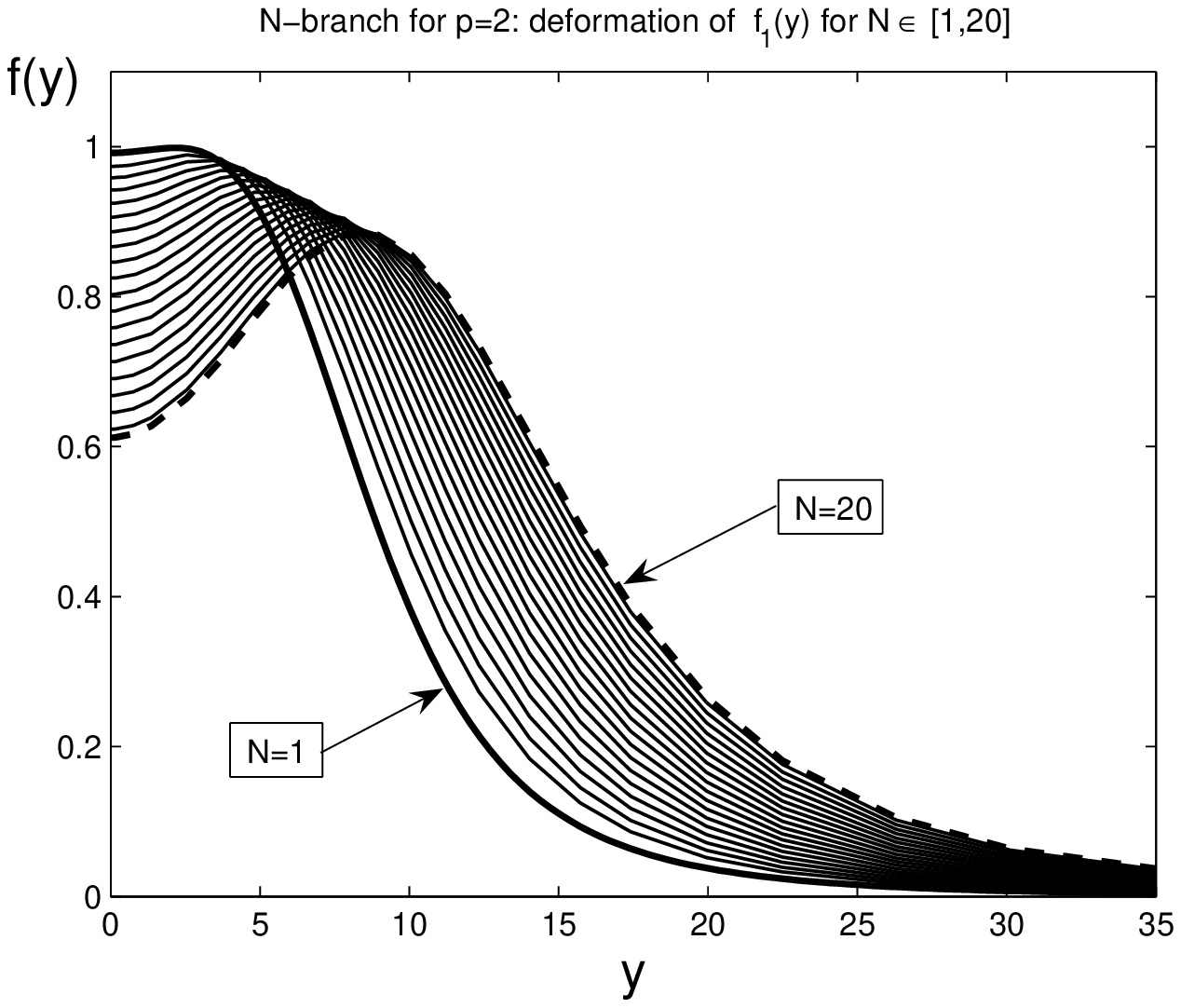}
} \subfigure[$N$-branch]{
\includegraphics[scale=0.52]{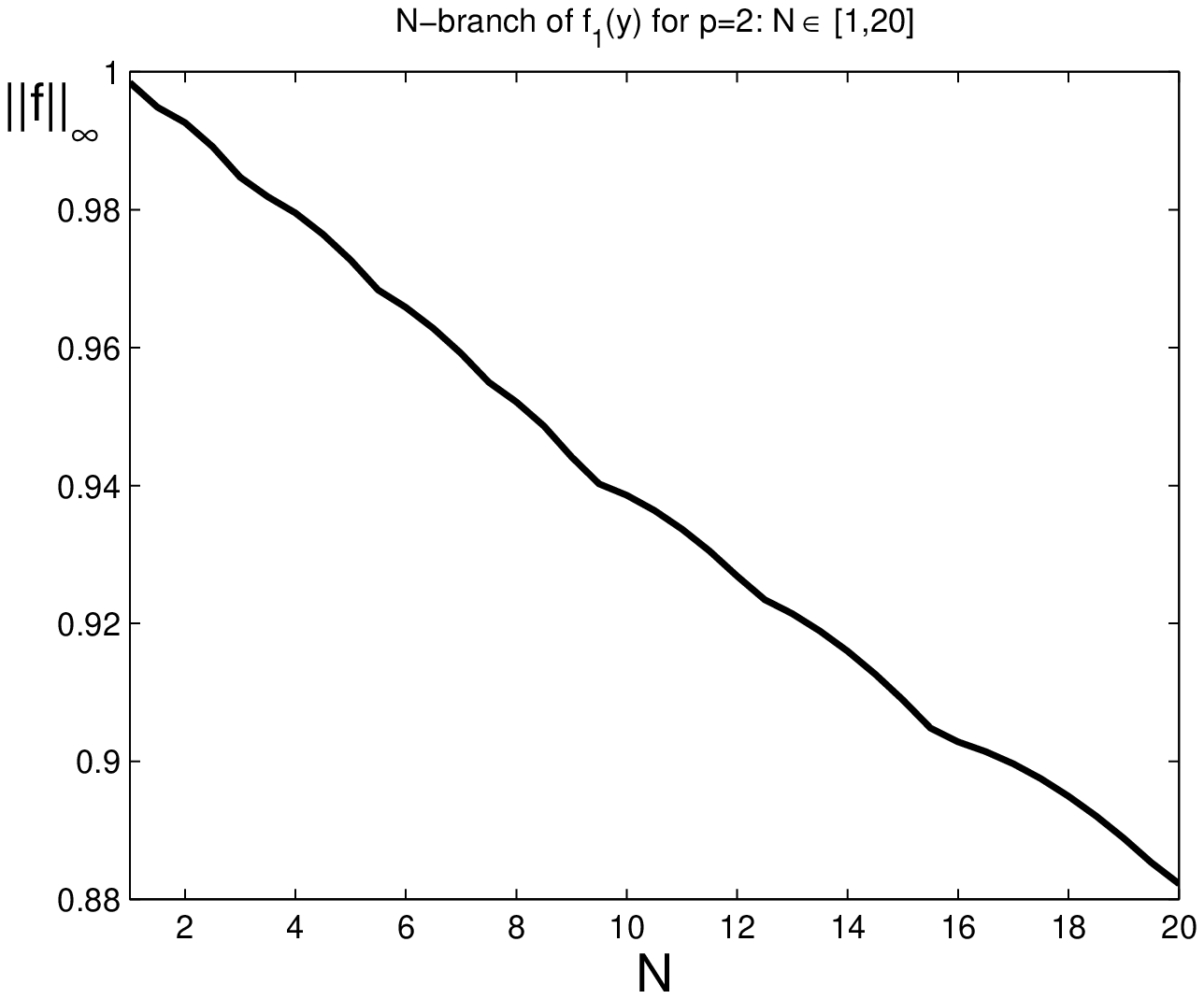}
}
 \vskip -.2cm
\caption{\rm\small The $N$-branch of $f_1$ for $p=2$ and $N \in
[1,20]$.}
%%% for $p=5$ (a) and for $p=2$
%%(b).}
 %% \vskip -.3cm
 \label{Ff1N}
\end{figure}

Finally, in Figure \ref{Ff1N}, we present the numerical results
confirming $N$-branching for $p=2$ of the second blow-up profile
$f_1$, where, as usual, (a) describes smooth deformation of
$f_1(y)$, while (b) shows  the $N$-branch. It seems that
$N$-branches of $f_1$ are global and do not suffer from a
saddle-node bifurcations.

%%FIG%%%%%%%%%%%%%%%%%%%%%%%
%%%%%%%%%%%%%%%%%%%%%%%%%%%%%%%%%%%%%%%%%%%%%%%%%%%%%%%%%%%%%%%%%%%

%%%%%%%%%%%%%%%%%%%%%%%%%%%%%%%%%%%%%%%%%%%%%%%%%%%%%%%%%%%%%%%%%%%%%%%%%
\subsection{On sign changes of $f_0(y)$ and $f_1(y)$}
%% final time measure-like Type I(ss) blow-up}

We now study some particularly important properties of blow-up
similarity profiles. We begin with the easier property of sign
changes. We have seen already several strictly positive profiles
$f_0(y)$ for some $p$'s, which is rather surprising since the
equations do not obey the Maximum Principle. However, we will show
that, for smaller $p$, the similarity profiles can gain extra
zeros as sign changes. Such $p$'s , when a zero is gained, we
denote by $p_0(N)$.

\begin{figure}
%%\vskip -.3cm
\centering \subfigure[$f_0$ profiles, $N=1$]{
\includegraphics[scale=0.52]{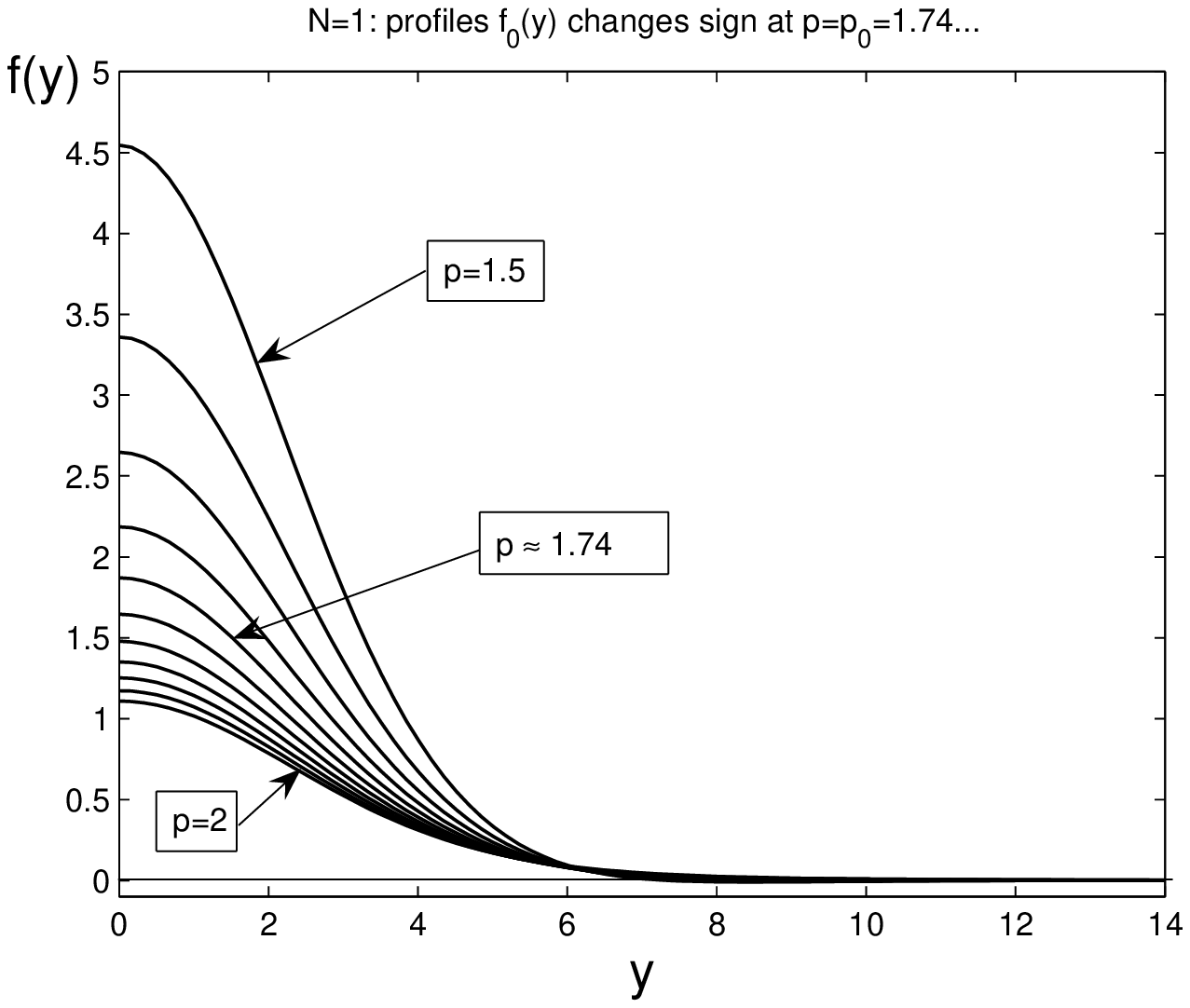}
} \subfigure[shooting $p_0(1)$, $N=1$]{
\includegraphics[scale=0.52]{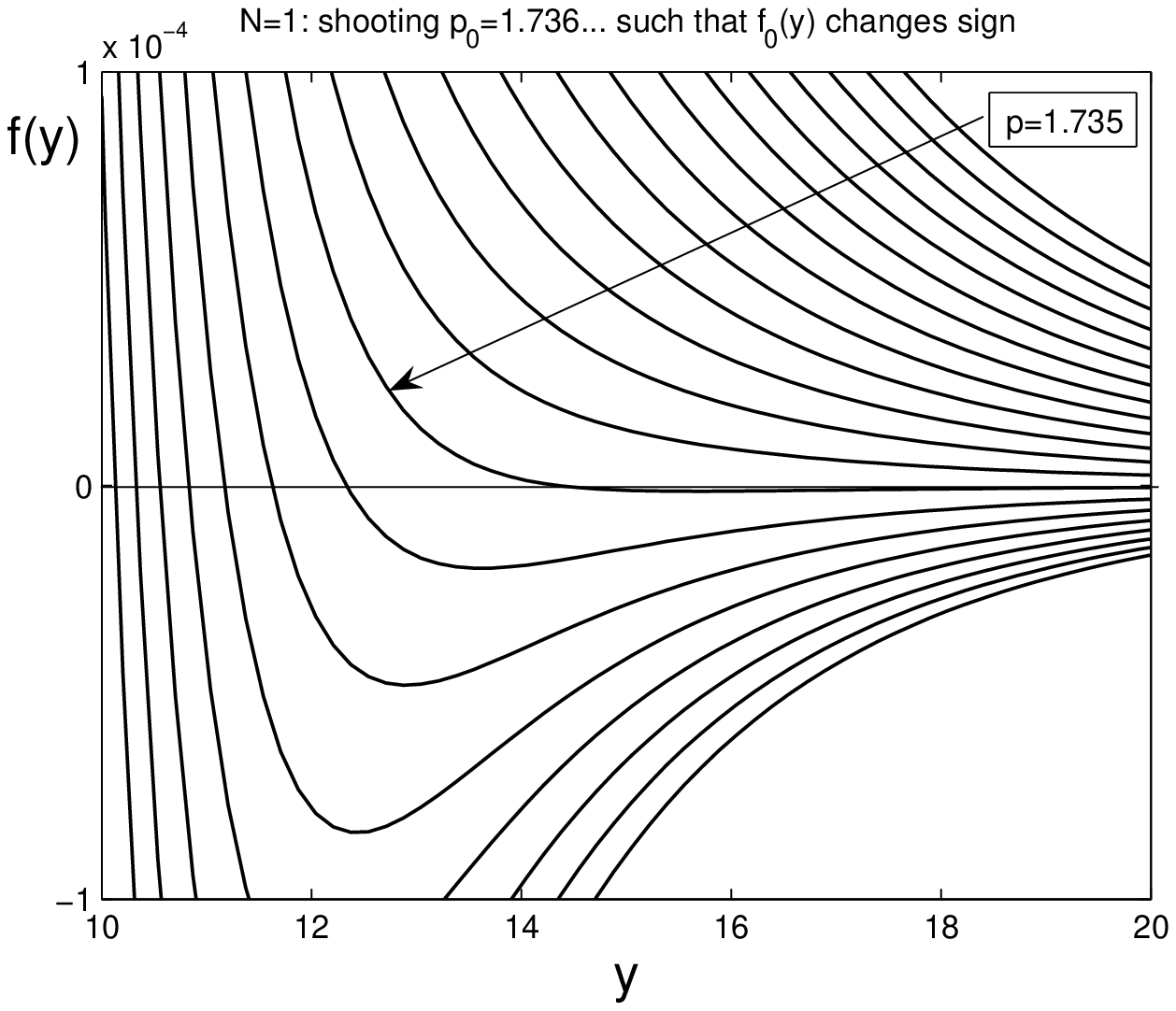}
}
 \vskip -.2cm
\caption{\rm\small Shooting the root $p_0(1)$  for $N=1$.}
%%% for $p=5$ (a) and for $p=2$
%%(b).}
 %% \vskip -.3cm
 \label{FD1}
\end{figure}

\begin{figure}
%%\vskip -.3cm
\centering \subfigure[$f_1$ profiles, $N=1$]{
\includegraphics[scale=0.52]{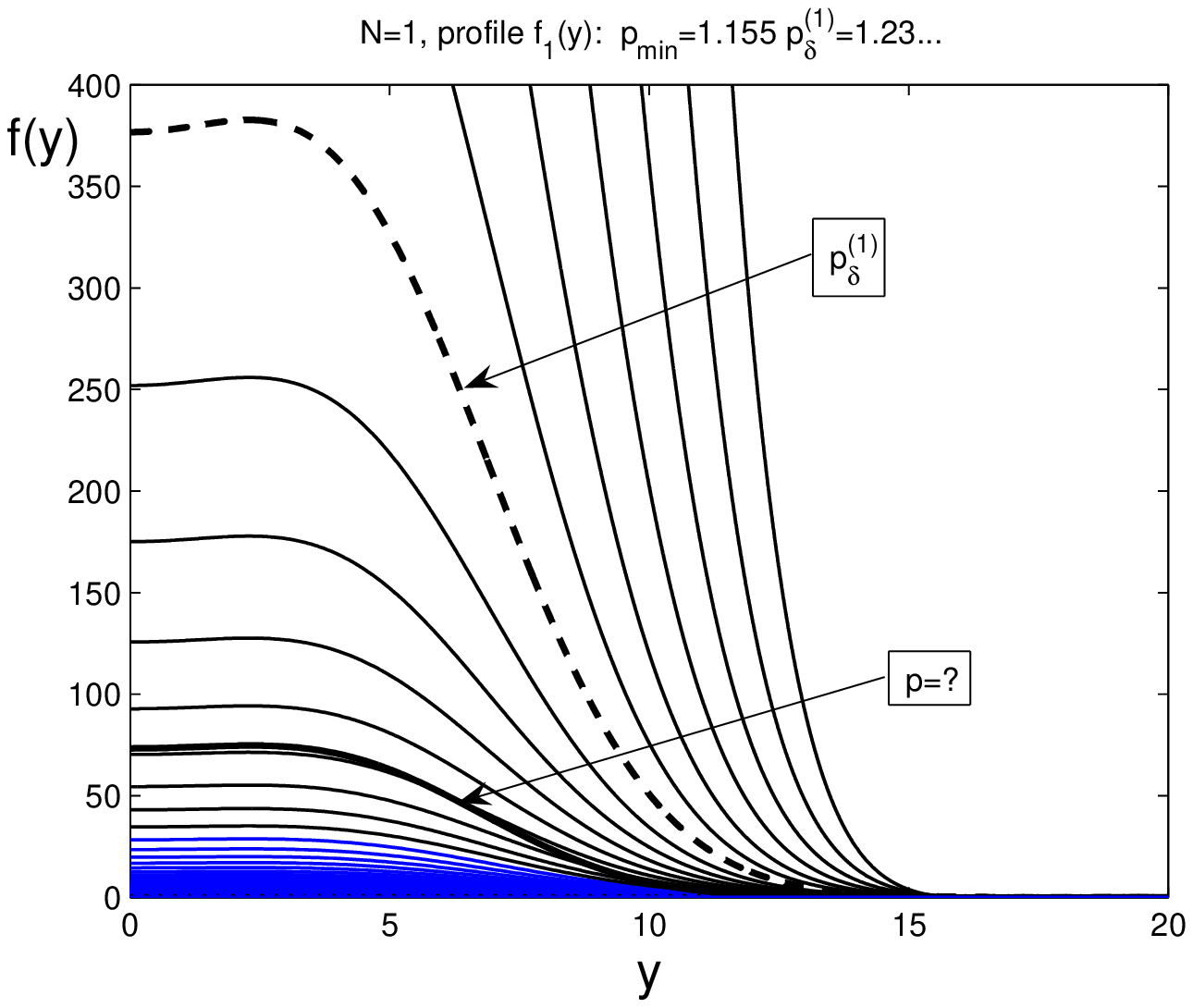}
} \subfigure[shooting $p_0^{(1)}(1)$, $N=1$]{
\includegraphics[scale=0.52]{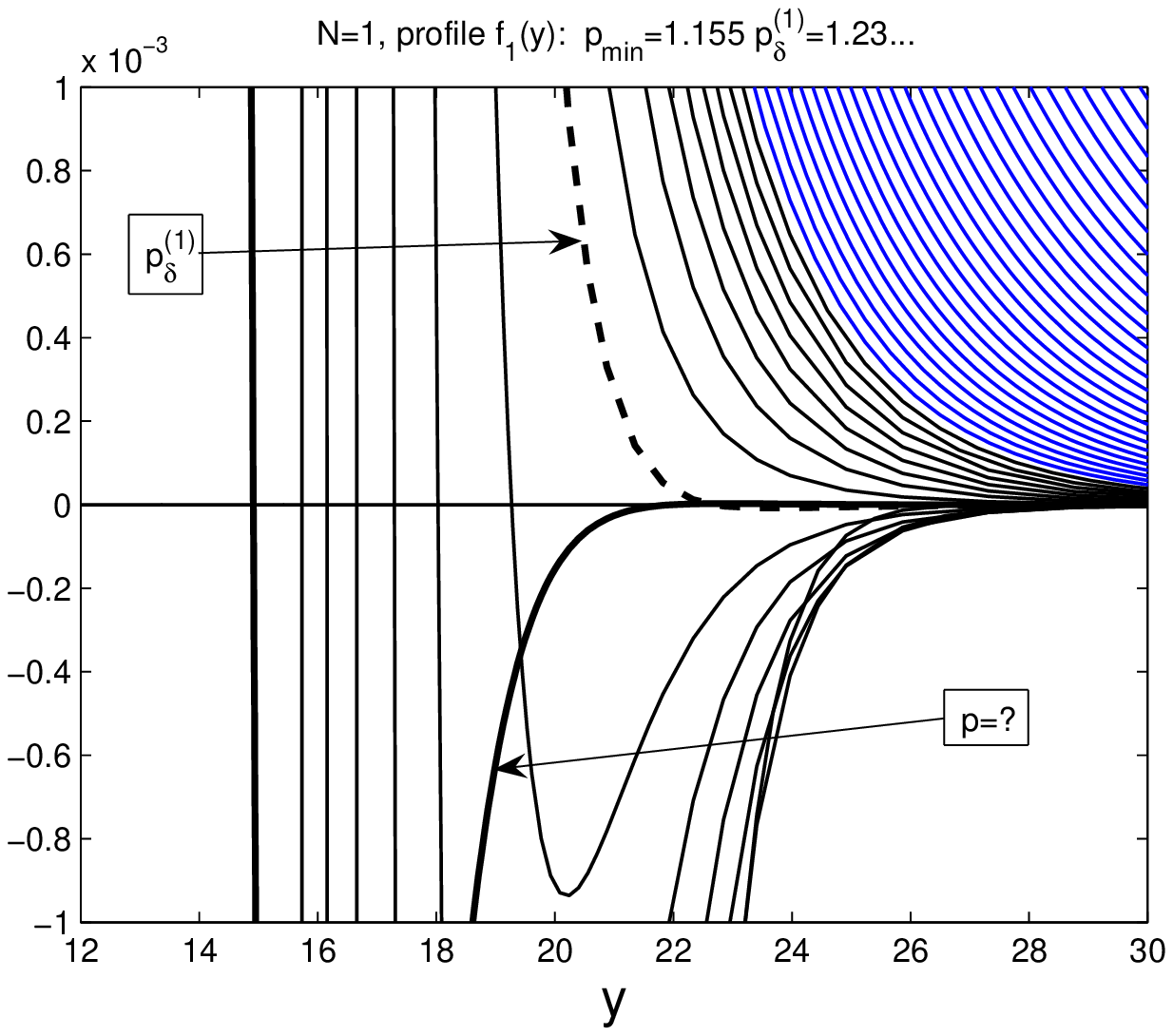}
}
 \vskip -.2cm
\caption{\rm\small Shooting the root $p_0^{(1)}(1)$  for $N=1$.}
%%% for $p=5$ (a) and for $p=2$
%%(b).}
 %% \vskip -.3cm
 \label{FD11}
\end{figure}

 %%We now consider the particularly difficult problem \ef{dd4} and
%%begin with
Consider  one dimension $N=1$. Firstly, the attentive Reader can
see that in Figure \ref{F1}(a), already for $p= \frac 32$ the
profile $f_0(y)$  changes sign, while for $p=2$, it is positive.
Hence,
 $p_0(1) \in (\frac 32,2)$. Secondly, more thorough numerics are presented
 %%  The numerical
%%results
 in Figure \ref{FD1}, where (a) shows $f_0$'s in a vicinity of
  \be
  \label{pp1ss}
  p_0(1)=1.7358...\, ,
   \ee
while (b) shows a sharp shooting of the critical value \ef{pp1ss}.
Figure \ref{FD11} shows shooting%% the root
 $$
 p_0^{(1)}(1)=1.23...\, ,
  $$
  at which the {\em second}
profile $f_1(y)$ gets a new zero. By the boldface line we denote a
new ``$p$-undetected" solution with  extra zeros gained at another
$p_{01}$, showing that such roots are not unique.
%%% seems also satisfying \ef{dd4}.
 Since $f_1(y)$ is expected to be less stable, we will
concentrate on the roots $p_0$ for the generic blow-up profile
$f_0(y)$.

\begin{figure}
%%\vskip -.3cm
\centering \subfigure[shooting $p_0(3)$]{
\includegraphics[scale=0.52]{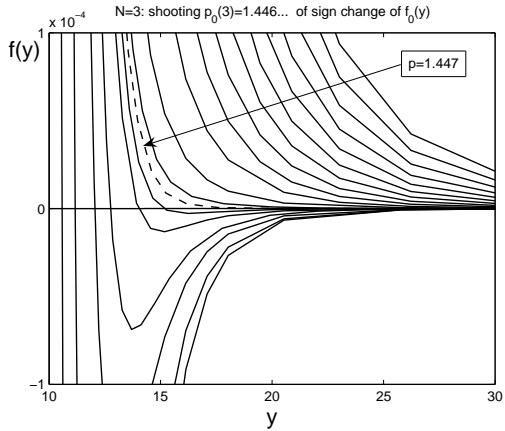} } \subfigure[$N=9$]{
\includegraphics[scale=0.52]{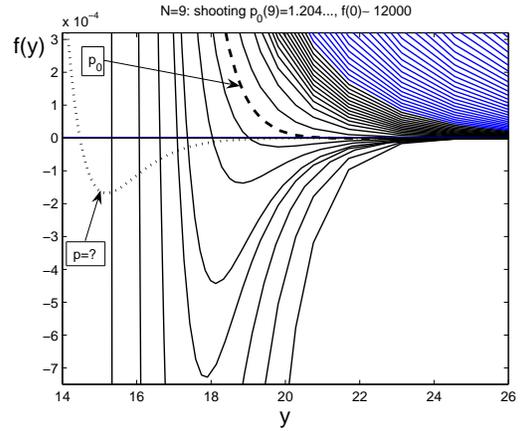} %%%%%%%%{ff1.eps} , N=8             %%{delN8R2.eps}
}
 \vskip -.2cm
\caption{\rm\small Shooting the root $p_0(1)$: for $N=3$ (a) and
for $N=9$ (b).}
%%% for $p=5$ (a) and for $p=2$
%%(b).}
 %% \vskip -.3cm
 \label{FD2}
\end{figure}

\begin{figure}
 %%\vskip -.2cm
   \centering \subfigure[shooting $p_0$ for $N=10$]{
\includegraphics[scale=0.52]{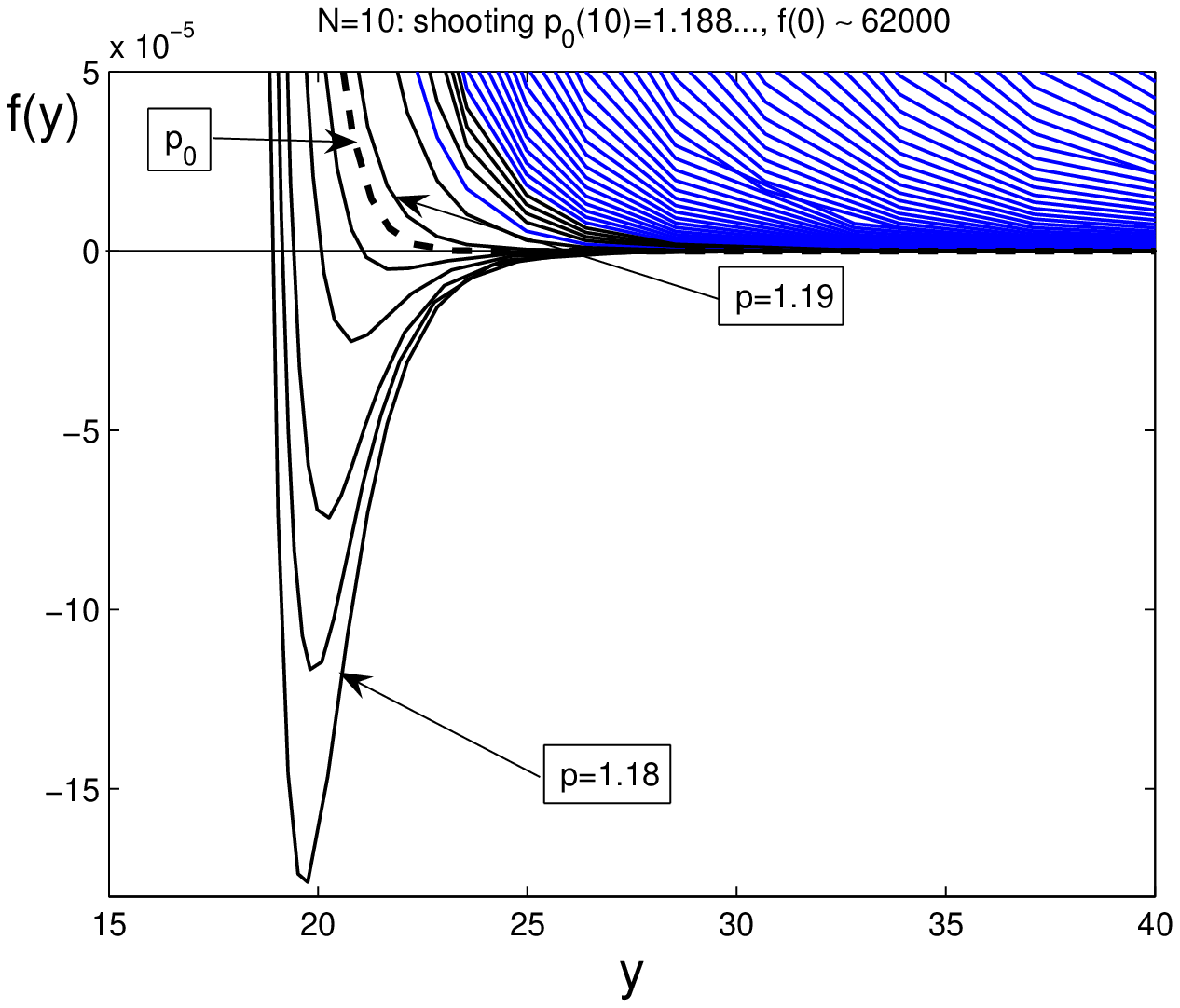}
} \subfigure[$N=12$]{
\includegraphics[scale=0.52]{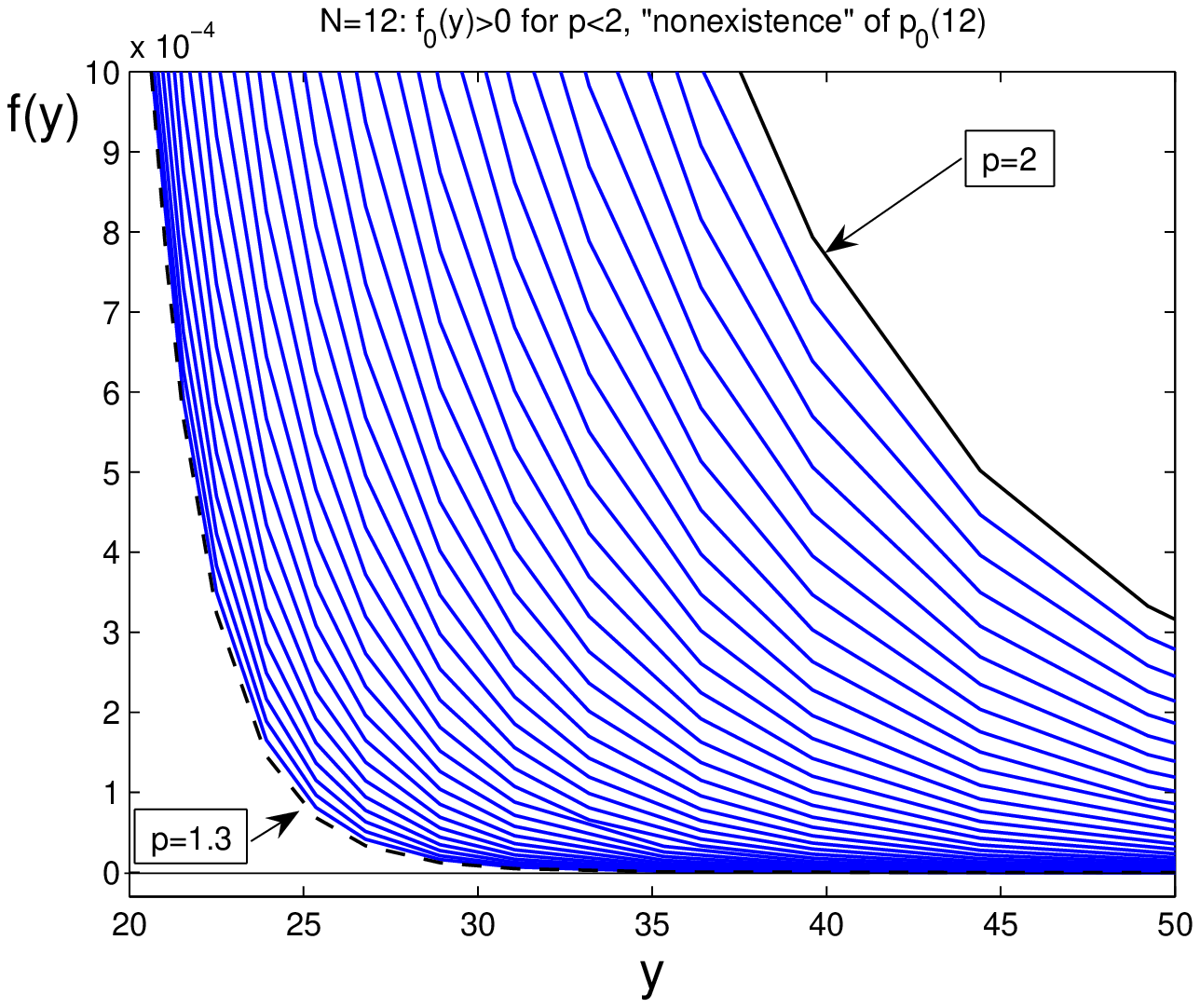}
}
 \vskip -.2cm
\caption{\rm\small Shooting the root $p_0(1)$: for $N=10$ (a) and
no sign changes of $f_0(y)$ for $N=12$ for  $p \ge 1.3$ (b).}
%%% for $p=5$ (a) and for $p=2$
%%(b).}
 %% \vskip .2cm
 \label{FD3}
\end{figure}

Thus, similarly, Figure \ref{FD2}(a) yields
%% , we show a similar shooting of the
%%root of \ef{dd4}
\be
  \label{pp2}
  p_0(3)=1.446... \andA p_0(9)=1.204...\,.
  %%%\quad(p_\d(8)=1.226...). %%%%% 1.226... for N=8
   \ee
    In Figure \ref{FD3}(a), a similar phenomenon is checked  for $N=10$, with
    $p_0(10)=1.188...\,.$
    In (b), we see no sign changes of $f_0(y)$ for $N=12$, but
    this happens for smaller $p \approx 1^+$, when $f(0)$ gets $10^{5}$--$10^{6}$,
     while their negative counterparts take
values $\sim -10^{5}$,
    and numerics become rather unreliable.
  The  overall numerical results for shooting $p_0(N)$ are shown in Table 1.

 %%\vskip -.2cm
%%%%%%%%%%%%%%%%%%%%%%%%%%%%%%%%%%%%%%%%%%%%%%%%%%%%%%%%%%%%%%%%%%%%
\begin{table}[h]
\caption{Roots $p_0(N)$  for profiles $f_0(y)$ to get a new zero}
\begin{tabular}{@{}ll}
 $N$ & $p_0(N)$
 \\\hline
  $1$ & 1.7358...
  \\
   $2$ & 1.53...
 \\
    $3$ & 1.446...
   \\
    $4$ & 1.377...
   \\
    $5$ & 1.320...
   \\
    $6$ & 1.28...
   \\
    $7$ & 1.25...
   \\
    $8$ & 1.226...
   \\
    $9$ & 1.204...
   \\
    $10$ & 1.188...
   \\
    $11$ & 1.16...
   \\
\end{tabular}
\end{table}
%%\vskip -.2cm
%%%%%%%%%%%%%%%%%%%%%%%%%%%%%%%%%%%%%%%%%%%%%%%%%%%%%%%%%%%%%%%%%%%%%%%%%

%%More precisely, we observed $f_0(y)$ of changing sign for $N \le
%%8$ only, so we expect that
%% \be
%% \label{pp3}
%% \mbox{for $N \ge 9$, the algebraic equation (\ref{dd4}) does not
%% have a solution $p_\d$.}
%%  \ee
%%It was difficult to check this numerically, since suspicious
%%profiles $f_0(y)$ for $p \approx 1^+$ can take huge values
%%%$\|f\|_\iy \sim 10^{4}$,
 A proper asymptotic theory for $p \approx
1^+$ involving expansions such as \ef{dd1} and $(p-1)^{-\frac
1{p-1}}$-scaling of the ODE \ef{2.4} (cf. \cite[\S~5]{GHUni})
would be fruitful. Note that we have observed some numerical
evidence for existence of the second root $p_{0 1}$ for $N=8$ and
$9$ (see the dotted line in Figure \ref{FD2}), which is surprising
in view of non-oscillating of the exponential term in \ef{dd1},
but numerics were too difficult and rather poor to identify the
new root if any.
%%% to support or discharge the claim \ef{pp3}.

%%FIG%%%%%%%%%%%%%%%%%%%%%%%
%%%%%%%%%%%%%%%%%%%%%%%%%%%%%%%%%%%%%%%%%%%%%%%%%%%%%%%%%%%%%%%%%%%

%%FIG%%%%%%%%%%%%%%%%%%%%%%%
%%%%%%%%%%%%%%%%%%%%%%%%%%%%%%%%%%%%%%%%%%%%%%%%%%%%%%%%%%%%%%%%%%%

%%FIG%%%%%%%%%%%%%%%%%%%%%%%
%%%%%%%%%%%%%%%%%%%%%%%%%%%%%%%%%%%%%%%%%%%%%%%%%%%%%%%%%%%%%%%%%%%

%%FIG%%%%%%%%%%%%%%%%%%%%%%%
%%%%%%%%%%%%%%%%%%%%%%%%%%%%%%%%%%%%%%%%%%%%%%%%%%%%%%%%%%%%%%%%%%%

%%%%%%%%%%%%%%%%%%%%%%%%%%%%%%%%%%%%%%%%%%%%%%%%%%
\subsection{On final time measure-like Type I(ss) blow-up}

This is a much more difficult problem, which we resolve
numerically for $N=1$ only. For $N \ge 2$, we got no sufficiently
reliable results (rather plausibly, such a self-similar phenomenon
may be unavailable in some higher dimensions).

\begin{figure}
%  \vskip -.3cm
\centering
\includegraphics[scale=0.75]{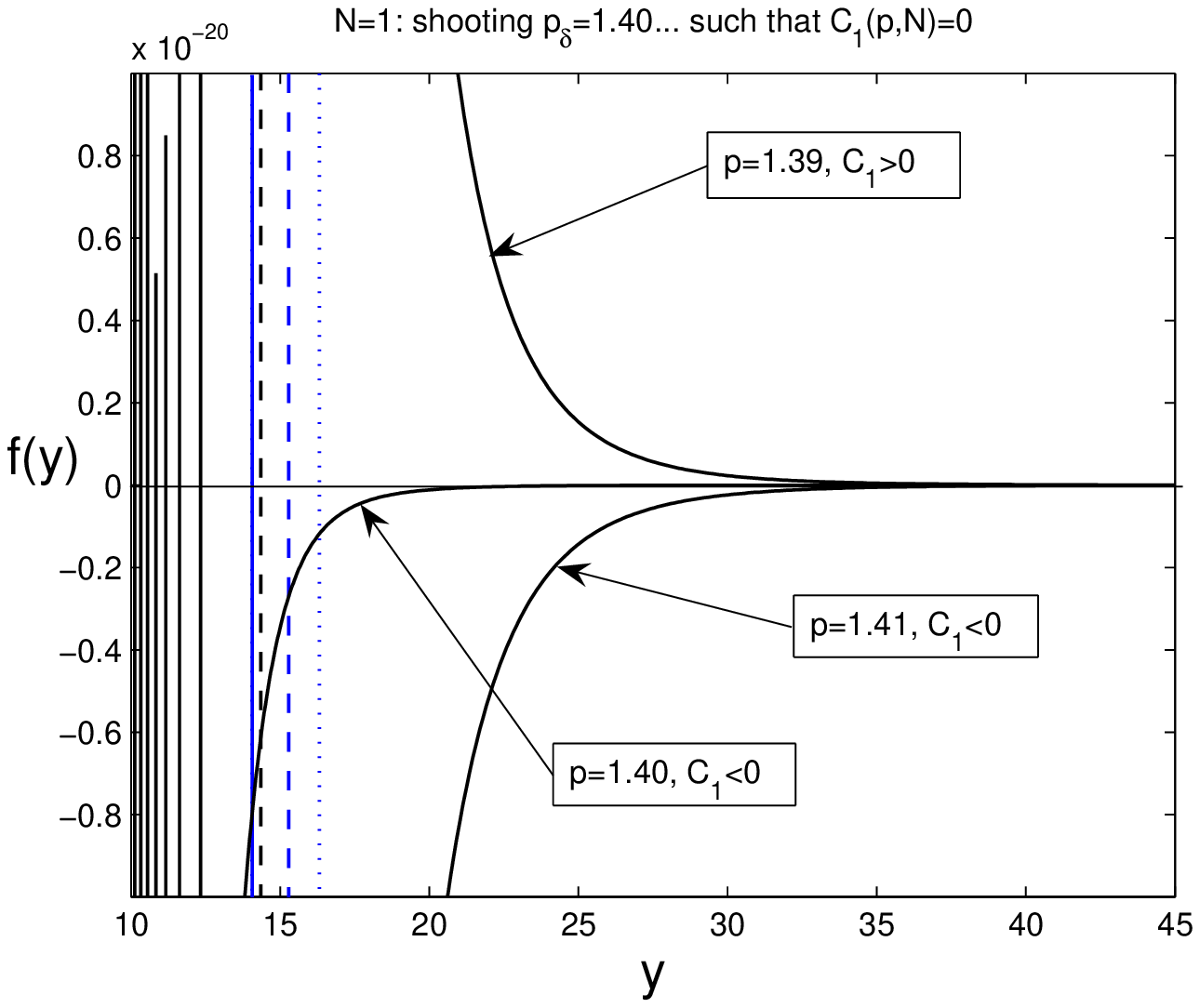} %%%%{AmFBP1.eps}  %%%%%%%{4het5.eps}   %%%%%{4het.eps}  Old
\vskip -.3cm \caption{\small Shooting the root \ef{rr1} of the
algebraic equation \ef{dd4} for $N=1$.}
  %%% \vskip -.3cm
 \label{FF1gg}
\end{figure}

We will refine Figure \ref{FD1}. We claim that the equation
\ef{dd4} has the following root:
 \be
 \label{rr1}
 p_\d(1)=1.40... < p_0(1)=1.7358... \, ,
  \ee
  at which the coefficient $C_1(p,N)$ vanishes, so that $f_0(y)$
  has exponential decay at infinity. To see this, we show in
  Figure \ref{FF1gg} with the scale of $10^{-20}$ how the
  coefficient $C_1(p,1)$ changes sign  around \ef{rr1}:
   \be
   \label{rr2}
  C_1(1.39,1)>0, \quad \mbox{while} \quad C_1(1.40,1)<0.
   \ee
Non vanishing of these two profiles in smaller scaled up to
$10^{-40}$ was checked in the logarithmic scale (we do not present
here a number of such numerics).

%%%%%%%%%%%%%%%%%%%%%%%%%%%%%%%%%%%%%%%%%%%%%%%%%%%%%%%%%%%%

%%%%%%%%%%%%%%%%%%%%%%%%%%%%%%%%%
\subsection{On non-radial self-similar blow-up patterns in dimensions $N \ge 2$}

This question was not studied in the literature at all and indeed
is very difficult. We make a slight observation only: the
performed below linearization \ef{4.1} about the constant
equilibrium in the elliptic equation \ef{2.2} leads to the
perturbed linear elliptic equation
 \be
 \label{lin11}
 (\BB^* + I)Y + {\bf D}(Y)=0
  \ee
(on spectral properties of $\BB^*$, see
%%% are listed in
 Lemma \ref{lemSpec2}).
Then, $\BB^*+I$ has a large unstable subspace
 \be
 \label{ww1}
 E^{\rm u}(0)= {\rm Span}\{ \psi_\b: \,\,\, \l_\b+1 >0\},
  \ee
 %%$W^{\rm u}(0)$,
so that the corresponding eigenfunctions may characterize possible
shapes of various similarity solutions (actually, this is true for
$N=1$ \cite{BGW1}). Roughly speaking, we claim that:
 \be
 \label{kk1}
  \fbox{$
  \begin{matrix}
 \mbox{the dimension $M(p,N)={\rm dim}\, E^{\rm u}(0)-(N+1)$
 %%of the unstable manifold of $\BB^*+I$
  can
 characterize}\ssk\\
 \mbox{the total number of blow-up similarity patterns as
 solutions of (\ref{2.2}).}
  \end{matrix}
  $}
  \ee
Note that we subtract $(N+1)$-dimensions corresponding to natural
instabilities relative to shifting the blow-up point $0 \in \ren$
($N$ dimensions) and blow-up time $T$ (1 dimension). These
unstable modes are not available if the blow-up point $(0,T)$ is
fixed. The dimension $M=M(p,N)$ can characterize the total number
of solutions $f(y)$ of the elliptic problem \ef{2.2} including
many non-radial ones. In other words, we expect that those
unstable $M$ modes initiate heteroclinic connections through the
corresponding unstable manifold $ W^{\rm u}(0)$ to the set of
steady solutions $\{f_k(y), \, k=1,2,...,M\}$. In Section
\ref{S.4}, we show that stable modes from $E^{\rm s}(0)$ with
$\l_\b+1>0$ and the centre ones $E^{\rm c}(0)$ with $\l_\b+1=0$
will lead to other ``linearized" blow-up patterns, so that
$\{f_k\}$ are ``nonlinear eigenfunctions".

Proving any part of the claim \ef{kk1} is a  difficult
%% and
%%challenging problem that is entirely open
open problem for any $N \ge 2$. Note also that, for the
second-order {\em quasilinear} counterpart \ef{mm1} ($m>1$ is
essential!), non radially symmetric self-similar blow-up patterns
have been known for more than thirty years; see \cite{KurdPot84}
and a survey \cite{Kurd90} for extra details.

%%%%%%%%%%%%%%%%%%%%%%%%%%%%%%%%%%%%%%%%%%%%%%%%%%%%%
\section{{\bf Type
I(log)}: self-similar patterns with angular logTW swirl}
 \label{S.3}

This is a simple idea for producing non-radial blow-up patterns,
but its consistency is quite questionable.

\subsection{Nonstationary rescaling}

 Dealing with non-self-similar blow-up,
instead of \ef{2.1}, we use the full similarity scaling:
 \be
 \label{3.1}
  \tex{
 u(x,t)=(T-t)^{-\frac 1{p-1}}v(y,\t), \quad y= \frac
 x{(T-t)^{1/4}}, \quad \t = - \ln(T-t) \to +\iy, \,\, t \to T^-.
}
  \ee
 Then $v(y,\t)$ solves the following parabolic equation:
  \be
  \label{3.2}
   \tex{
  v_\t= {\bf A}(v) \equiv - \D^2 v - \frac 14\, y \cdot \n v-\frac 1{p-1}\, v +
  |v|^{p-1}v\inB \ren \times (\t_0,\iy), \,\,\, \t_0=-\ln T,
  }
   \ee
   where ${\bf A}$ is the stationary elliptic operator in
   \ef{2.2}, so that similarity profiles (if any) are just
   stationary solutions of \ef{3.2}.

%%%%%%%%%%%%%%%%%%%%%%%%%%%%%%%%%%%%%%%%%%%%%%%%%%%%%
\subsection{Blow-up angular swirling mechanism}

   We begin with $N=2$, where $y=(y_1,y_2)$, and, in the corresponding polar coordinates $\{\rho,
   \var\}$, with $\rho^2=y_1^2+y_2^2$,
    \be
    \label{3.3}
     \tex{
     \D= \D_\rho + \frac 1{\rho^2}\, D^2_\var \whereA  \D_\rho = D^2_{\rho} +
     \frac 1\rho \, D_\rho \andA y \cdot \n=
     \rho D_\rho.
     }
     \ee
We next consider a TW in the angular direction by fixing the
angular dependence
 \be
 \label{3.4}
\var=  \s\, \t + \mu \equiv - \s \ln(T-t) + \mu, \quad \mu
\in(0,2\pi),
 \ee
 where $\s \in \re$ is a constant (a nonlinear eigenvalue).
 In the original independent variables $\{x,t\}$, \ef{3.4}
 represents a {\em blowing up logarithmic TW} in the angular direction with
 unknown wave speeds $\s$. In other words, \ef{3.4} assumes that blowing up as
 $t \to T^-$ is accompanied by
  a focusing TW-angular behaviour also in a logarithmic  blow-up
  manner.

%%%\ssk

Thus,  assuming the logTW angular  dependence \ef{3.4} of the
solution $v=v(y,\mu,\t)$, $\var=\s \t + \mu$, yields the equation
\be
  \label{3.5}
   \tex{
  v_\t= {\bf A}(v) - \s v_\mu \equiv - \D^2 v - \frac 14\, y \cdot \n v-\frac 1{p-1}\, v +
  |v|^{p-1}v -\s v_\mu \inB \ren \times (\t_0,\iy),
  }
   \ee
   where $ \t_0=-\ln
  T$.
In particular, this non-radial self-similar blow-up may be
generated by bounded steady profiles satisfying
\be
  \label{3.6}
   \tex{
  {\bf A}(f) - \s f_\mu \equiv - \D^2 f - \frac 14\, y \cdot \n f-\frac 1{p-1}\, f +
  |f|^{p-1}f -\s f_\mu=0 \inB \re^2.
  }
   \ee
For $\s \ne 0$, which, as we have mentioned, plays  the role of a
nonlinear eigenvalue, the blow-up behaviour with swirl corresponds
to periodic orbits as $\o$-limit sets; see a discussion in
\cite{GalJMP} to the NSEs \ef{NS1}.
 As a first approach to solvability of \ef{3.6}, one can assume
 branching of a solution $f(\rho,\mu)$ from the radial one $f_0$
 at $\s=0$. Then setting $f=f_0 + \s \psi^*+...$ yields that $\psi^*(\rho,\mu)$
 must be a nontrivial non-radial eigenfunction for $\l=0$:
  \be
   \label{3.9}
  {\bf A}'(f_0) \psi^* =0.
   \ee
On the other hand, branches of solutions $f$ of \ef{3.6} may occur
at a saddle-node bifurcation $\s=\s_* \ne 0$, where $\s_*$ belongs
to spectrum of the linear pencil
 $
 {\bf A}'(f) - \s D_\mu$. Both  eigenvalue problem  are very
 difficult, and we do not exclude the possibility that, overall,
the problem \ef{3.6} for $\s \ne 0$ may admit such solutions only
that are singular at the origin $y=0$. Anyway, even in this
unfortunate case, we believe that introducing such rather unknown
types of non-radial blow-up with swirl deserves mentioning among
other more practical patterns.
 %%% In general, existence for some $\s \ne 0$ of a bounded solution
%%\ssk
Let us also mention that, in $\ren$, one can distribute the
variables as
 $$
 y=(y_1,y_2,y') \in \ren \whereA y' \in \re^{N-2},
 $$
 and arrange a $\s_1$-logTW in variables $(y_1,y_2)$ only to get periodic blow-up behaviour.
 Choosing other disjoint  pairs $(y_k,y_{k+1})$ and constructing
 the corresponding periodic swirl in these variables, in particular,
 it is formally possible to
 produce a {\em quasi-periodic blow-up swirl} with arbitrary
 number $\s_1$,\,...,\,$\s_n$,  $n \le \big[\frac N2\big]$, of fundamental
 frequencies. Of course, this leads to complicated nonlinear
 eigenvalue problems, which are open even for $n=1$, i.e., for the
 periodic motion introduced above first.

%%\ssk

%%%%%%%%%%%%%%%%%%%%%%%%%%%%%%%%%%%%%%%%%%%%%%%%%%%%%%%%%
%%\noi{\bf
\subsection{Remark: on the origin of logTWs and invariant
solutions}
 The scaling group-invariant nature of such logTWs seems
was first obtained by Ovsiannikov in 1959 \cite{Ov59}, who
performed a full group classification of the nonlinear heat
equation
 $$
 u_t=(k(u)u_x)_x,
  $$
  for arbitrary functions $k(u)$.
  In particular, such invariant solutions appear for the porous
  medium and fast diffusion equations for $k(u)=u^n$, $ n \ne 0$:
   $$
   \tex{
   u_t=(u^n u_x)_x \LongA \exists \,\,\, u(x,t)=t^{-\frac 1n} f(x+ \s \ln
   t)\whereA
    - \frac 1n \, f + \s f'=(f^n f')'.
    }
    $$
 Blow-up angular dependence  as $t \to T^-$  such
as in \ef{3.4} was  studied later on in \cite{Bak88}, where the
corresponding similarity solutions for the reaction-diffusion
equation with source \ef{mm1} in $\re^2 \times (0,T)$
%% with anisotropic nonlinear coefficients
 %%\be
%% \label{an1}
%% u_t=\n \cdot (u^\s \n u)+ u^\b
%%% u_t= (u^{\s_1}u_x)_x + (u^{\s_2}u_y)_y + u^\b
%% \inB \re^2 \times (0,T) \quad (\b>1,\,\,\,\s>0),
 %%% \s_1>\s_2>0),
%%  \ee
were indicated by reducing the PDE to a quasilinear elliptic
problem (it seems, there is no still a rigorous proof of existence
of such patterns). For parabolic models such as \ef{mm1}, that are
 order-preserving
via the MP and
 do
not have a natural ``vorticity" mechanism, such ``spiral waves" as
$t \to T^-$ must be generated by large enough initial data
%%%that are
specially ``rotationally"  distributed in $\re^2$. For
 the biharmonic operator as in \ef{m2} with no MP,
 such a swirl blow-up dependence may be more relevant; see below.

 \ssk

 %%%%%%%%%%%%%%%%%%%%%%%%%%%%%%%%%%%%%%%%%%%%%%%%%%%%%
\section{{\bf Type
I(Her)}: non self-similar ``linearized" patterns with a local
generalized Hermite polynomial structure}
 \label{S.4}

 For the classic R--D equation \ef{FK1p}, a countable set of non self-similar  blow-up patterns
 of a similar structure
 was
first formally introduced  in \cite{GHPV}, though
 the history of such  non-self-similar blow-up asymptotics goes back
 to  Hocking--Stuartson--Stuart
in 1972, \cite{HSS}, who invented an interesting novel formal
technique of analytic expansions (in fact, an analogy of a centre
manifold analysis) to confirm that blow-up occurs on subsets
governed by the ``hot spot" variables, as $t \to T^-$:
 \be
 \label{hs1}
 \tex{
u(x,t) = (T-t)^{-\frac 1{p-1}}\, \big[f_*(\xi)+o(1)\big] \whereA
\xi = \frac x{\sqrt{(T-t) | \ln (T-t)|}}
   }
   \ee
   and $f_*$ is a unique solutions of a Hamilton-Jacobi equation
   of the form \ef{4.72} (with $\frac 14 \mapsto \frac 12$).
  A justified
construction of such patterns and other applications were
performed a few years later mainly in dozens of  papers by Herrero
and Vel\'azquez; see \cite{Herr93, Vel} as a guide together with
other papers traced by the {\tt MathSciNet}. It is curious that
earlier, in 1987,  a sharp upper bound of such a non-similarity
blow-up evolution \ef{hs1} (``first half of blow-up") was proved
in \cite{GPos87} by a modification of Friedman--McLeod gradient
estimate \cite{FM86}, though the ``second half of blow-up" took
extra ten years co complete along similar lines,
\cite[\S~7]{GalDS}.

For the RDE--4, there is no  hope to get an easy and fast rigorous
justification of such non-self-similar blow-up scenario, though
the main idea remains the same. We follow \cite{Gal2m} and also
 \cite{GalCr}, where such a construction applied to
non-singular absorption phenomena (regular flows with
 no blow-up), so a full mathematical justification is available therein.

%%%%%%%%%%%%%%%%%%%%%%%%%%%%%%%%%%%%%%%%%%%%%%%%%%%%%%%%
\subsection{Linearization and spectral properties}

The construction of such blow-up patterns is as follows.
Performing the standard linearization about the constant
equilibrium in the equation \ef{3.2} yields the following
perturbed equation:
 \be
 \label{4.1}
 \tex{
 v=f_*+Y, \,\, f_*=(p-1)^{-\frac 1{p-1}}  \LongA
Y_\t= (\BB^*+I) Y  +\DD(Y),
%% \,\, \mbox{where}\,\, \DD(Y)=c_0
%%Y^2+...
 }
 \ee
where $\DD(Y)=c_0 Y^2+...$, $c_0=\frac {p}2(p-1)^{\frac{1}{p-1}}$,
is a quadratic perturbation as $Y \to 0$ and
  %% which
 %%does not play an essential role for further formal calculus.
 %%% In \ef{4.1},
  \be
  \label{4.2}
  \tex{
  \BB^*=-\D^2 - \frac 1{4}\, y \cdot \n \inB L^2_{\rho^*}(\ren),
  \,\,\,\rho^*(y)={\mathrm e}^{-a|y|^{4/3}}, \quad a \in \big(0,3 \cdot
 2^{- \frac {8}3}\big),
  }
  \ee
 %% where $a \in (0,3 \cdot
 %%2^{- {8}/3})$ is a constant,
  is the adjoint Hermite operator with some good spectral
  properties \cite{Eg4}:

  \begin{lemma}
\label{lemSpec2} %%(i)
  $\BB^*: H^{4}_{\rho^*}(\ren) \to
L^2_{\rho^*}(\ren)$ is a bounded linear operator with the spectrum
 \be
\label{SpecN}
 \tex{
 \s(\BB^*) = \{ \l_\b =  -\frac{|\b|}4, \,\,\,
|\b|=0,1,2,...\} \quad \big(= \s(\BB), \,\, \BB=-\D^2 + \frac
1{4}\, y \cdot \n+ \frac N4 \, I \big). }
   \ee
     Eigenfunctions
$\psi_\b^*(y)$ are $|\b|$th-order generalized Hermite polynomials:
\be
\label{psidec}
 \tex{
\psi_\b^*(y) =  \frac {1}{ \sqrt{\b !}}\big[ y^\b +
\sum_{j=1}^{[|\b|/4]} \frac 1{j !}(\Delta)^{2 j} y^\b \big],
 \quad |\b|=0,1,2,...\, ,
 }
  \ee
%%    (ii)
and the subset $\{\psi_\beta^*\}$ is complete
%%and closed
 in
$L^2_{\rho^*}(\ren)$.
%% (iii) $\BB^*$ is sectorial in
%%%%$l^2_{\rho^*}$.

\end{lemma}

As usual, if $\{\psi_\b\}$ is the adjoint basis of eigenfunctions
of the adjoint operator
 \be
 \label{BBB}
 \tex{
 \BB= - \D^2 + \frac 14\, y \cdot \n + \frac N4\, I \inB
 L^2_\rho(\ren) \withA \rho= \frac 1{\rho^*},
 }
\ee with the same spectrum \ef{SpecN}, the bi-orthonormality
condition holds in $L^2(\ren)$:
%%% \cite{Eg4}:
 \be
 \label{Ortog}
\langle \psi_\mu, \psi_\nu^* \rangle = \d_{\mu\nu} \quad \mbox{for
any} \quad \mu, \,\, \nu.
 \ee

\subsection{Inner expansion}

Thus, in the {\em Inner Region} characterized by compact subsets
in the similarity variable $y$, we assume a centre or a stable
subspace behaviour as $\t \to +\iy$ for the linearized operator
$\BB^*+I$:
 \be
 \label{4.5}
  \begin{matrix}
 \mbox{centre:} \quad Y(y,\t)=a(\t) \psi_\b^*(y) + w^\bot \quad (\l_\b=-1, \,\, |\b|=4),\quad\,\,\ssk\ssk\\
 \mbox{stable:} \quad  Y(y,\t) =-C {\mathrm e}^{\l_\b \t} \psi_\b^*(y)+w^\bot \quad (\l_\b <-1, \,\, |\b|>4).
    \end{matrix}
    \ee
For the centre subspace behaviour in \ef{4.5}, substituting the
eigenfunctions expansion into equation \ef{4.1} yields the
following coefficient:
 \be
 \label{4.6}
  \tex{
  \dot a=a^2 \g_0+... \whereA \g_0= c_0 \langle (\psi_\b^*)^2,
\psi_\b \rangle
  \LongA a(\t) =- \frac {1}{\g_0\, \t}+... \, .
  }
  \ee
 Note that for the matching purposes, we have to
 assume that (see details in \cite{Gal2m}):
 \be
 \label{4.7}
\mbox{if $\psi_\b^*(0)>0$, then} \quad  \g_0>0 \andA C>0.
  \ee
Actually, for $k=4$ and $N=1$,  it is calculated explicitly that
 \be
 \label{4.71d}
  \tex{
  \g_0=-c_0 136 \sqrt 6<0,
   }
   \ee
 so that the centre manifold patterns with the positive eigenfunction
 \be
 \label{4.72d}
  \tex{
   \psi_4^*(y)= \frac 1{\sqrt{24}}\, (y^4+ 24)
    }
    \ee
correspond to solutions that blow-up on finite interfaces; see
\cite[\S~3]{GW1}. A full justification of such a behaviour can be
done along the lines of classic invariant manifold theory (see
e.g. \cite{Lun}), though can be very difficult.

Actually, we can construct more general asymptotics  by taking an
arbitrary linear combination of eigenfunctions from the centre
subspace. Overall, the whole variety of such asymptotics is
characterized as follows:
 \be
 \label{YY1}
  %%\tex{
  Y(y,\t) = a(\t) |y|^\b \chi(\var)+...\, \whereA a(\t) = -
 \left\{
  \begin{matrix}
 \,\,\,
  \frac
  1{\g_0 \t}+... \,\,\,
  \quad
  \mbox{for}\quad |\b|=4,\ssk\\
  %%% \,\, \mbox{and} \,\, a(\t) \sim
  C{\mathrm e}^{\l_\b \t}+... \quad \mbox{for} \quad |\b|>4.
 \end{matrix}
 \right.
 %% }
  \ee
 In general,  here, $\chi(\var)>0$ is an arbitrary
 smooth function on the sphere ${\mathbb
  S}^{N-1}$, where its positivity is induced by matching issues to
  be revealed below.

%%%%%%%%%%%%%%%%%%%%%%%%%%%%%%%%%%%%%%%%%%%%
\subsection{Outer region: matching}

We follow \cite{Gal2m}, where it is shown that the asymptotics
\ef{4.5} admit matching with the {\em Outer Region}, being a
Hamilton--Jacobi (H--J) one. More precisely, in the centre case
with $|\b|=4$, according to \ef{4.5}, \ef{4.6}, we introduce the
outer variable and obtain from \ef{3.2}  the following perturbed
H--J equation:
 \be
 \label{4.71}
 \tex{
 \xi= \frac y {\t^{1/4}} \LongA v_\t= - \frac 1 4\, \xi \cdot \n v- \frac 1{p-1}\, v + |v|^{p-1}
 v + \frac 1 \t \, \big( \frac 14\, \xi \cdot \n v -\D^2 v \big).
 }
 \ee
Passing to the limit as $\t \to +\iy$ in such singularly perturbed
PDEs is not easy at all even in the second-order case (see a
number of various  applications in \cite{AMGV}). Though currently
not rigorously (this looks being completely illusive), we assume
stabilization to the stationary solutions $f(\xi)$ satisfying the
unperturbed H--J equation:
 \be
 \label{4.72}
 \tex{
 - \frac 14\, \xi \cdot \n f- \frac 1{p-1}\, f + |f|^{p-1}
 f=0 \inB \ren.
 }
 \ee
 This is solved via characteristics, where we
 have to choose the solution satisfying
 %%by \ef{4.5}
 %%and
 \ef{YY1}:
%%% (other non-symmetric patterns  also occur):
  \be
  \label{4.73}
   \tex{
  f(\xi)= f_*  - \frac 1{\g_0}\, |\xi|^4 \, \chi(\var)
  +... \asA
  \xi \to 0 \LongA f(\xi)=f_*(1+c_* |\xi|^4 \chi(\var))^{-\frac
  1{p-1}},
  }
  \ee
  where $c_*=\frac 1{\g_0}\,  (p-1)^{p/(p-1)}$.
   Since $\g_0<0$ according to \ef{4.71d}, the resulting
   profile satisfies
   $f(\xi) \ge f_*$ and  blows up on the surface $\{c_* |\xi|^4 \chi(\var) =
   -1\}$.
    Note that this actual nonexistence of a bounded centre
    subspace pattern plus the known unstable eigenspace of $\BB^*+I$
    in the linearized equation \ef{4.1}
  somehow reflect existence of two self-similar solutions
  $f_0(y)$ and $f_1(y)$ as ``nonlinear eigenfunctions";
    see %%%% analogous phenomena
    \cite{GW1}.

   %%We claim that such an unusual behaviour serves as
   %%a  heteroclinic connection $f_1 \mapsto 0$.

Thus, these centre subspace patterns are not bounded and should be
excluded from the consideration. On the other hand,
    for the $2m$th-order PDE \ef{mm98} with odd
   $m=3,5,...\,$, we have $\g_0>0$ and then \ef{4.73} can represent
   standard blow-up patterns, \cite{Gal2m}.

Similarly, for the stable behaviour for $|\b|>4$ in \ef{4.5}, we
use the following change:
\be
 \label{4.71S}
 \tex{
 \xi=  {\mathrm e}^{ \frac{4-|\b|}{4|\b|}\, \t}\, y  \LongA v_\t=
 - \frac 1{|\b|}\, \xi \cdot \n v- \frac 1{p-1}\, v + |v|^{p-1}
 v - {\mathrm e}^{- \frac {|\b|-4}{|\b|}\, \t} \, \D^2 v.
 }
 \ee
Passage to the limit $\t \to +\iy$ and matching with Inner Region
are analogous and lead to truly existent blow-up patterns for any
$m \ge 2$ in \ef{mm98}, \cite{Gal2m}.

%%It is worth mentioning that,
Overall,  according to matching conditions \ef{YY1} and
\ef{4.71S},
%% not being proved rigorously,
the whole set of  possible blow-up patterns of Type II(Her) is
composed from a countable
 set for $|\b|=4$ ($m$ odd),\,5,\,6,\,... of continuous (uncountable) families
 induced by smooth functions $\chi$ on ${\mathbb S}^{N-1}$.

 %%%%%%%%%%%%%%%%%%%%%%%%%%%%%%%%%%%%%%%%%%%%%%%%%%%%%
\section{{\bf Type
II(sing)}: linearization about the SSS and matching}
 \label{S.5}

 The idea of such Type II blow-up patterns for the RDE--2
 \ef{FK1p} is due to Herrero--Vel\'azquez \cite{HVsup}, where
 a justification of existence was achieved (see  \cite{Miz07}
 for extra details).
 We apply this method to the RDE--4 \ef{m2} and, by the same reasons, we are not
 obliged to
 concentrate  on a proof.
 Thus, instead of the linearization \ef{4.1} about the constant
 equilibrium, we perform it about a singular one.

%%%%%%%%%%%%%%%%%%%%%%%%%%%%%%%%%%%%%%%%%%
\subsection{Singular stationary solution (SSS)}

Consider the stationary equation
 \be
 \label{5.1}
 -\D^2 U + |U|^{p-1}U=0 \inB \ren \setminus \{0\}.
  \ee
  The explicit radial SSS has the standard scaling invariant form
   \be
   \label{5.2}
    \begin{matrix}
    U(y)=
C_* \, |y|^{-\mu} \whereA \mu= \frac 4{p-1}, \quad
  C_*= D^{\frac 1{p-1}}, \ssk\ssk\\
 \mbox{and} \quad
  D=\mu(\mu+2)[(\mu+1)(\mu+3)+(N-1)(N-5-2\mu)].
   \end{matrix}
   \ee
It follows that such SSS exists, i.e., $D>0$,
 in the following
parameter ranges:
 \be
 \label{5.3}
 \tex{
  p>
 \frac N{N-4}, \,\,\,N>4, \quad \mbox{or} \quad p < \frac {N+2}{N-2},\,\,\,N>2.
 %%%\end{matrix}
 }
 \ee

%%%%%%%%%%%%%%%%%%%%%%%%%%%%%%%%%%%%%%%%%%
\subsection{Linearization in Inner Region: discrete
spectrum by  Hardy--Rellich inequality}

Thus, we perform linearization in \ef{3.2} about the SSS:
 \be
 \label{5.4}
  \tex{
v=U+Y \LongA v_\t= \hat \BB^* Y + \DD(Y), }
 \ee
 where, as usual, $\DD(Y)$ is a quadratic perturbation as $Y \to
 0$ and
  \be
  \label{5.5}
   \tex{
  \hat \BB^*= \HH^* - \frac 14\, y \cdot \n - \frac 1{p-1} \, I
  \andA
  \HH^*= -\D^2 + \frac {p D}{|y|^4}\, I.
   }
   \ee
Similar to Lemma \ref{lemSpec2}, the operator $\hat \BB^*$ at
infinity admits a proper functional setting in the same metric of
$L^2_{\rho^*}$. However, it is also singular at the origin $y=0$,
where its setting depends on the principal part $\HH^*$.

\begin{proposition}
\label{Pr.Har1} The symmetric operator $\HH^*$ admits a
Friedrich's self-adjoint extension  with the domain $H^4_0(B_1)$,
discrete spectrum,  and compact resolvent in $L^2(B_1)$, where
$B_1 \subset \ren$ is the unit ball, iff
 \be
 \label{5.6}
  \tex{p D \le c_{\rm H}= \frac{[N(N-4)]^2}{16}.
  }
  \ee
 %% \noi{\rm (ii)}
  \end{proposition}

  \noi{\em Proof.} Indeed, \ef{5.6} is just a corollary of the
  classic
  {\em Hardy--Rellich}-type {\em  inequality}\footnote{This was derived by
  Rellich already in 1954;  see \cite{GGM} and \cite{Yaf99} for
further references and full history.}
   \be
   \label{5.7}
   \tex{
 \frac {[N(N-4)]^2} {16}  \int\limits_{B_1} \frac {u^2}{|y|^4} \le
 \int\limits_{B_1}
    |\D u|^2 \forA u \in H^2_0(B_1),
    }
    \ee
    where the constant is sharp. For compact embedding of the corresponding spaces,
    see Maz'ya \cite[p.~65, etc.]{Maz}. $\qed$

    \ssk

The necessary inequality \ef{5.6} takes the form
 \be
 \label{5.8}
  %%\tex{
   \begin{matrix}
  G_p(N) \equiv
  \frac{[N(N-4)]^2}{16} -
  \frac{4p}{p-1}\big(2+\frac{4}{p-1}\big)
  \qquad\qquad\qquad\qquad\qquad\,\,
  \ssk\ssk\\
  \times \big[\big(1+\frac{4}{p-1}\big)
   \big(3+\frac{4}{p-1}\big)+(N-1)(N-5- \frac{8}{p-1}\big)\big]
   \ge 0
  \end{matrix} %%% }
   \ee
  and does not admit an easy analytic analysis. In Figure
  \ref{Fp1}, numerics show that
  %%% we present numerical results showing that
   \be
   \label{5.9}
   \mbox{(\ref{5.8}) holds for $N \ge 24$ if $p=2$, and $N \ge 19$ if $p=3$.}
    \ee
  In particular,   checking \ef{5.8} at $p=+\iy$ yields the inequality:
 \be
 \label{YY2}
  \tex{
  G_\iy(N) \equiv \frac{[N(N-4)]^2}{16} -8[3+(N-1)(N-5)] >0.
  }
  \ee
If this is true, then \ef{5.8} holds for all $p \gg 1$, so:

    %% the sign of the last term on the right-hand side of
    %%%\ef{5.8} yields:

%%FIG%%%%%%%%%%%%%%%%%%%%%%%

\begin{figure}
%%\vskip -.3cm
\centering \subfigure[$p=2$]{
\includegraphics[scale=0.52]{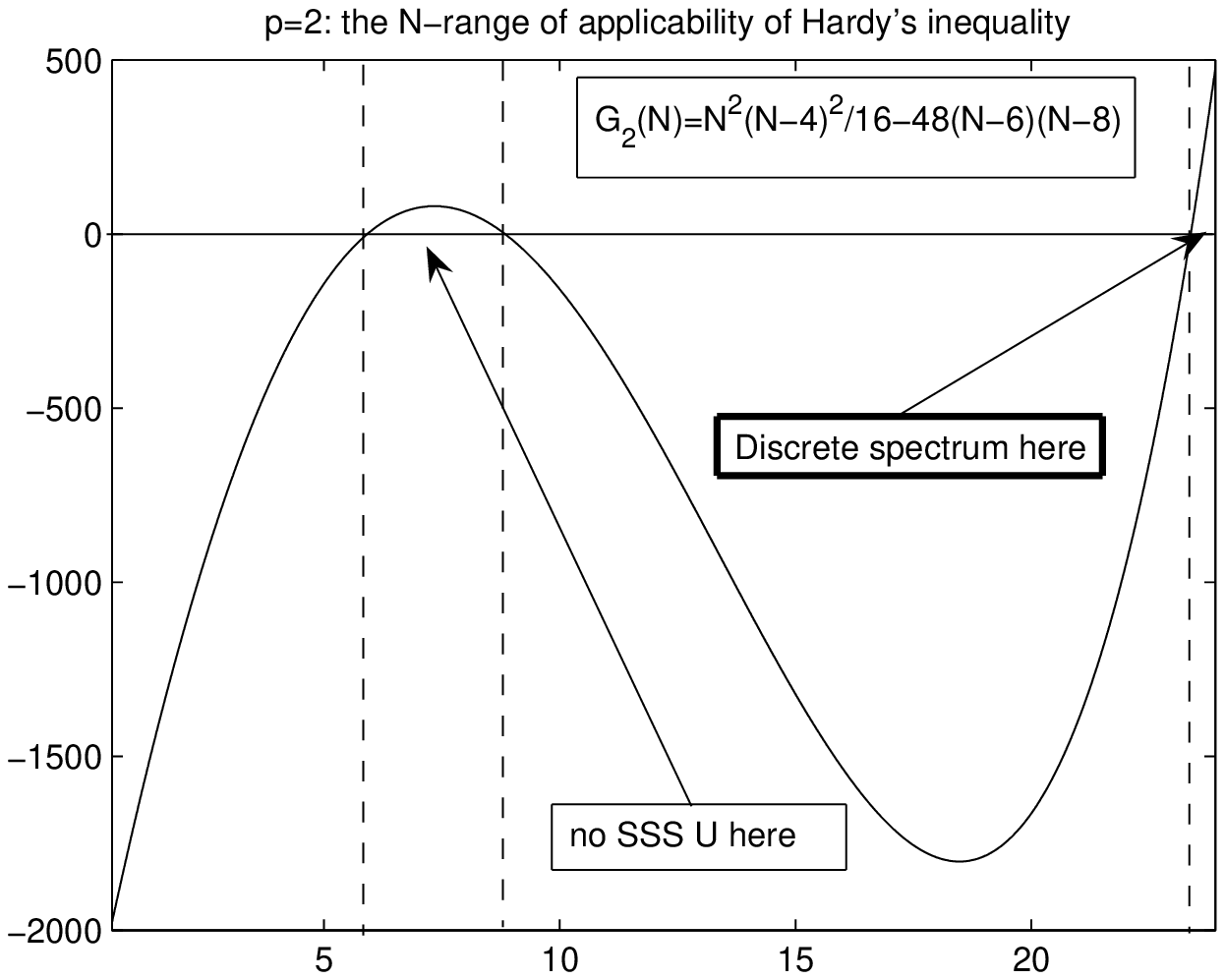}
} \subfigure[$p=3$]{
\includegraphics[scale=0.52]{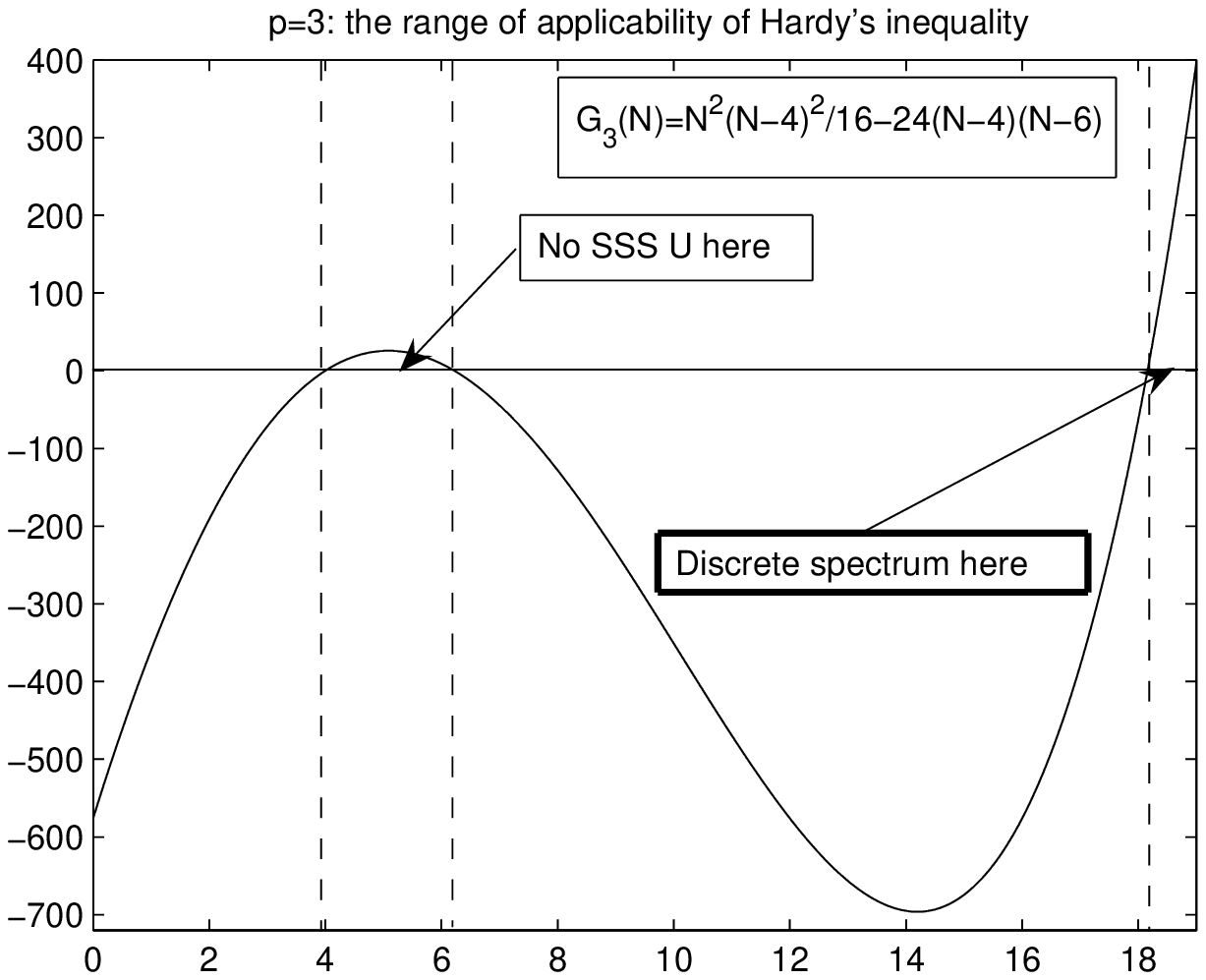}
}
 \vskip -.2cm
\caption{\rm\small Checking inequality \ef{5.8}: $p=2$ (a) and
$p=3$ (b).}
 %% \vskip -.3cm
 \label{Fp1}
\end{figure}
%%%%%%%%%%%%%%%%%%%%%%%%%%%%%%%%%%%%%%%%%%%%%%%%%%%%%%%%%%%%%%%%%%%

 \begin{proposition}
  \label{Pr.Har2}
  For any $N \ge 13$, there exists a $p_{\rm H}(N)>1$ such that
 \be
 \label{5.10}
  \mbox{$(\ref{5.8})$
 holds for all $p \ge p_{\rm H}(N)$},
  \ee
   and hence the operator $\hat \BB^*$ in $(\ref{5.5})$ has a
 discrete spectrum in $L^2_{\rho^*}(\ren)$.
 \end{proposition}

%%%%%%%%%%%%%%%%%%%%%%%%%%%%%%%%%%%%%%%%
\subsection{Inner Region I}

Thus, we assume that, under certain conditions, \ef{5.10} holds
and $\s(\hat \BB^*)=\{\hat \l_k\}$ is discrete, with the
eigenfunctions $\{\hat \psi^*_\b, \, |\b|=k\}$. Furthermore, it is
also convenient to assume that the spectrum is (at least
partially) {\em real}. To justify such an assumption for this
non-self-adjoint operator, we rewrite \ef{5.5} in the form
 \be
 \label{5.11}
  \tex{
 \hat \BB^*= \BB^* + \frac c{|y|^4}\, I - \frac 1{p-1}\, I \whereA c=p D
 }
  \ee
 and $\BB^*$ is the previous operator \ef{4.2} with the real
 spectrum shown in Lemma \ref{lemSpec2} (actually, this
 means that $\BB^*$ admits a natural self-adjoint representation
 in the space $l^2_{\rho^*}$ of sequences, where it is also sectorial,
 \cite{2mSturm}). Therefore, the real spectrum of \ef{5.11} can
 be obtained by branching-perturbation theory (see Kato \cite{Kato}) from that $\{\l_\b=- \frac {k}{4}- \frac 1{p-1}, \, k=|\b| \ge 0\}$ of
 $\BB^*-\frac 1{p-1}\, I$ at $c=0$. Next, the branch must be  extended to $c=p D$, which is
 also a difficult mathematical problem; see \cite[\S~6]{GK1} for some extra details, which are not
  necessary here in such a formal blow-up analysis.

Thus, we fix a certain exponentially decaying pattern in {\em
Inner Region} I:
 \be
 \label{5.12}
 Y(y,\t) = C{\mathrm e}^{\hat \l_\b \t} \hat \psi_\b^*(y)+... \asA
 \t \to +\iy \quad (\hat \l_\b < 0).
  \ee
If there exists $\hat \l=0 \in \s(\hat \BB^*)$, the expansion will
mimic that in \ef{4.5} for the centre subspace  case. Note that
\ef{5.12} includes all the non-radial linearized blow-up patterns.

%%%%%%%%%%%%%%%%%%%%%%%%%%%%%%%%%%%%%%%%%%%%%%%%%%%%%%%%
\subsection{Matching with Inner Region II close to the origin}

In order to match \ef{5.12} with  a smooth bounded flow close to
$y=0$, which we call {\em Inner Region} II, one needs the
behaviour of the eigenfunction $\hat \psi_\b^*(y)$ as $y \to 0$.
To get this, without loss of generality, we assume the radial
geometry. Then,
%% as $y \to 0$,
 the principal operator in the
eigenvalue problem
 \be
 \label{5.131}
  \tex{
  \HH^* \hat \psi^*+... = \l \hat \psi^* \asA y \to 0
  }
  \ee
yields the following characteristic polynomial (see \cite{GalH1}):
 \be
 \label{5.13}
  \hat \psi^*(y) = |y|^\g+... \LongA
 H_c(\g) = -\g(\g-2)(\g+N-2)(\g+N-4) + c=0.
 %% \equiv -[\g+\mbox{$\frac {(N-4)}2$}
 %%]^2 [\g^2+(N-4)\g -
 %%\mbox{$\frac {N^2}4$}].
 \ee

 Consider the most interesting critical and extremal case
 \be
 \label{5.14}
  \tex{
c \equiv  p D = c_{\rm H}= \frac{[N(N-4)]^2}{16} \LongA
 H_c(\g)  \equiv -[\g+\mbox{$\frac {(N-4)}2$}
 ]^2 \big[\g^2+(N-4)\g -
 \mbox{$\frac {N^2}4$}\big].
 }
 \ee
  There exists the double root
 $
 \g_{1,2}= - \frac{N-4}2< 0,
 $
which generates two $L^2$-behaviours:
%% as $y \to 0$:
 \begin{equation}
 \label{3.61}
\hat \psi_1^*(y) = |y|^{-\frac{N-4}2}\ln |y|(1+o(1)) \quad
\mbox{and} \quad \hat \psi_2^*(y) = |y|^{-\frac{N-4}2}
(1+o(1))\quad \mbox{as} \,\,\, y \to 0.
 \end{equation}
Note that $H^2_0$-approximations of $\hat \psi_2^*$ establish that
$c_{\rm H}$ is the best constant in (\ref{5.6}). Other two roots
of the characteristic equation in \ef{5.14} are
 \begin{equation}
 \label{3.7}
 \g_{3,4} =  \mbox{$\frac 12$}\, \big[\, 4-N \mp \sqrt{(N-4)^2+N^2}\,\big],
 \end{equation}
 where $\g_3 < \g_{1,2}< 0$ and $\g_4>0$ corresponds to $L^2$-solutions.
 We have
 \begin{equation}
 \label{3.8}
\hat \psi_3^*(y)= |y|^{\g_3}(1+o(1)) \not \in L^2,
 \end{equation}
 so that in $L^2$ the deficiency indices
 of $\BB$ are  $(3,3)$ and cannot be equal to $(4,4)$.
Unlike the second-order case,  the straightforward conclusion on
the discreteness of the spectrum in the case $(4,4)$
\cite[p.~90]{Nai1} does not apply, so Friedrich's extension of
$\HH^*$ is constructed by other arguments \cite{GalH1} and include
settings, where two most singular behaviour in \ef{3.61} and
\ef{3.8} are excluded.

Overall, this gives the following behaviour of the proper
eigenfunctions at the origin:
 \be
 \label{5.15}
 \tex{
 \hat \psi_\b^*(y) =- \nu_\b |y|^{-\frac{N-4}2} +... \asA y \to 0
 \quad (\nu_\b >0 \,\,\, \mbox{are normalization constants}).
  }
  \ee
 %%First of all,
This allows to detect the rate of blow-up of such
patterns by estimating the maximal value of the expansion near the
origin:
 \be
 \label{5.17}
  \tex{
  v_\b(y,\t)= C_*|y|^{-\frac 4{p-1}}- \nu_\b C {\mathrm e}^{\hat
  \l_\b \t} |y|^{-\frac{N-4}2} +... \asA y \to 0 \andA \t \to
  +\iy,
  }
  \ee
  where we observe the natural condition of matching:
   \be
   \label{5.18}
   \nu_\b C>0.
    \ee
 Calculating the absolute maximum in $y$ of the function on the
 right-hand side of \ef{5.18} (this is a standard and justified trick in some R--D
 problems; see e.g., \cite{Dold1})  yields an exponential
 divergence:
 \be
 \label{5.19}
 \tex{
 \|v_\b(\cdot,\t)\|_\iy = d_\b {\mathrm e}^{\rho_\b\t}+...\,
 \whereA \rho_\b= \frac {8|\hat \l_\b|}{(N-4)(p-p_{\rm S})}>0
 \quad \big(p > p_{\rm S}\big),
 }
 \ee
 where $d_\b>0$ are some constants. Note that, depending on the
 spectrum $\{\hat \l_\b<0\}$, \ef{5.19} can determines a countable
 set of various Type II blow-up asymptotics.

 Let us define more
%%Of course, we then need to show more
clearly the necessary matching procedure. In a standard manner,
%%(cf. \ef{4.71S}),
we return to the original rescaled equation \ef{3.2} and perform
the rescaling in Region II according to \ef{5.19}:
 \be
 \label{5.20}
  \tex{
  v(y,\t) = {\mathrm e}^{\rho_\b \t} w(\xi,s), \quad%%% \whereA
  \xi={\mathrm e}^{\mu_\b \t}y, \quad \mu_\b=
  \frac{(p-1)\rho_\b}4, \quad %%\andA
  s= \frac 1{(p-1)\rho_\b} \, {\mathrm e}^{(p-1)\rho_\b \t}.
  }
  \ee
 Then $w$ solves the following exponentially perturbed uniformly parabolic
 equation:
  \be
  \label{5.21}
   \tex{
   w_s= - \D^2 w + |w|^{p-1}w - \frac 1{(p-1)\rho_\b} \, \frac 1
   s\, \big[\big(\frac 14 + \mu_\b\big) \xi \cdot \n w+ \big(
   \frac 1{p-1} + \rho_\b\big) w \big].
    }
    \ee
As above, we arrive at a stabilization problem to a bounded
stationary solution, which is widely used in blow-up applications
(see examples in \cite{AMGV}). In general, once the uniform
boundedness of the orbit $\{w(s), \,\, s>0\}$ is established, the
passage to the limit in \ef{5.21} as $s \to +\iy$ is a standard
issue of asymptotic parabolic theory.

Our blow-up patterns correspond to the stabilization uniformly on
compact subsets:
 \be
 \label{5.22}
 w(\xi,s) \to W(\xi), \,\, s \to +\iy, \,\,\mbox{where} \,\, -\D^2 W+ |W|^{p-1}W=0, \,\, \xi \in\ren, \,\, W(0)=d_\b,
  \ee
for all admissible $|\b|=0,1,2,...\,$. We next discuss a crucial
issue on such a matching.

%%%%%%%%%%%%%%%%%%%%%%%%%%%%%%%%%%%%%%%%%%%%%%%%%%%
\subsection{Matching: on necessary structure of global bounded stationary solutions}

There are two issues associated with the stationary problem
\ef{5.22}.

\begin{figure}
%  \vskip -.3cm
\centering
\includegraphics[scale=0.6]{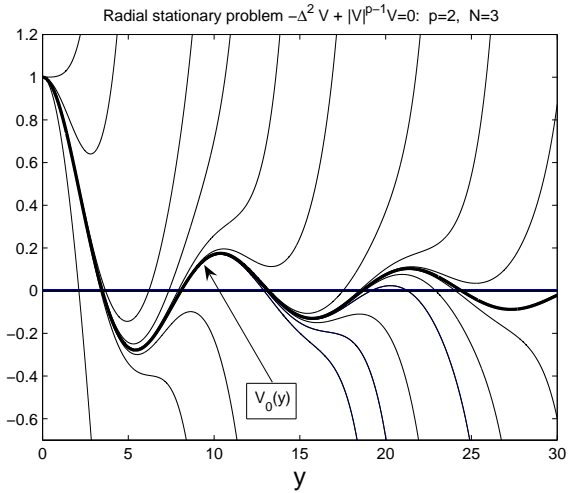} %%%%{AmFBP1.eps}  %%%%%%%{4het5.eps}   %%%%%{4het.eps}  Old
\vskip -.3cm \caption{\small  Shooting an oscillatory solution of
\ef{5.22} for $p=2$ and $N=3$.}
  %%% \vskip -.3cm
 \label{Fos1}
\end{figure}

{\bf 1.} Firstly and elementary, one can see that, bearing in mind
the matching of Regions I and II, the bounded stationary solutions
$W(\xi)$ defined by \ef{5.22} must be {\em positive} and {\em
non-oscillatory} as $\xi \to \iy$. Otherwise, such a matching with
{\em positive} SSS $U(\xi)$ is impossible. There exists a definite
negative result in the subcritical Sobolev range (there is a
diverse literature on this popular nowadays subject, so we refer
to a recent paper \cite{Guo08} as a guide):
 \be
 \label{5.30}
 \mbox{a solution $W>0$ of (\ref{5.22}) is nonexistent for
 $p \in(1,p_{\rm S})$.}
  \ee
Actually, this means that all the entire (i.e., without
singularities) solutions of \ef{5.22} are oscillatory, as Figure
\ref{Fos1} shows for $p=2$, $N=3$. Note that we are restricted by
\ef{5.6}.

%%%%%%%%%%%%%%%%%%%%%%%%%%%%%%%%%%%%%%%%%%%%%%%%%%%%%%%%%%%%
%%%%%%%%%%%%%%%%%%%%%%%%%%%%%%%%%%%%%%%%%%%%%%

\ssk

{\bf 2.} Secondly and fortunately, existence of such positive
solutions $W(\xi)$ is well established already
\cite[p.~908]{Gaz06}:
 \be
 \label{nn1}
 \mbox{for $p > p_{\rm S}$, for any $d_\b>0$, there exists a unique
 positive solution $W(\xi)$.}
  \ee

\begin{figure}
%%\vskip -.3cm
\centering \subfigure[$N=13>12=N_2$]{
\includegraphics[scale=0.52]{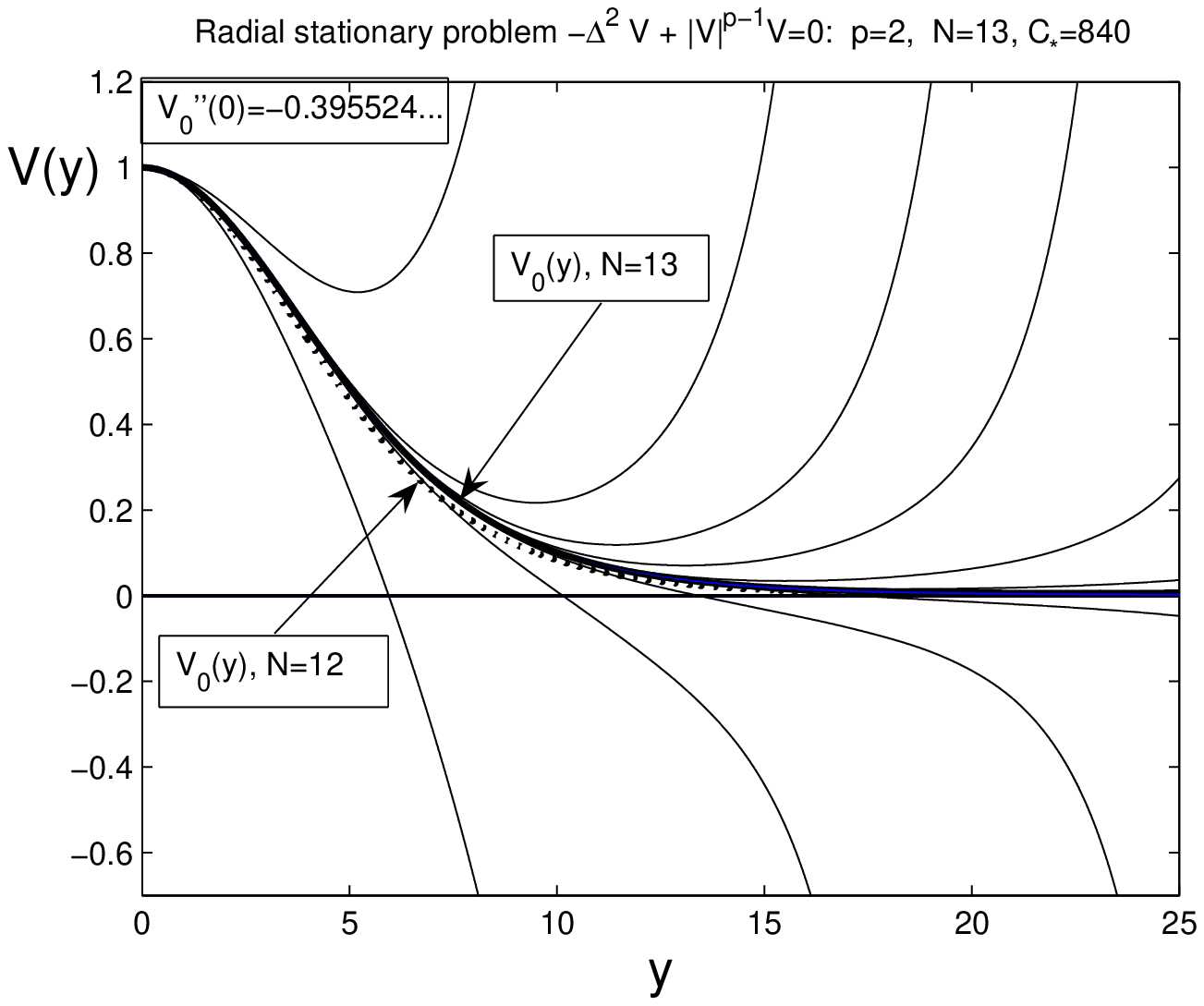}
} \subfigure[$N=24$]{
\includegraphics[scale=0.52]{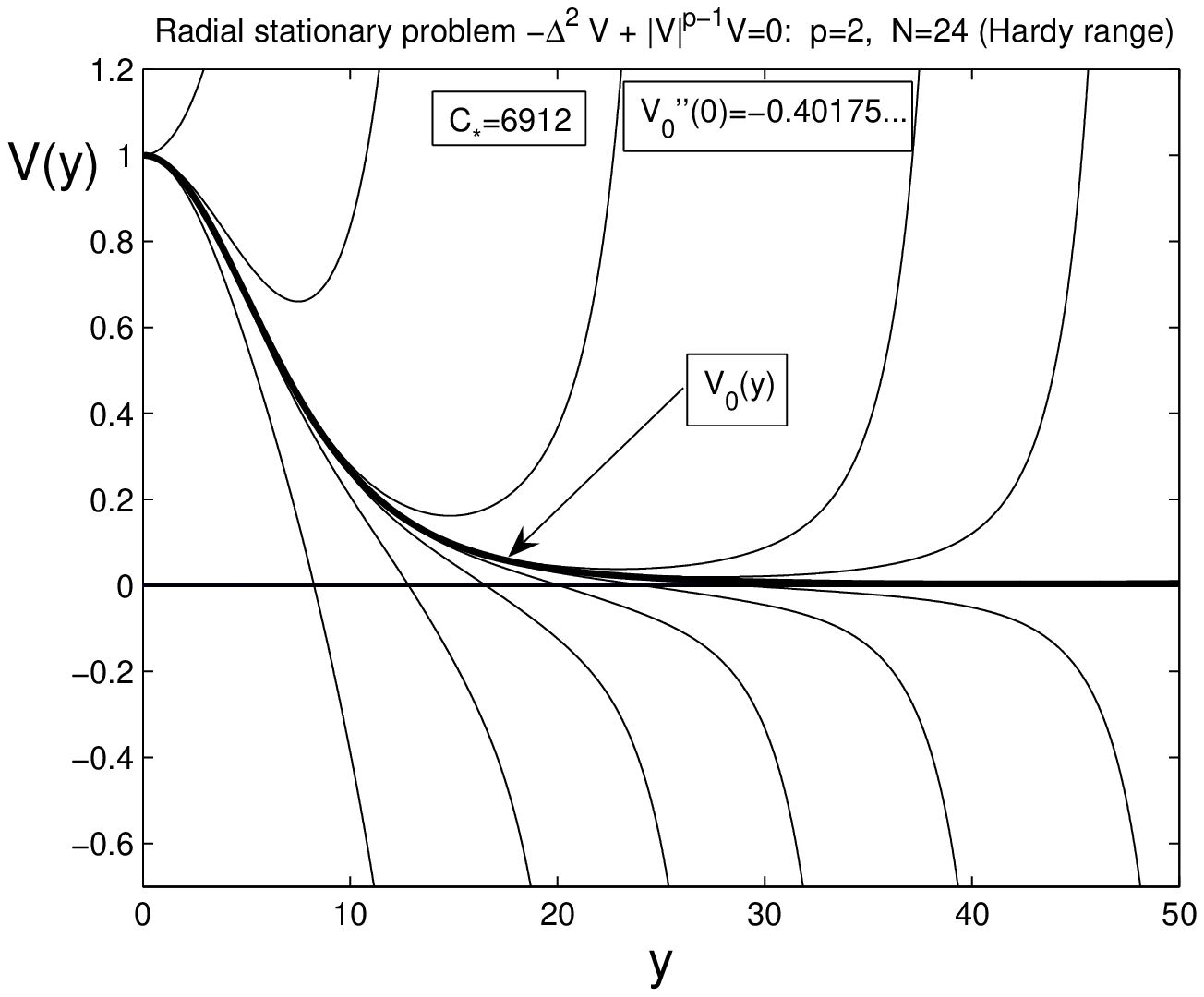}
}
 \vskip -.2cm
\caption{\rm\small Checking the positivity of  the solutions of
\ef{5.22} for $p=2$: $N=13$ (a) and $N=24$ (b).}
 %% \vskip -.3cm
 \label{Fos2}
\end{figure}

Here we exclude the critical case $p=p_{\rm S}$, where exact
positive solutions exist to be  used in Section \ref{S.6}. As a
numerical illustration, Figure \ref{Fos2} shows two such results
for $N=13$ (a), where the dotted line denotes the explicit
solution for $N=12$. It is clearly seen that $W(\xi)$ for $N=13$
lies above this, so remains positive.  In (b), we show the
positivity of the solution $W$ for $N=24$, where by \ef{5.9}, the
spectrum is guaranteed to be discrete.

\ssk

{\bf 3.} Of course, the above is not sufficient for matching of
Inner Regions I and II to get a  blow-up pattern. More
importantly, we have the following:

\begin{proposition}
  \label{Pr.Har3}
 %% In the same range $(\ref{5.10}$
   %%$p \ge p_{\rm H}$, $N \ge 13$,
 The entire solutions $W(\xi)$ of the radial ODE $\ef{5.22}$ are not oscillatory
 as $y \to +\iy$ about  the SSS
 $(\ref{5.2})$
 iff  $(\ref{5.10})$    holds,
  and then:
  %%%% in particular,
  \be
  \label{5.10N}
p \ge p_{\rm H}: \quad \mbox{$W(\xi)$ has at most finite
intersections with $U(\xi)$ on $\xi \in (0,+\iy)$}.
 \ee
 \end{proposition}

 \noi{\em Proof.}  It suffices to observe that, as customary,
 the oscillatory behaviour as $y \to +\iy$ is governed by the linearized
 operator therein, which is \ef{5.5} (the limit \ef{ss1} below justifies the linearization).
 Hence, in the critical Hardy
 case, the characteristic polynomial \ef{5.14} has real roots
 only (actually, all of them, and this is  quite a general property \cite{GalH1, GalH2}),
  and obviously the same holds in the subcritical range $pD<c_{\rm
  H}$, meaning that $W(\xi)$ is not oscillatory about $U(\xi)$ as
  $\xi \to +\iy$. Clearly, if $p D > c_{\rm H}$, \ef{5.13} and \ef{5.14} imply existence
  of a proper
   root $\g \in {\mathbb C}$ with a not that large negative real part.  $\qed$

\ssk

Thus, we have concluded that, for the present problem:
 \be
 \label{jj1}
  \mbox{discrete spectrum and non-oscillation occur in the same
   range $p \ge p_{\rm H}(N)$.}
   \ee
Indeed, this has some natural roots in general spectral theory of
ordinary differential operators. For instance, for second-order
singular operators, the non-oscillating behaviour
%% (\ref{2.2}),
%%(\ref{2.3})
at singular endpoints always imply existence of a
 %% $r=0$ then implies that  spectra of
self-adjoint extension  in $L^2$ with a discrete spectrum; see
Lemma 3.1.1 in \cite[p.~74]{LS}. For higher-order symmetric
operators \cite{Nai1}, such a universal conclusion is not that
clear, though is easily observed in particular problems related to
 simpler homogeneous operators for Hardy's inequalities as in \cite{GalH1, GalH2}.

Proposition \ref{Pr.Har3} for $p > p_{\rm H}(N)$ was  proved in
\cite[p.~909]{Gaz06},
%% by a different comparison method,
where
other important properties of entire solutions $W(\xi)$
 of \ef{5.22} have been established. So we do not need to mention
 them here in detail and will use the following only (see also \cite{Dal91}):
  \be
  \label{ss1}
   \tex{
   \big(\frac{p+1}2 \big)^{\frac 1{p-1}} > \frac {W(\xi)}{U(\xi)}
   \to 1 \asA \xi \to +\iy.
   }
   \ee
  However, a number of  problems concerning \ef{5.2} remain
open. For instance, proving that (cf. Open Problem 3 in
\cite[p.~915]{Gaz06} on ordering of the family $\{W(\xi), \,
d_\b>0\}$)
 \be
 \label{gg1}
  \mbox{for $p \ge p_{\rm H}(N)$, $W(\xi)$ does not intersect
  $U(\xi)$.}
  \ee
 %%On the other hand,
%% The non-oscillatory property in Proposition
%%%\ref{Pr.Har3} settles Open Problem 2 in \cite[p.~915]{Gaz06}.
 Note that in view of inevitable using shooting techniques, the
 property \ef{gg1} is very difficult to check numerically.

%%\ssk

Fortunately, as a standard topology suggests,  solving the open
problem \ef{gg1} {\em is not necessary} for the validity of the
matching of Inner Regions I and II, since the non-oscillating of
$W(\xi)$ as $\xi \to +\iy$ is in  principal demand (one can see
that  existence of a {\em finite} number of intersections cannot
spoil matching). Thus, we conclude that:
 \be
 \label{gg2}
  \mbox{for $p \ge p_{\rm H}(N)$,
  %%% nothing prevents
  matching of two flows (\ref{5.12}) and (\ref{5.22}) is plausible,}
   \ee
 though a huge mathematical work is necessary to prove this (the author still believes that this can be
 done in a reasonably finite period of time, but its scale can be
 beyond any expectation).

%%%%%%%%%%%%%%%%%%%%%%%%%%%%%%%%%%%%%%%%%%%%%%
\subsection{On new blow-up similarity solutions in the oscillatory
range $p < p_{\rm H}$}

Thus, \ef{5.14} clearly shows that for $p<p_{\rm H}$ the solutions
$W(\xi)$ are oscillatory about the SSS $U(\xi)$. Such topology (as
in the second-order case, see \cite{GV} and later publications)
suggests that in this subcritical Hardy range there may be a
sequence of similarity profiles satisfying \ef{2.4} and exhibiting
arbitrary finite oscillations about $U(\xi)$ for sufficiently
small  radial $\xi>0$. Such self-similar blow-up profiles
concentrated in a neighbourhood of the unstable singular
equilibrium  $U$ (above $U$, a.a. solutions must blow-up), are
expected to be also highly unstable, at least in comparison with
the previous profiles $f_0$ and $f_1$ studied in Section
\ref{S.2}. Therefore, we ignore such new families (possibly
countable depending  on parameter ranges) of the s-s blow-up.

%%FIG%%%%%%%%%%%%%%%%%%%%%%%

%%%%%%%%%%%%%%%%%%%%%%%%%%%%%%%%%%%%%%%%%%%%%%%%%%%%%%%%%%%%%%%%%%%

\subsection{On related non-radial blow-up patterns}

These can be predicted in a various ways. Firstly, one can start
with a non-radial SSS solving the elliptic equation \ef{5.1}, but
surely such ones are unknown. Secondly, under the condition
\ef{5.6}, a non-radial eigenfunction $\psi_\b^*(y)$ (e.g.,
corresponding to an ``angular" swirl obtained by angular
separation of variables) of $\hat \BB^*$ can be taken into
account. Then matching will assume using non-radial entire
solutions of \ef{5.12}, which then deserves further study.

 %%%%%%%%%%%%%%%%%%%%%%%%%%%%%%%%%%%%%%%%%%%%%%%%%%%%%
\section{{\bf Type
I(LN)}: non self-similar blow-up evolution on a manifold of
generalized Loewner--Nirenberg stationary solutions}
 \label{S.6}

\subsection{Classic Loewner--Nirenberg (L--N) conformally invariant exact
solutions}

These are classic solutions obtained in Loewner--Nirenberg
\cite{LN74} in 1974  for the second-order elliptic equation
 \be
 \label{6.1}
  \tex{
 \D W + W^p=0 \inB \ren, \,\,\, W(0)=d>0, \forA p=p_{\rm
 S}=\frac{N+2}{N-2} \quad (N>2),
 }
  \ee
  which are invariant under conformal and projective
  transformations (symmetries of \ef{6.1} were earlier studied by
  Ibragimov in 1968 \cite{Ib68}). These solutions are given by
 \be
 \label{6.2}
  \tex{
  W_0(\xi) = d\, \Big[ \frac{N(N-2)}{N(N-2)+ d^{4/(N-2)}|\xi|^2}
  \Big]^{\frac{N-2}2}>0 \inB \ren
  }
   \ee
 and exhibit a number of uniqueness and other exceptional
 properties of the equation \ef{6.1}.

%%%%%%%%%%%%%%%%%%%%%%%%%%%%%%%%%%%%%%%%%%%%%
 \subsection{Generalized L--N solutions for the biharmonic equation}

For the critical biharmonic counterpart of \ef{6.1}
 \be
 \label{6.1V}
  \tex{
 -\D^2 W + |W|^{p-1}W=0 \inB \ren, \,\,\, W(0)=d>0; \,\,\, p=p_{\rm
 S}=\frac{N+4}{N-4} \quad (N>4),
 }
  \ee
the corresponding exact solutions are known from the 1980s at
least, which we call the {\em generalized L--N ones}:
%%These
%%solutions are given by
 \be
 \label{6.2V}
  \tex{
  W_0(\xi) = d \, \Big[ \frac{\sqrt{(N+2)(N^2-4)(N-4)}}{\sqrt{(N+2)(N^2-4)(N-4)}+ d^{4/(N-4)}|\xi|^2}
  \Big]^{\frac{N-4}2}>0 \inB \ren.
  }
   \ee
 The earliest references to the exact expressions \ef{6.2V} we have
 found are \cite[p.~1057]{Gal85} in 1985 and \cite{Nus92, Swan92} in 1992, where
 in the latter one important properties of $W_0$ have been proved
 (see also \cite{Guo08} for further references).
 %%\ssk
 %%\noi{\bf Remark: polyharmonic equation}.
 Note that, for the $2m$th-order polyharmonic
 extension,  the corresponding positive entire solutions look
 similarly:
  $$
   \tex{
  -(-\D)^m W + |W|^{p-1} W=0, \,\,\,p=\frac{N+2m}{N-2m}
   } %% \inB \ren
   \LongA
%%  $$
%%%   $$
   \tex{
   W_0(\xi)= d \, \Big[ \frac{B}{B+ d^{4/(N-2m)}|\xi|^2}
  \Big]^{\frac{N-2m}2},
  }
  $$
  where $N>2m$ and  $
  B^m=\frac{(N+2(m-1))!!}{(N-2(m+1))!!}$. See
 Svirshchevskii in 1993 \cite{SVR93} (in a preprint, the solutions were
  published as earlier as in 1989 \cite{SVR89}), and more related exact solutions
  of other critical elliptic PDEs (e.g., with a $p$-Laplacian)
 and extra references
  in \cite[\S~5]{GKSob}.

%%%%%%%%%%%%%%%%%%%%%%%%%%%%%%%%%%%%%%%%%%%%
\subsection{Formal construction of Type II(LN) blow-up patterns for $p=p_{\rm S}$}
 \label{SRD1}

Let $v(y,\t)$ be the rescaled solution of \ef{3.2} in, say, radial
geometry at the moment. Let us assume that $v(y,\t)$ behaves for
$\t \gg 1$ being
 close to the stationary manifold composed of the explicit  equilibria
\ef{6.2V}, i.e., for some unknown function  $\var(\t) \to +\iy$ as
$\t \to +\iy$:
%%% the following holds:
 \be
 \label{sc12}
  v(y,\t) = \var(\t) W_0\big(\var^{\frac{p-1}4}(\t)y\big)+...
 \ee
on the corresponding shrinking compact subsets in the new variable
$\zeta=\var^{\frac{p-1}4}(\t)y$. It then follows that, on the
solutions \ef{6.2V} in terms of the original rescaled variable $y$
(cf. computations in \cite[p.~2963]{Fil00}; our notations have
been slightly changed)
 \be
 \label{sc13}
  \tex{
 |v(y,\t)|^{p-1}v(y,\t) \to \frac{e_N}{\var(\t)} \,
 \d(y) \quad \mbox{as\,\, $\t \to +\iy$}
 }
  \ee
  in the sense of distributions, where $e_N>0$ is some constant.
 Therefore, on this manifold of solutions, the rescaled equation \ef{3.2}
  takes asymptotically the form
  \be
  \label{tt127}
   \tex{
 v_\t= \AAA(v) \equiv -\D^2 v- \frac 14\, y\cdot \n v - \frac
  {N-4}{8}\, v +\frac{e_N}{ \var(\t)} \,
 \d(y) +...\forA  \t \gg 1.
 }
  \ee
 According to Lemma \ref{lemSpec2}, we are looking for Type II
 patterns
 of the form
  \be
  \label{gg4}
 %%% \begin{matrix}
  \tex{
  v_\b(y,\t)= c_\b(\t) \psi_\b^*(y) +...
  %%, \quad \mbox{where} \quad
  %%\ssk\\
 \LongA  \dot c_\b=-\a_\b c_\b +
  h_\b \frac 1{\var(\t)}+...\,,
  }
 %%  \ssk\ssk\\
 %%   \mbox{with} \quad \a_\b= \frac {2|\b|+N-4}8>0
 %%  \andA h_\b= e_N \psi_\b(0).
 %% \end{matrix}
  \ee
  where  $\a_\b= \frac {2|\b|+N-4}8>0$ and
   $h_\b= e_N \psi_\b(0)$.
Simple particular ``resonance" solutions correspond to an
exponential divergence:
 \be
 \label{gg5}
  \var(\t)={\mathrm e}^{\a_\b \t}+... \andA c_\b(\t)=h_\b \t {\mathrm
  e}^{-\a_\b \t}+... \forA \t \gg 1, \quad |\b| \ge 0.
  \ee

 Bearing in mind the scaling in \ef{sc12}, this yields a countable family of distinct complicated
blow-up structures, where most of them are not radially symmetric.
To reveal the actual space-time and changing sign structures of
such Type II patterns, special matching procedures apply. In
\cite{Fil00}, this analysis has been performed in the radial
geometry for \ef{FK1p}, though
%% (and this is also key for us to
%%recognize)
 still no rigorous justification of the existence of
such blow-up scenarios is available.
%% We just mention that
Thus,  the first Fourier coefficient in \ef{gg4}
 implies a complicated structure of the pattern around the formed
Dirac's  $\d(y)$ according to \ef{sc13}. However, since these
expansions are given by generalized Hermite polynomials
$\{\psi_\b^*\}$ this matching is expected not to impose more
difficulties as those similar in Section \ref{S.4}. In any case,
%% Therefore,
more matching details for the  much harder PDE \ef{m2} seem then
excessive here.

 %%%  where $\l>0$ is arbitrary.

%%{\bf Acknowledgement. }

%%%%%%%%%%%%%%%%%%%%%%%%%%%%%%%%%%%%%%%%%%%%%
%%\newpage
% \begin{thebibliography} {444}

\begin{appendix}
%%%%%\section*{\bf Appendix} \setcounter{section}{1}
\setcounter{equation}{0}

%%%%%%%%%%%%%%%%%%%%%%%%%%%%%%%%%%%%%%%%%%%%%%%%
%%%%%%%%%%%%%%%%%%%%%%%%%%%%%%%%%%%%%%%%%%%%5

\setcounter{section}{1} \setcounter{equation}{0}
\setcounter{subsection}{0}
 \begin{small}

\section*{Appendix A:  On universality
 of the open $L^p \Longrightarrow L^\iy$ problem in PDE theory}
 \label{SState}

%%%%%%%%%%%%%%%%%%%%%%%%%%%%%%%%%%%%%%%%%%%%%%%%%
%%%%\subsection{
%%On ``non-uniqueness" of the
%% Millennuim Problem:
%%}

The Millennium Prize Problem, posed specially for the NSEs
\ef{NS1}, is, in a loose sense, ``non-unique", since similar open
regularity problems (or not that lighter significance) occur for
many evolution PDEs of various types. We list a few of them, where
the difficult open mathematical aspects of global existence and/or
blow-up are associated with the following factors:

(i) supercritical Sobolev parameter range of the principal
operator (hence, standard or very enhanced embedding-interpolation
techniques fails), and, in fact, as a corollary,

(ii) multi-dimensional space $x \in \ren$, with $N \ge 3$, at
least (this leaves a lot of room for constructing various $L^\iy$
blow-up patterns via self-similarity, angular swirl, axis
precessions, linearization, matching, etc.).

\ssk

We now list those PDEs, where we give a few recent basic
references to feel the subject.

\noi (I) Supercritical defocusing {\em nonlinear Schr\"odinger
equation}\footnote{The author would like to thank I.V.~Kamotski,
who first attracted his attention to this problem.} (NLSE) (see
\cite{Mer05}--\cite{MerR052}, \cite{Vis07, Plan07})
\be
\label{Sr1}
   \tex{
   %%%\quad
-{\rm i}\, u_t= \D u - |u|^{p-1} u \withA p>p_{\rm S}(2)=
\frac{N+2}{N-2} \quad (N \ge 3);
 }
 \ee
 (II) $2m$th-order supercritical {\em semilinear heat equation with
 absorption} ($m=1$ is covered by the MP; see \cite {GW2} and \cite{ChGal2m},
  where the result in \S~4
 for $p>p_{\rm S}(2m)$ applies to small
  solutions only):
\be
\label{Sr2}
   \tex{
 u_t= - (-\D)^m u - |u|^{p-1} u \withA p>p_{\rm S}(2m)=
 \frac{N+2m}{N-2m}\quad
 (N>2m, \,\,\,m \ge 2);
 }
 \ee
 (III) The semilinear {\em supercritical wave equations} (see
 \cite{Ikeh08, Yang08}, as most recent guides)
 %%%% to the subject)
\be
\label{Sr3}
   \tex{
 u_{tt}= \D u - |u|^{p-1} u \withA p>p_{\rm S}(2)=
 \frac{N+2}{N-2}\quad
 (N \ge 3).
 }
 \ee
 Possibly, here the Maximum Principle kind  arguments associated
 with the single Laplacian $\D$ can still play a role; then it is to be
 replaced by $-\D^2$; see below.

One can add to those ``supercritical" PDEs some others of a
different structure such as the {\em Kuramoto--Sivashinsky
equations} for $l=1,2,...$ \cite{G3}
 \be
 \label{hSr1}
  \tex{
 u_t=-(-\D)^{2l} u + (-\D)^l u+ \frac 1{p} \, \sum_{(k)} d_k
 D_{x_k}(|u|^p), \,\,\, |{\bf d}|=1,  \,\,\,  p> p_0= 1+ \frac{2(4l-1)}{N}.
 %%%%\quad (l=1,2,...).
 }
 \ee
Here, $p_0$ is not the Sobolev critical exponent, though precisely
for $p>p_0$, $L^2 \not \Rightarrow L^\iy$ by blow-up scaling,
\cite[\S~5]{G3}. On the other hand, a more exotic applied models
exhibit similar fundamental difficulties such as the following
{\em nonlinear dispersion equation} (see \cite{GalNDE5, GPnde} for
references and some details)
 \be
 \label{hS2}
  \tex{
  u_t= - D_{x_1}[(-\D)^m u]- D_{x_1} (|u|^{p-1}u)
   \withA p>p_{\rm S}(2m)=
 \frac{N+2m}{N-2m}.
 %%\quad
 %%%(N>2m, \,\,\,m \ge 2);
 }
 \ee
 %%Note that,
 In view of the  conservation properties for the
models \ef{hSr1} and \ef{hS2}, these, though being local, can be
more adequate to the nonlocal NSEs \ef{NS1}, than the others
above.

In most of the cases, the operator on the right-hand sides
satisfying for  $u \in C_0^\iy(\ren)$
 \be
 \label{Sr4}
 \tex{
 \AAA(u) = - (-\D)^m u - |u|^{p-1}u
 %% \withA
 \LongA
 \langle
 \AAA(u),u\rangle=- \int |D^m u|^2 - \int |u|^{p+1} \le 0,
 }
 \ee
is indeed coercive and monotone in the metric of $L^2(\ren)$,
which always helps for global existence-uniqueness of sufficiently
smooth solutions of these evolution PDEs. For the NLS \ef{Sr1},
this gives a stronger conservation laws than for the focusing
equation with the ``source-like" term $=|u|^{p-1}u$.  Evidently,
replacing $\D$ in \ef{Sr1} and \ef{Sr3} by $-(-\D)^m$, $m \ge 2$
moves the supercritical range to that in \ef{Sr2}. On the other
hand, introducing quasilinear differential operators $-(-\D)^m
|u|^\s u$ with $\s>0$ moves the critical exponent to $p_{\rm
S}(2m,\s)=(\s+1) \frac {N+2m}{N-2m}$. Similar supercritical PDEs
can contain $2m$th-order $p$-Laplacian operators, such as the one
for $m=2$, with $\s>0$,
\be
 \label{Sr4NN}
 \tex{
 \AAA(u) = -
\D(|\D u|^\s \D u) - |u|^{p-1}u, \quad
 %% \withA
 %%%\LongA
 \langle
 \AAA(u),u\rangle=- \int |\D u |^{\s+2} - \int |u|^{p+1} \le 0.
 }
 \ee
%%% where the absence of the compact embedding demands

%% $$ for $m=2$,  where the critical exponent
%%$p_{\rm S}(2m)(\s)$ changes according to Sobolev's embedding
%%theorem.

 However, the
lack of embedding-interpolation techniques to get $L^\iy$-bounds,
which can be expressed as the lack of
  compact Sobolev embedding of the corresponding spaces
%%(no Sobolev
%%embeddings)
  for bounded domains $\O \subset \ren$ (this analogy is not
  straightforward and is used as a certain illustration only)
 \be
 \label{Sr5}
 \tex{
 H^m(\O) \not \subset L^{p+1}(\O) \forA p
> p_{\rm S}(2m),
 }
 \ee
 actually presents the core of the problem: it is not clear how
 and when bounded solutions can attain in a finite blow-up time a
 ``singular blow-up component" in $L^\iy$.
 For the operator in \ef{Sr4NN}, a similar supercritical demand reads
\be
 \label{Sr5NN}
 \tex{
 W^2_{\s+2}(\O) \not \subset L^{p+1}(\O) \forA p
> p_{\rm S}(4,\s) = \frac{(\s+1)N+2(\s+2)}{N-2(\s+2)}, \quad
N>2(\s+2).
 }
 \ee
 In the given supercritical Sobolev
 ranges, finite mass/energy blow-up patterns for \ef{Sr1}--\ef{hS2} are unknown, as
 well as global existence of arbitrary (non-small) solutions.

It is curious that for the NSEs with the same absorption mechanism
as above,
 \be
 \label{Sr6}
  \tex{
  \uu_t +(\uu \cdot \n)\uu= - \n p + \D \uu -|\uu|^{p-1}\uu,
  \quad {\rm div}\, \uu=0 \inB \re^3 \times \re_+,
  }
  \ee
  by the same reasons and similar to \ef{Sr2},
the global existence of smooth solutions is guaranteed
\cite{Cai08} in the subcritical Sobolev range only: for
 \be
 \label{Sr7}
  \tex{
  p \le 5= \frac{N+2}{N-2}\big|_{N=3} \quad \big(\mbox{and $p \ge
  \frac 72$ by another natural reason}\big).
   }
   \ee

 We thus  claim that, even for the PDEs with local nonlinearities
 \ef{Sr1}--\ef{Sr3} (and similar higher-order others), the study
 of the admissible types of possible blow-up patterns can
 represent an important and constructive problem, with the results
 that can be key also for the non-local parabolic flows such as
 \ef{NS1},
 %%\ef{NS1m},
  \ef{Sr6}, etc. Moreover, it seems reasonable
 first to clarify the blow-up origins in some of looking similar and simpler (hopefully, yes, since
 \ef{NS1} is both nonlocal and vector-valued unlike the others)
 %%%it is not that clear)
local supercritical PDEs, and next to extend the approaches to the
non-local NSEs \ef{NS1}; though, obviously, the former ones are
not that attractive and, unfortunately, are not related to
``millennium" issues (however, many PDE experts very well
recognize how important these are for general PDE theory).

\end{small}
\end{appendix}

\end{document}